\documentclass[12pt]{elsarticle}
\usepackage{lineno,hyperref}
\usepackage[margin=1in]{geometry}
\usepackage[title]{appendix}
\usepackage{amsmath}
\usepackage{amssymb}
\usepackage{amsthm}
\usepackage{mathrsfs}
\usepackage{booktabs}
\usepackage{commath}
\usepackage{color}
\usepackage{soul}
\usepackage{multirow}
\usepackage[dvipsnames]{xcolor}
\usepackage{lineno}
\usepackage{float}
\usepackage{bm}
\restylefloat{table}

\newtheorem{example}{Example}[section]
\usepackage{caption}
\usepackage{subcaption}
\DeclareMathOperator{\sech}{sech}
\newcommand{\RomanLabelFont}{\textbf}
\newcommand{\RomanLabelSize}{\normalsize}
\newcommand{\TitleFont}{\color{blue}}
\newcommand{\TitleSize}{\large}
\newcommand{\MinMaxFont}{\relax}
\newcommand{\MinMaxSize}{\footnotesize}
\newcommand{\ValueHSpace}{0.5em}
\newcommand{\ColumnSpace}{0.5em} 
\newcommand{\mympwidth}{0.4\textwidth}
\biboptions{square,sort&compress}
\begin{document}
\begin{frontmatter}
\title{
Spectral Analysis of Node- and Cell-Centered Higher-Order Compact Schemes for Fully Discrete One and Two-Dimensional Convection–Dispersion Equation
}
\author[1]{Lavanya V Salian}
\ead{lavanya\_vs@iipe.ac.in}
\author[2,3]{Vivek S Yadav}
\ead{vyadav06@mail.ubc.ca }
\author[1]{Rathan Samala \corref{cor}}
\ead{rathans.math@iipe.ac.in}
\author[4]{Rakesh Kumar}
\ead{rakesh.kumar@mahindrauniversity.edu.in }
\address[1]{Department of Humanities and Sciences, Indian Institute of Petroleum and Energy, Visakhapatnam, 530003, Andhra Pradesh, India} \address[2]{Department of Mechanical Engineering, University of British Columbia, V6T1Z4 Vancouver, BC, Canada}
\address[3]{School of Mathematics, IISER Thiruvananthapuram, 695551, Kerala, India}
\address[4]{Department of Mathematics, École Centrale School of Engineering, Mahindra University, Hyderabad, 500043, Telangana, India}
\cortext[cor]{Corresponding author}


\begin{abstract}
\noindent 
In this study, we present a comprehensive global spectral analysis of the convection-dispersion equation, which is also referred to in specific contexts as the Korteweg-de Vries (KdV) equation, to investigate the behavior of high-order numerical schemes across a wide range of nondimensional parameters. The motivation for this analysis stems from the equation’s importance in modeling wave propagation and transport phenomena, where accurate resolution of dispersive effects is critical, and traditional numerical schemes often suffer from spurious artifacts. We analyze one sixth-order and two eighth-order compact spatial discretization schemes, encompassing both node-centered and cell-centered formulations, combined with a third-order strong stability-preserving Runge-Kutta (SSPRK3) time integrator. The analysis is performed in terms of key nondimensional parameters such as the wavenumber, Courant–Friedrichs–Lewy number (\( N_c \)), and dispersion number (\( D_{\alpha} \))-over the full spectral plane for both one- and two-dimensional cases. Key numerical indicators, including the amplification factor, normalized phase speed, and normalized group velocity, are evaluated to characterize stability, dispersion error, errors in energy transport, and directional anisotropy. Critical dispersion thresholds and Courant numbers are identified, beyond which numerical instability and nonphysical phenomena such as spurious \( q \)-waves and reversed phase or energy transport arise. Theoretical predictions are validated through numerical experiments involving linear and nonlinear one- and two-dimensional test problems, including cases with exact solutions and established benchmark results. This comprehensive analysis uncovers subtle numerical errors and offers practical guidance for selecting reliable discretization parameters, ensuring accurate and stable simulations of convection–dispersion systems.
\end{abstract}

\end{frontmatter}

\section{Introduction}
\label{sec:1}
Nonlinear partial differential equations (PDEs) with dispersive effects arise naturally in many physical contexts such as nonlinear acoustics~\cite{guo2018study, tam1993dispersion, visbal2001very}, plasma physics~\cite{rahman2021dust}, and nonlinear optics~\cite{iordache2010study, hasegawa1973transmission}. These equations often exhibit a rich variety of wave phenomena resulting from the interplay between nonlinearity and dispersion. Among such PDEs, the Korteweg-de Vries (KdV) and modified Korteweg-de Vries (mKdV) equations serve as canonical models for the study of unidirectional wave propagation in dispersive media. A general framework for such convection-dispersion equations is given by:
\begin{equation}
\begin{aligned}
u_t + \nabla \cdot \mathbf{f}(u) + \nabla \cdot L (u) &= 0, \quad (\mathbf{x}, t) \in \mathbb{R}^d \times \mathbb{R}^+, \\
u(\mathbf{x}, 0) &= u_0(\mathbf{x}), \quad \mathbf{x} \in \mathbb{R}^d,
\end{aligned}
\end{equation}
where $ u(\mathbf{x}, t) $ denotes a scalar field, $ \mathbf{f}(u) $ represents the nonlinear flux, and $ L $ is a linear dispersion operator. In one spatial dimension, setting $\mathbf{f}(u) = \frac{\mu}{2} u^2$ and $L = \partial_{xx}$ yields the classical KdV equation~\cite{korteweg1895xli}:
\begin{equation}
u_t + \mu\, u u_x + u_{xxx} = 0, \quad (x,t) \in \mathbb{R} \times \mathbb{R}^+,
\label{eq:kdv}
\end{equation}
which models the propagation of long, weakly nonlinear waves in shallow water. The term $ u u_x $ captures nonlinear steepening, while $ u_{xxx} $ introduces dispersion that leads to the formation of solitary waves. A notable feature of the KdV equation is its support for solitons, which are localized, stable waves that maintain a permanent form and speed, decay at infinity, and preserve their identity after nonlinear interactions~\cite{drazin1989solitons, wazwaz2010partial}. The mKdV equation~\cite{tanaka1972modified}, given by:
\begin{equation}
u_t + \mu\, u^2 u_x + u_{xxx} = 0, \quad (x,t) \in \mathbb{R} \times \mathbb{R}^+,
\label{eq:mkdv}
\end{equation}
replaces the quadratic nonlinearity with a cubic one and appears in diverse contexts such as modeling nonlinear wave propagation in plasmas~\cite{khater1998backlund}, traffic flow~\cite{komatsu1995kink, song2017tdgl}, fluid mechanics~\cite{helal2002soliton}, and nonlinear optics~\cite{leblond2009few}. The enhanced nonlinearity in the mKdV equation leads to more complex soliton dynamics and richer wave interactions. Accurate numerical simulation of wave-dominated problems requires discretization schemes that not only resolve spatial and temporal variations effectively but also replicate the underlying physical properties~\cite{sengupta2013high,hirsch1990numerical}. 

\par
A wide range of numerical techniques have been developed to solve convection-dispersion equations, driven by their extensive applications and intricate mathematical structure. Numerous analytical~\cite{karunakar2019differential, ablowitz1991solitons, wazwaz2007variational, fu2004new} and numerical methods~\cite{levy2004local, saucez2004method} have been proposed to study these models. Due to their structural resemblance to hyperbolic conservation laws, particularly the emergence of sharp fronts and finite-speed wave propagation, high-resolution methods are often adopted. For instance, Ahmat and Qiu~\cite{ahmat2023direct} implemented a fifth-order WENO scheme based on polynomial reconstruction, while Salian and Rathan~\cite{salian2024exponential} developed a WENO method utilizing exponential bases with an adjustable tension parameter. High-order compact finite difference schemes~\cite{li2006high, ashwin2015kdv, lele1992compact, salian2024compact} have also demonstrated strong performance. A recent scheme~\cite{salian2024novel} achieves enhanced accuracy by equating the weighted average of third derivatives at cell nodes to a weighted combination of function values at both nodes and cell centers. Notably, many prior studies rely on techniques such as von Neumann analysis and matrix-based approaches to examine the semi-discretized formulations.

\par
Global spectral analysis (GSA) provides a unified framework for evaluating numerical schemes by analyzing the fully discretized governing equations in the spectral domain. Unlike conventional approaches such as von Neumann or matrix-based analyses, GSA captures the combined influence of spatial and temporal discretizations, and extends to both non-periodic problems~\cite{suman2017spectral, sengupta2020global} and non-uniform grids~\cite{sharma2017hybrid,sengupta2016new,maurya2020new}. A key advantage of GSA is its ability to reveal how physical parameters, such as convection speed, diffusion, and dispersion coefficients, though constant in the continuous formulation, become functions of discretization-dependent length scales in numerical computations. These dependencies are often characterized within non-dimensional parameter spaces involving the wave number, CFL number, and Peclet number, where critical conditions can be identified~\cite{trefethen1982group,sengupta2013high,sagaut2023global}. GSA also enables the evaluation of numerical phase speed and group velocity, which typically differ from their physical counterparts in dispersive systems, making it particularly useful for designing accurate schemes like \textit{dispersion-relation-preserving (DRP)} methods. It has been successfully applied to various linear model equations, including convection~\cite{sengupta2007error}, diffusion~\cite{maurya2021new, sengupta2014error}, convection–diffusion~\cite{suman2017spectral, sengupta2022global, yadav2022new},  convection–diffusion–reaction systems~\cite{sengupta2020global} and reaction-diffusion equation~\cite{yadav2022spatiotemporal}. However, to the best of the authors' knowledge, a comprehensive spectral analysis that incorporates both spatial and temporal discretizations for the convection-dispersion equation in both one-dimensional (1D) and two-dimensional (2D) has not yet been addressed in the existing literature.

\par
In this work, we conduct a comprehensive global spectral analysis (GSA) of high-order finite difference schemes for the linear convection–dispersion equation in both one and two spatial dimensions. Our objective is to investigate the numerical stability and dispersion properties of these schemes using higher-order compact finite-difference methods. For time discretization, we employ the third-order strong-stability-preserving Runge-Kutta method (SSPRK3), while for spatial discretization, we consider three different compact schemes with distinct grid placements: the sixth-order node-centered compact scheme (CNCS6), the eighth-order node-centered compact scheme (CNCS8), and the eighth-order cell-centered compact scheme (CCS8). This selection allows for a direct comparison of how scheme order and the placement of degrees of freedom, whether node-centered or cell-centered, affect spectral resolution, numerical accuracy, and stability in multidimensional simulations. To assess the performance of these schemes, we compute key numerical indicators such as the amplification factor, normalized phase speed, and normalized group velocity over the spectral plane for the 1D and 2D cases. These quantities allow us to characterize the schemes' ability to resolve physical wave modes, quantify directional dispersion, and detect the emergence of numerical artifacts such as spurious \( q \)-waves, phase reversal, and artificial energy transport. Our analysis also identifies the critical values of the dispersion parameter \( D_\alpha \) and Courant number \( N_c \) for each scheme, beyond which numerical instability or nonphysical behavior arises. Notably, the spectral domain is extended to \( kh \in [0, 2\pi] \) in the case of the CCS8 due to its cell-centered structure, which supports resolution of higher-frequency modes compared to node-centered formulations. To validate the theoretical findings, we validate our spectral predictions through numerical experiments involving both linear and nonlinear partial differential equations. These include benchmark problems that exhibit wave propagation, dispersive interactions, and nonlinear steepening, such as 1D and 2D linear convection-dispersion equations, the KdV equation for single and double soliton interactions, and problems with small third-derivative coefficients leading to dispersive shocks under continuous and discontinuous initial conditions. The study also examines the mKdV equation to analyze single and two-soliton interactions, including cases without analytical solutions. The results demonstrate strong agreement between spectral predictions and numerical behavior, confirming the effectiveness of our analysis in capturing both the stability boundaries and the dispersive fidelity of the schemes in realistic computational settings. 

\par
The structure of the paper is organized as follows. Section~\ref{Sec:NM} reviews existing compact finite difference schemes and time discretization methods applicable to the convection-dispersion equation. In Section~\ref{Sec:1D_GSA}, global spectral analysis is conducted for both 1D and 2D linear cases to evaluate the effectiveness of various numerical schemes. The analysis employs both cell-node and cell-centered compact spatial discretizations, coupled with the third-order Strong Stability Preserving Runge-Kutta (SSPRK3) method for temporal integration. Theoretical characteristics are discussed, accompanied by property plots that illustrate the numerical behavior of the schemes. Section~\ref{Sec:NE} presents numerical experiments for both 1D and 2D problems, comparing the numerical solutions with exact ones under different parameter settings and validating the results through spectral analysis. Additionally, cases with no known exact solutions are investigated. Finally, Section~\ref{Sec:Conclusion} summarizes the main findings and concludes the paper.

\section{Numerical Discretization Methods}\label{Sec:NM}
In this section, we revisit the compact finite difference methods used to solve Eq.~\eqref{eq:mkdv} on a uniform mesh with a grid size of $ h_x $ in the $ x $-direction and $ h_y $ in the $ y $-direction. The mesh widths are defined as  
$h_x = \Delta x = \dfrac{x_{N_x} - x_0}{N_x},\, h_y = \Delta y = \dfrac{y_{N_y} - y_0}{N_y},$ 
for positive integers $ N_x $ and $ N_y $, respectively. For the time domain $[0, T]$, time step length is given by  $\tau = \dfrac{T}{M}$ for a positive integer $M$. The mesh points are defined as follows:  $x_i = x_0 + i h_x,\, \text{for } i = 0,1,\dots,N_x,$\, $y_j = y_0 + j h_y,\, \text{for } j = 0,1,\dots,N_y,$ $t_n = n\tau,\, \text{for } n = 0,1,\dots,M.$ The discrete mesh is then given by  $\Omega_h = \{ (x_i, y_j) \mid i = 1, \dots, N_x, \, j = 1, \dots, N_y \}.$ 

 \subsection{Spatial Discretization}
 We analyze numerical schemes for solving the two-dimensional convection-dispersion equation given by
\begin{equation}\label{eq:main_equation}
    u_t + g_1(u)_x + g_2(u)_y + f_1(u)_{xxx} + f_2(u)_{yyy} = 0.
\end{equation}
To discretize this equation in space, we employ a finite difference approach, leading to the semi-discrete formulation:
\begin{equation}
    \frac{du_{i,j}}{dt} = - \Bigl( \bigl(g_{1,x}^{\prime}(u)\bigr)_{i,j} + \bigl(g_{2,y}^{\prime}(u)\bigr)_{i,j} + \bigl(f_{1,x}^{\prime\prime\prime}(u)\bigr)_{i,j} + \bigl(f_{2,y}^{\prime\prime\prime}(u)\bigr)_{i,j} \Bigr).
\end{equation}
Here, $\bigl(g_{1,x}^{\prime}(u)\bigr)_{i,j}$ and $\bigl(g_{2,y}^{\prime}(u)\bigr)_{i,j}$ approximate $g_1(u)_x$ and $g_2(u)_y$, while $\bigl(f_{1,x}^{\prime\prime\prime}(u)\bigr)_{i,j}$ and $\bigl(f_{2,y}^{\prime\prime\prime}(u)\bigr)_{i,j}$ approximate $f_1(u)_{xxx}$ and $f_2(u)_{yyy}$ at the grid node $x_{i,j}$, respectively. We revisit compact schemes for the first- and third-order derivatives in the following subsections. Formulas for the $x$-direction are provided, and analogous formulas for the $y$-direction are given in \ref{Appendix1}.
\subsubsection{Cell-Node Compact Scheme (CNCS)}
As a first spatial discretization method, we consider the cell-node compact finite difference scheme (CNCS), following the approach in~\cite{lele1992compact,li2006high,salian2024novel}.
For approximating first-order derivatives in the $x$-direction, the compact finite difference scheme takes the form:
\begin{equation}\label{FDCNCS}
\begin{split}   
    \alpha_1 \bigl(g_{1,x}^{\prime}\bigr)_{i-1,j} + \bigl(g_{1,x}^{\prime}\bigr)_{i,j} + \alpha_1 \bigl(g_{1,x}^{\prime}\bigr)_{i+1,j} &= a_1\dfrac{ \bigl(g_{1,x}\bigr)_{i+1,j} - \bigl(g_{1,x}\bigr)_{i-1,j} }{2h_x} + b_1 \dfrac{ \bigl(g_{1,x}\bigr)_{i+2,j} - \bigl(g_{1,x}\bigr)_{i-2,j} }{4h_x} \\
    &+ c_1 \dfrac{ \bigl(g_{1,x}\bigr)_{i+3,j} - \bigl(g_{1,x}\bigr)_{i-3, j} }{6h_x}.
\end{split}    
\end{equation}
A tridiagonal scheme with sixth-order accuracy is derived using the coefficients $\alpha_1 = 1/3,\, a_1 = 14/9,\, b_1=1/9,\, c_1=0,$ while an eighth-order accurate scheme is achieved using $\alpha_1 = 3/8,\, a_1 = 25/16,\, b_1=1/5,\, c_1=-1/80.$
\par
For third-order derivatives in the $x$-direction, the compact finite difference scheme is formulated as~\cite{li2006high, salian2024novel}:
\begin{equation}\label{TDCNCS}
\begin{split}
    \alpha_1 \bigl(f_{1,x}^{\prime\prime\prime}\bigr)_{i-1,j} + \bigl(f_{1,x}^{\prime\prime\prime}\bigr)_{i,j} + \alpha_1 \bigl(f_{1,x}^{\prime\prime\prime}\bigr)_{i+1,j} &= a_1\frac{ \bigl(f_{1,x}\bigr)_{i+2,j} - 2\bigl(f_{1,x}\bigr)_{i+1,j} + 2\bigl(f_{1,x}\bigr)_{i-1,j} - \bigl(f_{1,x}\bigr)_{i-2,j} }{2h_x^3}  \\
    &+ b_1 \frac{ \bigl(f_{1,x}\bigr)_{i+3,j} - 3\bigl(f_{1,x}\bigr)_{i+1,j} + 3\bigl(f_{1,x}\bigr)_{i-1,j} - \bigl(f_{1,x}\bigr)_{i-3,j} }{8h_x^3}\\
    &+ c_1 \frac{ \bigl(f_{1,x}\bigr)_{i+4,j} - 4\bigl(f_{1,x}\bigr)_{i+1,j} + 4\bigl(f_{1,x}\bigr)_{i-1,j} - \bigl(f_{1,x}\bigr)_{i-4,j} }{20h_x^3}.
\end{split}
\end{equation}
A sixth-order tridiagonal scheme is obtained with $\alpha_1=7/16,\, a_1 = 2,\, b_1=-1/8,\, c_1=0.$ An eighth-order tridiagonal scheme is achieved using $\alpha_1 = 205/472,\, a_1=2367/1180,\, b_1=-167/1180,\, c_1=1/236.$

\subsubsection{Central Compact Scheme (CCS)}
As an alternative spatial discretization approach, we employ the central compact finite difference scheme (CCS), as described in~\cite{liu2013new, salian2024novel}.
\par
For first-order derivatives in the $x$-direction, the compact finite difference scheme~\cite{liu2013new} is given by:
\begin{equation}\label{FDCCS}
    \begin{split}
       \alpha_1 \bigl(g_{1,x}^{\prime}\bigr)_{i-1,j} + \bigl(g_{1,x}^{\prime}\bigr)_{i,j} + \alpha_1 \bigl(g_{1,x}^{\prime}\bigr)_{i+1,j} &= a_1\dfrac{\bigl(g_{1,x}\bigr)_{i+\frac{1}{2},j} - \bigl(g_{1,x}\bigr)_{i-\frac{1}{2},j}}{h_x}+ b_1\dfrac{ \bigl(g_{1,x}\bigr)_{i+1,j} - \bigl(g_{1,x}\bigr)_{i-1,j} }{2h_x}\\
       &+ c_1 \dfrac{\bigl(g_{1,x}\bigr)_{i+\frac{3}{2},j} - \bigl(g_{1,x}\bigr)_{i-\frac{3}{2},j} }{3h_x}.
    \end{split}
\end{equation}
A sixth-order accurate tridiagonal scheme is obtained with coefficients $\alpha_1 = -1/12$, $a_1 = 16/9$, $b_1=-17/18$, $c_1=0$. The eighth-order accurate scheme is derived using $\alpha_1 = -3/20$, $a_1 = 2$, $b_1 =-61/50$, $c_1 =-2/25$.
\par
For third-order derivatives in the $x$-direction, the corresponding compact finite difference scheme is formulated as~\cite{salian2024novel}:
\begin{equation}\label{TDCCS}
    \begin{split}
       \alpha_1 \bigl(f_{1,x}^{\prime\prime\prime}\bigr)_{i-1,j} + \bigl(f_{1,x}^{\prime\prime\prime}\bigr)_{i,j} + \alpha_1 \bigl(f_{1,x}^{\prime\prime\prime}\bigr)_{i+1,j}
       &= a_1\dfrac{ 4\bigl(f_{1,x}\bigr)_{i+1,j} -8\bigl(f_{1,x}\bigr)_{i+\frac{1}{2},j} +8\bigl(f_{1,x}\bigr)_{i-\frac{1}{2},j} -4 \bigl(f_{1,x}\bigr)_{i-1, j} }{h_x^3} \\
       &+ b_1\dfrac{ 8\bigl(f_{1,x}\bigr)_{i+\frac{3}{2},j} -12\bigl(f_{1,x}\bigr)_{i+1, j} + 12\bigl(f_{1,x}\bigr)_{i-1,j}-8\bigl(f_{1,x}\bigr)_{i-\frac{3}{2},j}}{5h_x^3}\\
       &+ c_1\dfrac{ 8\bigl(f_{1,x}\bigr)_{i+\frac{5}{2}, j} -20\bigl(f_{1,x}\bigr)_{i+1, j} + 20\bigl(f_{1,x}\bigr)_{i-1, j}-8\bigl(f_{1,x}\bigr)_{i-\frac{5}{2}, j}}{35h_x^3}.
    \end{split}
\end{equation}
The tridiagonal compact scheme achieves sixth-order accuracy when the coefficients are chosen as $\alpha_1 = -\frac{1}{2} $, $a_1 = 5 $, $b_1 = -5 $, and $c_1 = 0 $. Since $\alpha_1 = -\frac{1}{2}$ results in the loss of strict diagonal dominance, which is required for numerical stability (i.e., the condition $|\alpha_1| < 0.5$), the scheme may become unstable, particularly for high wavenumbers. Therefore, an eighth-order tridiagonal scheme is instead adopted with $\alpha_1 = -1261/3530$, $a_1 = 58021/14120$, $b_1 = -109007/28240$, and $c_1 = 1029/28240$. Both the first- and third-order derivative formulas involve values at the cell centers, which are initially unknown. These values are treated as independent computational variables. To update both the cell nodes and the cell centers, the same scheme is applied by shifting the indices in Eqs.~(\ref{FDCCS}) and (\ref{TDCCS}) by $1/2$.
 \begin{figure}[htbp!]  
    \centering
    \begin{minipage}[b]{0.45\linewidth}
    \includegraphics[width=\linewidth, trim={0 0 0 0}, clip]{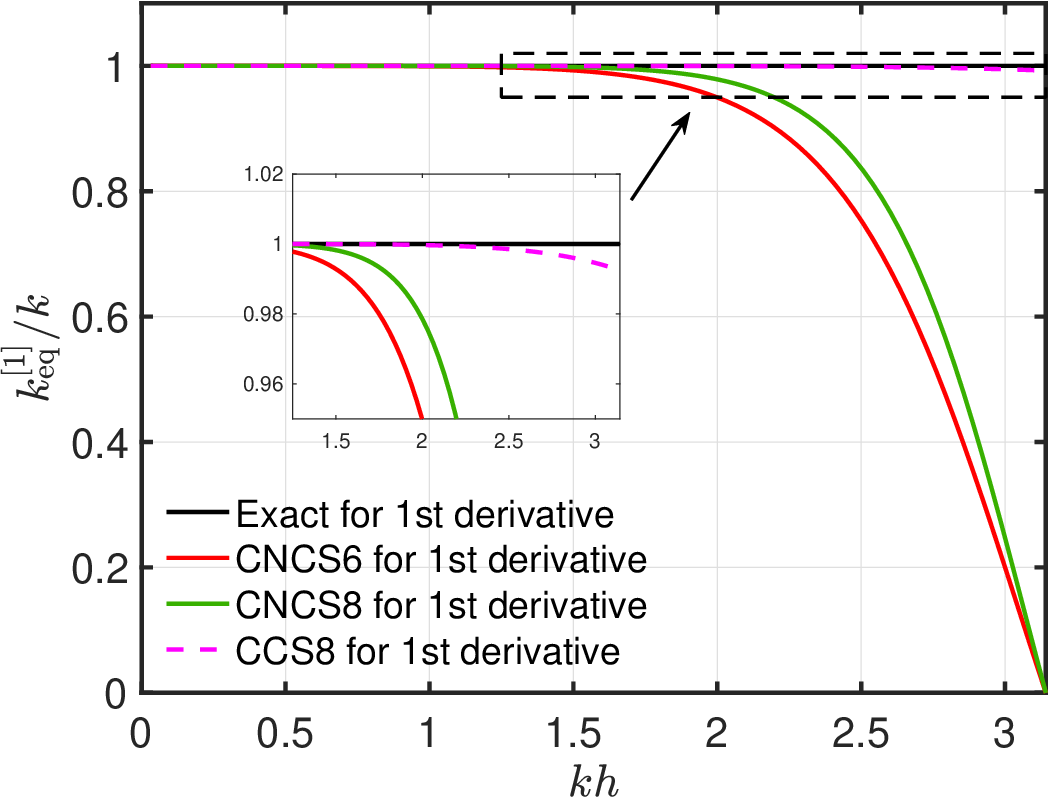}
      \subcaption{}
      \label{Fig:Spectrala}
    \end{minipage}\hfill
    \begin{minipage}[b]{0.45\linewidth}
    \includegraphics[width=\linewidth, trim={0 0 0 0}, clip]{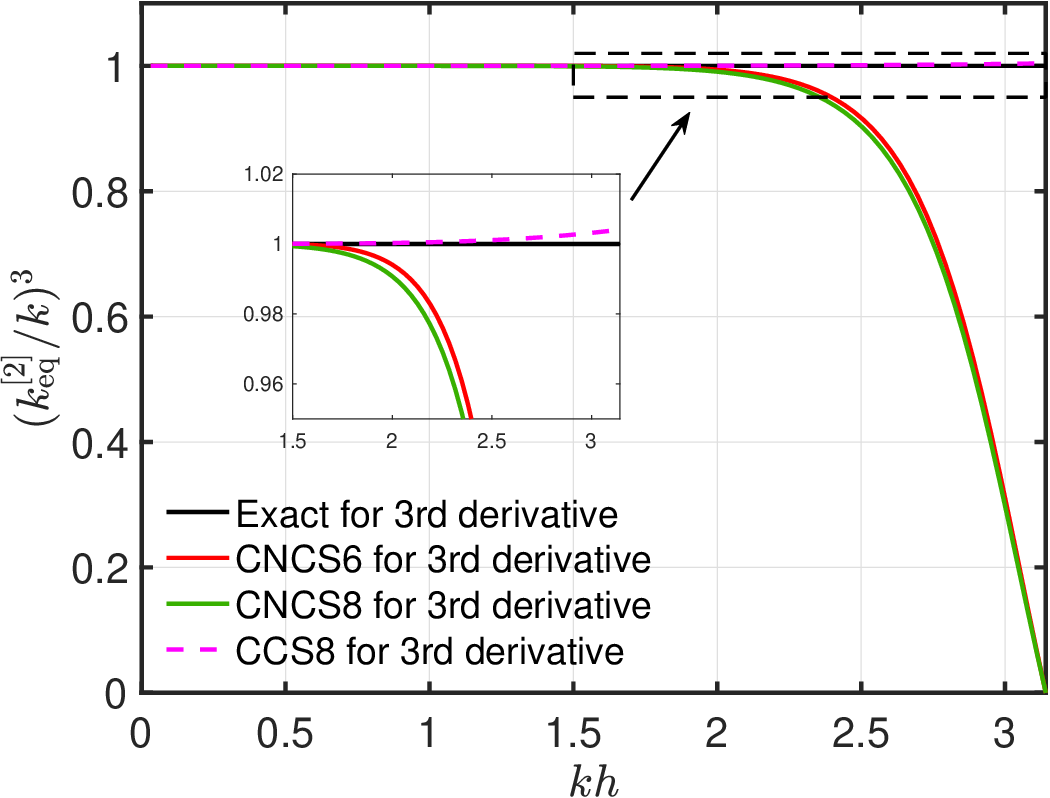}
      \subcaption{}
      \label{Fig:Spectralb}
    \end{minipage}\hfill
    \caption{Spectral resolution characteristics of indicated first- and third-order derivative schemes, plotted as a function of the non-dimensional wavenumber $ kh $.}
    \label{Fig:Spectral}
\end{figure}
\par
Fig.~\ref{Fig:Spectral} shows the comparison of the spectral resolution for the first- and third-order derivatives in one dimension using three different compact finite difference schemes: CNCS6, CNCS8, and CCS8. The expressions for the equivalent wavenumbers corresponding to the methods considered in this study are provided in~\ref{Appendix2}. In the case of the first derivative given in Fig.~\ref{Fig:Spectrala}, the CCS8 demonstrates the closest agreement with the ideal spectral resolution~\cite{rajpoot2012solution,liu2013new}, maintaining $k^{[1]}_{\mathrm{eq}}/k \approx 1$ up to a non-dimensional wavenumber of $ kh \approx 2.3 $,  which is also highlighted inside the inlarged inset figure. This high spectral resolution ensures minimal numerical dispersion and improved transport behavior. The CNCS8 shows an improved resolution over CNCS6 due to its higher order accuracy, maintaining good agreement up to $ kh \approx 1.4 $, but it still exhibits a visible departure from unity at higher wavenumbers. CNCS6, being the lowest order among the three, deviates the earliest, beginning around $ kh \approx 1.3 $, highlighting its limited capacity to resolve finer spatial structures. A similar behaviour is observed in the third derivative spectral resolution given in Fig.~\ref{Fig:Spectralb}, where CCS8 again outperforms both CNCS schemes. The CCS8 maintains a flatter and more accurate resolution curve, i.e., $(k^{[2]}_{\mathrm{eq}}/k)^3 \approx 1$, throughout the Nyquist range~\cite{sengupta2013high,salian2024novel}, but shows over-resolution around $kh \approx 2.3$ as it is enlarged inside the inset figure. Although CNCS8 exhibits better resolution than CNCS6, both under-resolve the higher wavenumber components. These deficiencies can result in non-physical numerical artifacts such as \( q \)-waves, especially in long-time or high-frequency simulations. 


 \subsubsection{Temporal Discretization}
 In this study, the governing equations are integrated in time using third-order explicit strong stability preserving Runge–Kutta scheme (SSPRK3)~\cite{gottlieb1998total, gottlieb2011strong}. Consider a general evolution equation of the form:
\begin{equation}
    \dfrac{\partial u}{\partial t} = \mathcal{S}(u),
\end{equation}
where $\mathcal{S}$ represents an operator involving spatial partial derivatives. Upon spatial discretization, the governing equation transforms into
\begin{equation}
    \dfrac{\partial u_{i,j}}{\partial t} = \mathcal{S}_h(u_{i,j}),
\end{equation}
where $\mathcal{S}_h$ denotes the spatial discretization operator. The time integration from $t_0$ (corresponding to the $n^{\text{th}}$ time step) to $t_0+\Delta t$ (i.e., the $(n+1)^{\text{th}}$ step) is performed using the following three-stage third-order SSPRK3 scheme:
\begin{align}
    u^{(0)} &= u^n, \notag \\
    u^{(1)} &= u^{(0)} + \Delta t\, \mathcal{S}_h(u^{(0)}), \notag \\
    u^{(2)} &= \dfrac{3}{4}u^{(0)} + \dfrac{1}{4}u^{(1)} + \dfrac{1}{4} \Delta t\, \mathcal{S}_h(u^{(1)}), \label{Time} \\
    u^{n+1} &= \dfrac{1}{3}u^{(0)} + \dfrac{2}{3}u^{(2)} + \dfrac{2}{3} \Delta t\, \mathcal{S}_h(u^{(2)}). \notag
\end{align}
Utilizing Eq.~\eqref{Time} along with the spatial discretization schemes, the numerical solution at each subsequent time step can be computed.

\section{Spectral Analysis of Fully Discrete Schemes}
\label{Sec:1D_GSA}
To assess the numerical behavior of the fully discrete convection-dispersion equation, we perform a global spectral analysis in both one and two dimensions. This analysis characterizes the dispersion and dissipation properties introduced by the numerical schemes, providing insight into their stability, accuracy, and suitability for capturing sharp gradients and convection-dominated features. The main objective is to understand how the combination of compact finite difference schemes in space and SSPRK3 in time affects the spectral characteristics of the numerical solution.
\subsection{Spectral analysis of fully discrete linear 1-D convection dispersion equation}
\label{Sec:GSA}
The theoretical and numerical analysis of the linear 1D convection-dispersion equation is presented first, followed by the application of GSA to the SSPRK3-CNCS and SSPRK3-CCS to demonstrate its contribution in the context of numerical method analysis. The primary focus regarding the scheme's response is on stability, dispersion, and group velocity characteristics. For simplicity, we denote $ k_x $ by $ k $ and $ h_x $ by $ h $ in the $ x $-direction. 

\par
Consider the one-dimensional linear convection-dispersion equation with periodic boundary conditions:
\begin{align}\label{eqn:Main}
    u_t + c\, u_x + \alpha\, u_{xxx} &= 0, \quad (x,t) \in \mathbb{R} \times \mathbb{R}^+,\notag\\
    u(x,0) &= u_0(x),
\end{align}
where the constants $c $ and $\alpha $ represent the convection speed and dispersion coefficient, respectively. To perform a global spectral analysis, the function $u(x,t) $ is represented using a hybrid spectral-plane formulation, where time remains in the physical domain while space is expressed in the wavenumber $k$ domain~\cite{sengupta2003analysis}. The function $u(x,t) $ can be expressed as:
\begin{equation}
    u(x,t) = \int \hat{U}(k,t) e^{\iota kx} \, dk,
\end{equation}
where $\hat{U}(k,t) $ is the Fourier amplitude, and the integral extends over the wavenumber range  $-k_{\max}$ to $k_{\max}$ defined by the Nyquist limit, given by $k_{\max} = \dfrac{\pi}{h} $. In the continuous limit as $h \to 0 $, the integration limits extend to infinity. Substituting this representation into Eq.~\eqref{eqn:Main}, we obtain its spectral form:
\begin{equation}
    \dfrac{\partial \hat{U}}{\partial t} + \iota c k \hat{U} - \iota \alpha k^3 \hat{U} = 0.
\end{equation}
For an initial condition $u(x,0) = f(x) = \int U_0(k) e^{\iota kx} dk$, the exact solution is given by:
\begin{equation}\label{exact_sol}
    \hat{U}(k,t) = U_0(k) e^{\iota(\alpha k^3 - c k)t}.
\end{equation}
To derive the physical dispersion relation, the unknown $u(x,t)$ is expressed using the bi-dimensional Fourier-Laplace transform as
\begin{equation}
    u(x,t) = \int \int \bar{U}(k, \omega) e^{i(kx - \omega t)} \, dk \, d\omega,
\end{equation}
where $\bar{U}(k,\omega)$ denotes the Fourier-Laplace amplitude, and $\omega$ is the angular frequency associated with the wave.
Substituting this into Eq.~(\ref{eqn:Main}) leads to the dispersion relation:
\begin{equation}\label{phy_dis_relation}
    \omega = c k - \alpha k^3.
\end{equation} 
This dispersion relation characterizes the wave propagation and dispersion within the medium, showing that the frequency $\omega $ depends on the wavenumber $k$, convection speed $c $, and dispersion coefficient $\alpha $. From this, we define the phase speed as:
\begin{equation}
    c_{\text{ph}} = \dfrac{\omega}{k} = c - \alpha k^2.
\end{equation}
The group velocity, representing the speed at which energy propagates, is given by:
\begin{equation}
    v_{\text{g}} = \dfrac{\partial \omega}{\partial k} = c - 3\alpha k^2.
\end{equation}
The physical amplification factor, denoted as $G_{\text{ph}}$, is obtained from Eq.~\eqref{exact_sol} as follows:
\begin{equation}
    G_{\text{ph}} = \dfrac{\hat{U}(k, t+\Delta t)}{\hat{U}(k, t)} = e^{\iota(\alpha k^3 - c k)\Delta t} = e^{-\iota \omega \Delta t} = e^{-\iota N_c kh} e^{\iota D_{\alpha}(kh)^3},
\end{equation}
where the parameters are defined as $N_c = \dfrac{c\Delta t}{h}$ (CFL number) and $D_{\alpha} = \dfrac{\alpha \Delta t}{h^3}$ (Dispersion number). Since the exponent is purely imaginary, the magnitude of $G_{\text{ph}}$ remains unity, i.e., $|G_{\text{ph}}| = 1$.

\par
By applying the discrete Fourier transform, the unknown function $u$ along with its spatial derivatives can be expressed in terms of their Fourier components~\cite{vichnevetsky1982fourier}.
\begin{equation}
    u(x, t)\Big|_{\text{num}} = \int_{-k_{\max}}^{k_{\max}} \hat{U}(k,t) e^{\iota kx} dk,
\end{equation}
\begin{equation}\label{eqn:numerical1}
     u^{\prime} (x, t)\Big|_{\text{num}} = \int_{-k_{\max}}^{k_{\max}} (\iota k^{[1]}_{\text{eq}}) \hat{U}(k,t) e^{\iota kx} dk,
\end{equation}
\begin{equation}\label{eqn:numerical2}
    u^{\prime\prime\prime} (x, t)\Big|_{\text{num}} = \int_{-k_{\max}}^{k_{\max}} (-\iota (k^{[2]}_{\text{eq}})^3) \hat{U}(k,t) e^{\iota kx} dk.
\end{equation}
Here \(k^{[1]}_{\text{eq}}\) and \(k^{[2]}_{\text{eq}}\) represent the equivalent wave number for first and third-order derivatives, respectively. In discrete computations, the discrepancy between the equivalent wavenumber $ k_{\text{eq}} $ and the actual wavenumber $k$ acts as a significant source of error. Many researchers have aimed to improve the accuracy of numerical schemes by minimizing the departure of $ k_{\text{eq}} $ from $ k $, often through the use of appropriate error norms, as demonstrated in~\cite{ashwin2015kdv,li2006high,salian2024novel}. In general, the terms $\iota k^{[1]}_{\text{eq}}$ and $-\iota (k^{[2]}_{\text{eq}})^{3}$ represent complex values, where the real component is associated with diffusion, while the imaginary part corresponds to undesirable numerical dispersion. To explore further numerical characteristics, we substitute Eqs.~(\ref{eqn:numerical1}) and (\ref{eqn:numerical2}) into the governing equation, Eq.~\eqref{eqn:Main}, resulting in:
\begin{equation}\label{eqn:IntLKdV}
    \int_{-k_{\max}}^{k_{\max}} \Biggl[\dfrac{\partial \hat{U}(k,t)}{\partial t} + c (\iota k^{[1]}_{\text{eq}}) \hat{U}(k,t) +\alpha (-\iota (k^{[2]}_{\text{eq}})^3) \hat{U}(k,t) \Biggr]e^{\iota kx} dk = 0.
\end{equation}
Since Eq.~\eqref{eqn:IntLKdV} holds for all wave numbers, the integrand must vanish for each $k$. Reformulating this condition yields:
\begin{equation}\label{eqn:DiffLKdV1}
    \begin{split}
    \dfrac{d\hat{U}}{\hat{U}} &= -c\,dt (\iota k^{[1]}_{\text{eq}}) - \alpha \,dt (-\iota (k^{[2]}_{\text{eq}})^3), \\
    &= -N_c (\iota k^{[1]}_{\text{eq}}h) - D_{\alpha} (-\iota (k^{[2]}_{\text{eq}})^3 h^3). \\
    \end{split}
\end{equation}
The left-hand side of the above equation can be expressed using the numerical amplification factor, $G_{\text{num}} = \hat{U}(k, t+\Delta t)/\hat{U}(k, t)$, which characterizes the growth or decay of the solution over a time step. For the three-stage third-order SSPRK3 time integration scheme, the numerical amplification factor is given by
\begin{equation}   
G_{\text{num}} = 1 - A + \frac{A^2}{2} - \frac{A^3}{6},
\end{equation}
where $ A = N_c (\iota k^{[1]}_{\text{eq}} h) + D_{\alpha} (-\iota (k^{[2]}_{\text{eq}})^3 h^3) $. 
The expressions for the equivalent wavenumbers for the different methods considered in this paper are provided in \ref{Appendix2}.
The numerical dispersion relation, analogous to the analytical relation Eq.~\eqref{phy_dis_relation}, is given by~\cite{sengupta2013high,sengupta2007error,sengupta2020global}:
\begin{equation}\label{num_dis_relation}
    \omega_{\text{num}} = c_{\text{num}}k - \alpha_{\text{num}} k^3,
\end{equation}
where $\omega_{\text{num}}$ denotes the numerical circular frequency, which deviates from the exact expression because both the effective numerical speed $c_{\text{num}}$ and the dispersion coefficient $\alpha_{\text{num}}$ are functions of the wavenumber $k$ and also depend on the numerical parameters $N_c$ and $D_{\alpha}$. The numerical phase shift per time step is given by~\cite{sengupta2007error, sengupta2013high}:
\begin{equation}
    \tan(\beta_{\text{num}}) = -\Biggl( \dfrac{(G_{\text{num}})_{Imag}}{(G_{\text{num}})_{Real}}\Biggr),
\end{equation}
from which the phase shift can be computed as
\begin{equation}
    \beta_{\text{num}} = \left(c_{\text{num}}k - \alpha_{\text{num}}k^3\right)\Delta t.
\end{equation}
The dimensionless numerical phase speed is expressed as:
\begin{equation}\label{phase_speed}
    \dfrac{c_{\text{num}}}{c_{\text{ph}}}(N_c, D_{\alpha}, kh) = \dfrac{\beta_{\text{num}}}{(kh) N_c - (kh)^3 D_{\alpha}}. 
\end{equation}
The numerical group velocity, $v_{\text{g,num}}$, can be computed from the numerical dispersion relation as $v_{\text{g,num}} = \dfrac{\partial \omega_{\text{num}}}{\partial k}$, which simplifies to:
\begin{equation}\label{group_vel}
    \dfrac{v_{\text{g,num}}}{v_{\text{g}}}(N_c, D_{\alpha}, kh)  = \dfrac{1}{N_{c} - 3 D_{\alpha}(kh)^{2}} \dfrac{d \beta_{\text{num}}}{dkh}.
\end{equation}

\par 
Accurate numerical simulation of convection-dispersion problems requires discretization schemes that resolve spatial and temporal features effectively and preserve the underlying physical characteristics~\cite{sengupta2013high,hirsch1990numerical}. For problems involving wave propagation and flow dynamics, the numerical method must replicate the physical dispersion behavior. For the one-dimensional convection-dispersion equation described in Eq.~\eqref{eqn:Main}, the numerical dispersion relation (Eq.~\eqref{num_dis_relation}) should closely match the exact physical dispersion relation (Eq.~\eqref{phy_dis_relation}) within the relevant parameter space. Schemes that satisfy this condition are termed \textit{dispersion-relation-preserving (DRP)} methods~\cite{sengupta2013high}. The derived expressions in Eqs.~(\ref{phase_speed}) and (\ref{group_vel}) reveal a strong dependence of numerical characteristics on key non-dimensional parameters such as the wave number ($kh$), Courant number ($N_c$), and dispersion number ($D_{\alpha}$), which influence the scheme's stability, dispersion, and transport behavior. Accuracy is typically assessed by how closely the numerical phase speed and group velocity ratios, $\dfrac{c_{\text{num}}}{c_{\text{ph}}}$ and $\dfrac{v_{\text{g,num}}}{v_{\text{g}}}$, approach unity. In the literature, the GSA has been successfully applied to PDEs with non-periodic problems, where boundary conditions are incorporated into the analysis~\cite{suman2017spectral, sengupta2020global}, as well as to non-uniform grids~\cite{sharma2017hybrid,sengupta2016new}. In this study, we employ GSA on a periodic and uniformly spaced mesh consisting of 2000 grid points to evaluate the performance of the SSPRK3-CNCS6, SSPRK3-CNCS8, and SSPRK3-CCS8 schemes. To assess their fully discrete spectral behavior, next, we have plotted \( |G_{\text{num}}| \), \( \dfrac{c_{\text{num}}}{c_{\text{ph}}} \), and \( \dfrac{v_{\text{g,num}}}{v_{\text{g}}} \) for each of the aforementioned schemes.

\subsubsection{SSPRK3-CNCS6} 
\begin{figure}[htbp!]
    \centering
    \begin{minipage}{\textwidth}
    \hspace*{0cm} 
    \begin{tabular}{ccc}
        \multicolumn{3}{c}{\makebox[0pt]{\hspace{0em}\textcolor{red}{\textbf{SSPRK3-CNCS6}}}} \\
        \multicolumn{1}{c}{\textbf{$\bm{D}_{\bm{\alpha}} = \bm{0.11}$}} & \textbf{$\bm{D}_{\bm{\alpha,\mathrm{cr}}} = \bm{0.12}$} \\[0.5em]

        \begin{minipage}{\mympwidth}
            \setlength{\fboxrule}{0.4pt} 
            \fcolorbox{black}{white}{
               \begin{minipage}{\textwidth}
                    \makebox[\linewidth][l]{\RomanLabelSize\RomanLabelFont (i)} \\[-0.6em]
                    {\centering \TitleSize\TitleFont {$|G_{\text{num}}|$} \par}
                    {\MinMaxSize\MinMaxFont \makebox[\linewidth]{\hspace{\ValueHSpace}Min = 0.94\hfill Max = 21.09\hspace{\ValueHSpace}}} \\[0.2em]
                    \includegraphics[width=\linewidth, trim=20 2 20 20, clip]{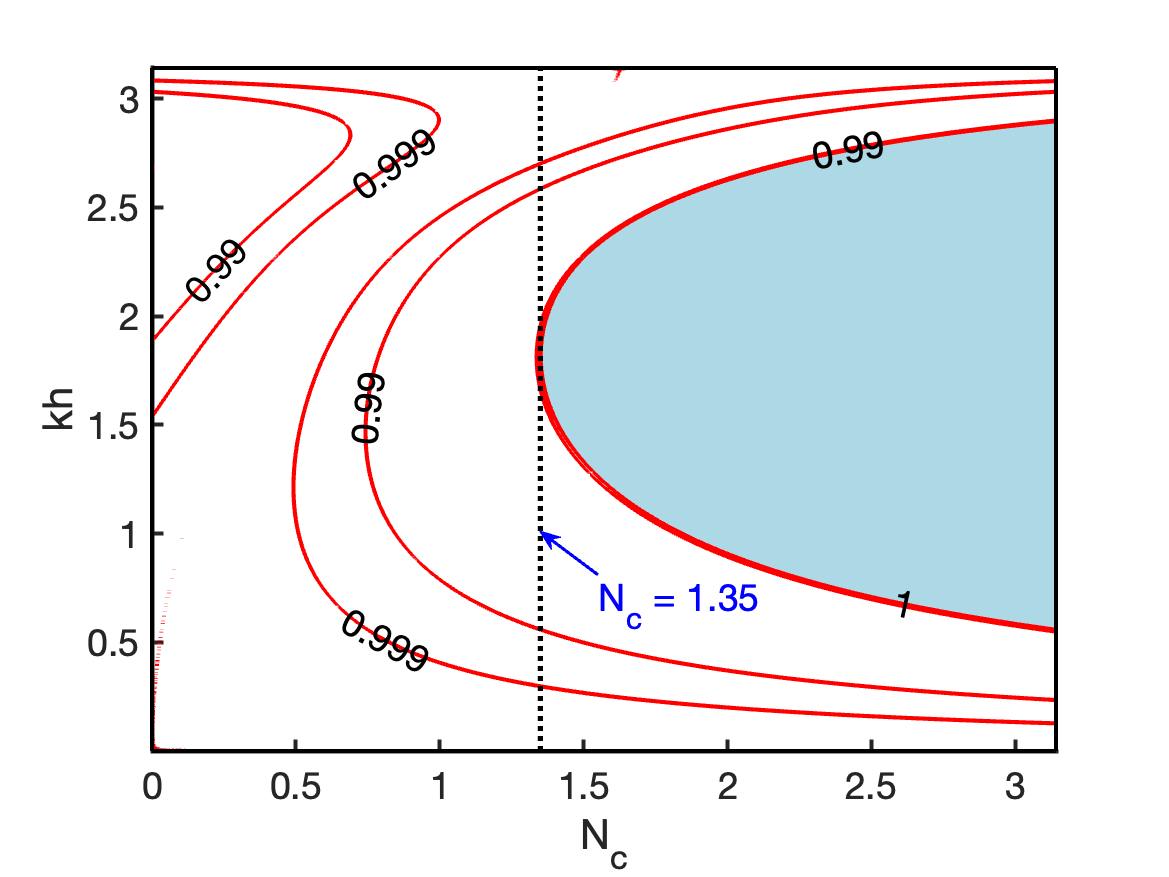}
               \end{minipage}
            }
        \end{minipage}
        &
        \hspace{\ColumnSpace} 
        \begin{minipage}{\mympwidth}
            \fbox{
                \begin{minipage}{\textwidth}
                    \makebox[\linewidth][l]{\RomanLabelSize\RomanLabelFont (iv)} \\[-0.6em]
                    {\centering \TitleSize\TitleFont {$|G_{\text{num}}|$} \par}
                    {\MinMaxSize\MinMaxFont \makebox[\linewidth]{\hspace{\ValueHSpace}Min = 0.94\hfill Max = 19.96\hspace{\ValueHSpace}}} \\[0.2em]
                    \includegraphics[width=\linewidth, trim=20 2 20 20, clip]{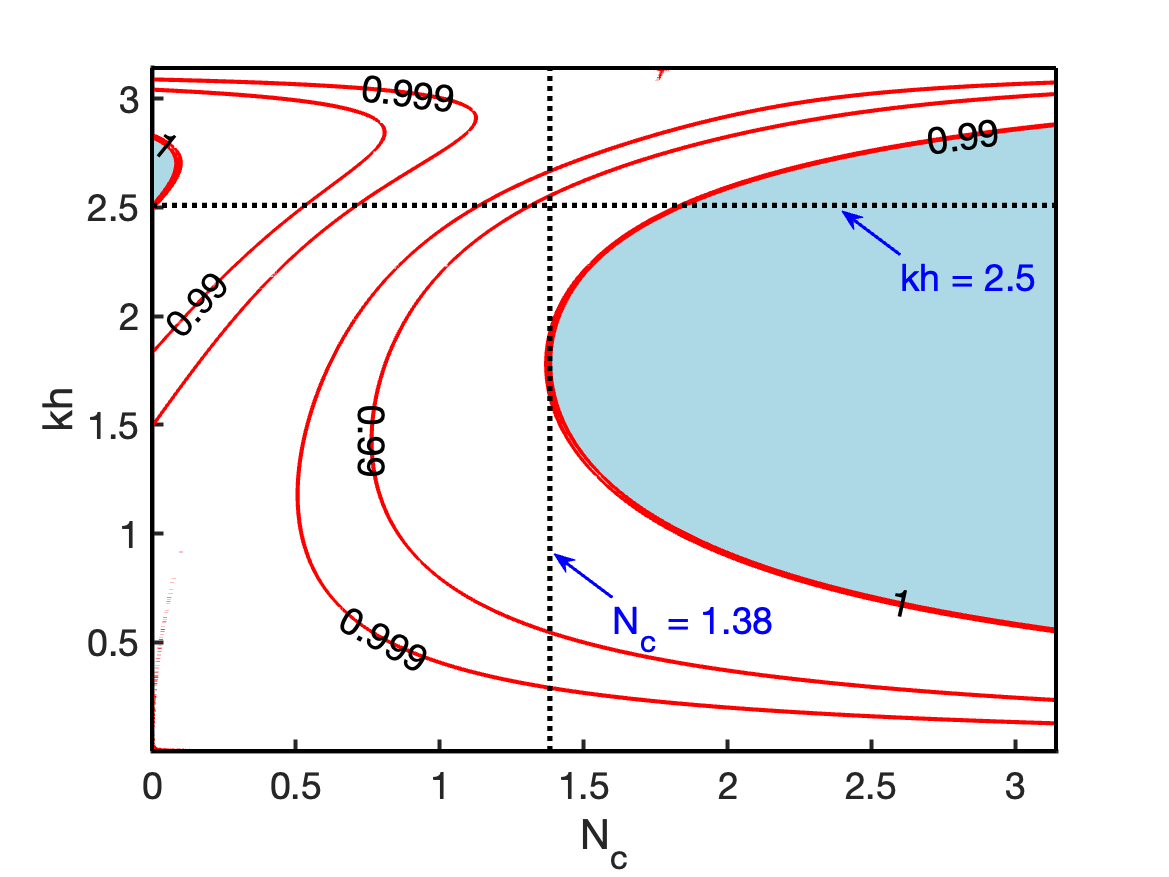}
                \end{minipage}
            }
        \end{minipage}
        &
        \hspace{\ColumnSpace} 
        \begin{minipage}{0.08\textwidth}
            \includegraphics[width=\textwidth]{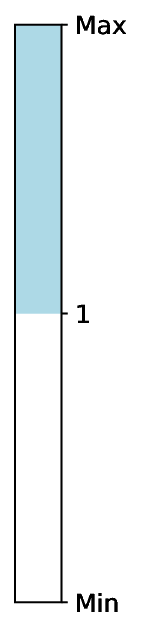}
        \end{minipage}
        \\[1em]

        \begin{minipage}{\mympwidth}
            \fbox{
                \begin{minipage}{\textwidth}
                    \makebox[\linewidth][l]{\RomanLabelSize\RomanLabelFont (ii)} \\[-0.6em]
                    {\centering \TitleSize\TitleFont {$c_{\text{num}}/c_{\text{ph}}$} \par}
                    {\MinMaxSize\MinMaxFont \makebox[\linewidth]{\hspace{\ValueHSpace}Min = -9.66e+04\hfill Max = 6.57e+04\hspace{\ValueHSpace}}} \\[0.2em]
                    \includegraphics[width=\linewidth, trim=20 2 20 20, clip]{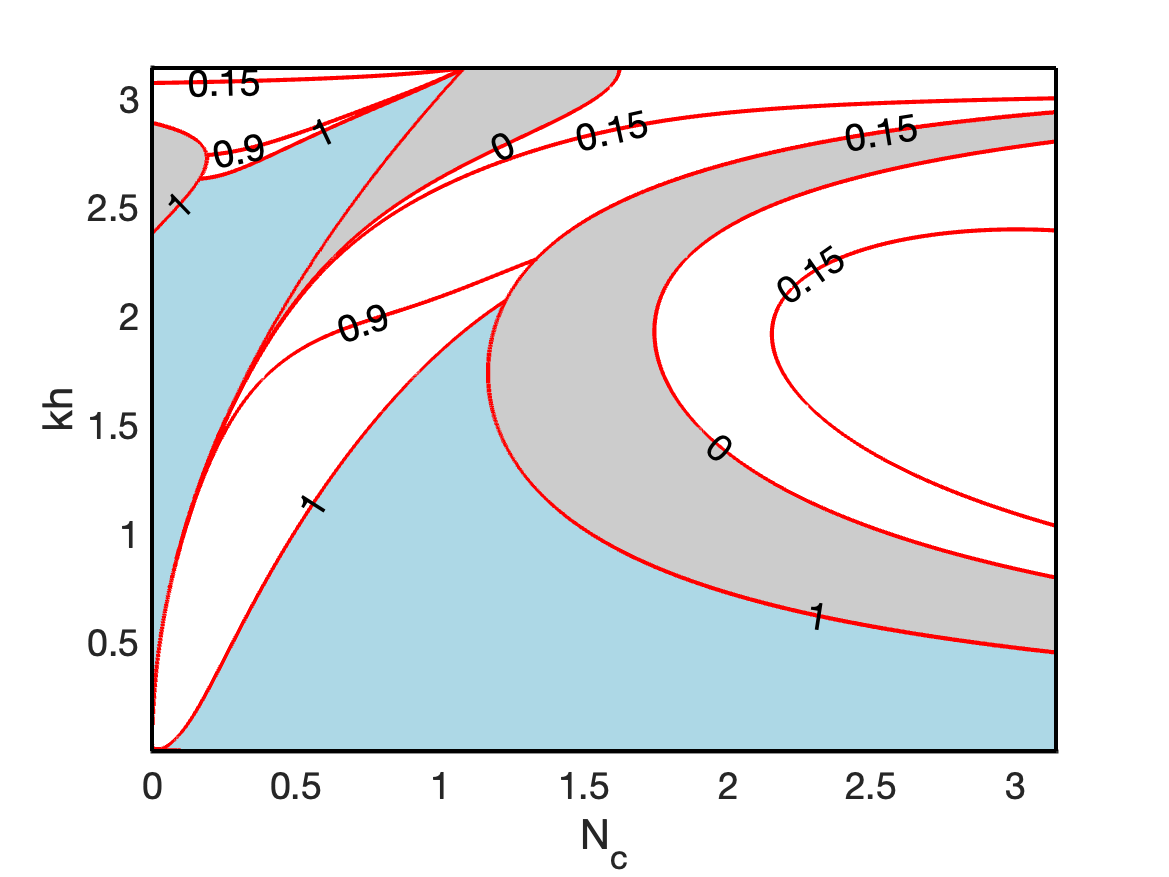}
                \end{minipage}
            }
        \end{minipage}
        &
        \hspace{\ColumnSpace} 
        \begin{minipage}{\mympwidth}
            \fbox{
                \begin{minipage}{\textwidth}
                    \makebox[\linewidth][l]{\RomanLabelSize\RomanLabelFont (v)} \\[-0.6em]
                    {\centering \TitleSize\TitleFont {$c_{\text{num}}/c_{\text{ph}}$} \par}
                    {\MinMaxSize\MinMaxFont \makebox[\linewidth]{\hspace{\ValueHSpace}Min = -2.91e+05\hfill Max = 5.07e+04\hspace{\ValueHSpace}}} \\[0.2em]
                    \includegraphics[width=\linewidth, trim=20 2 20 20, clip]{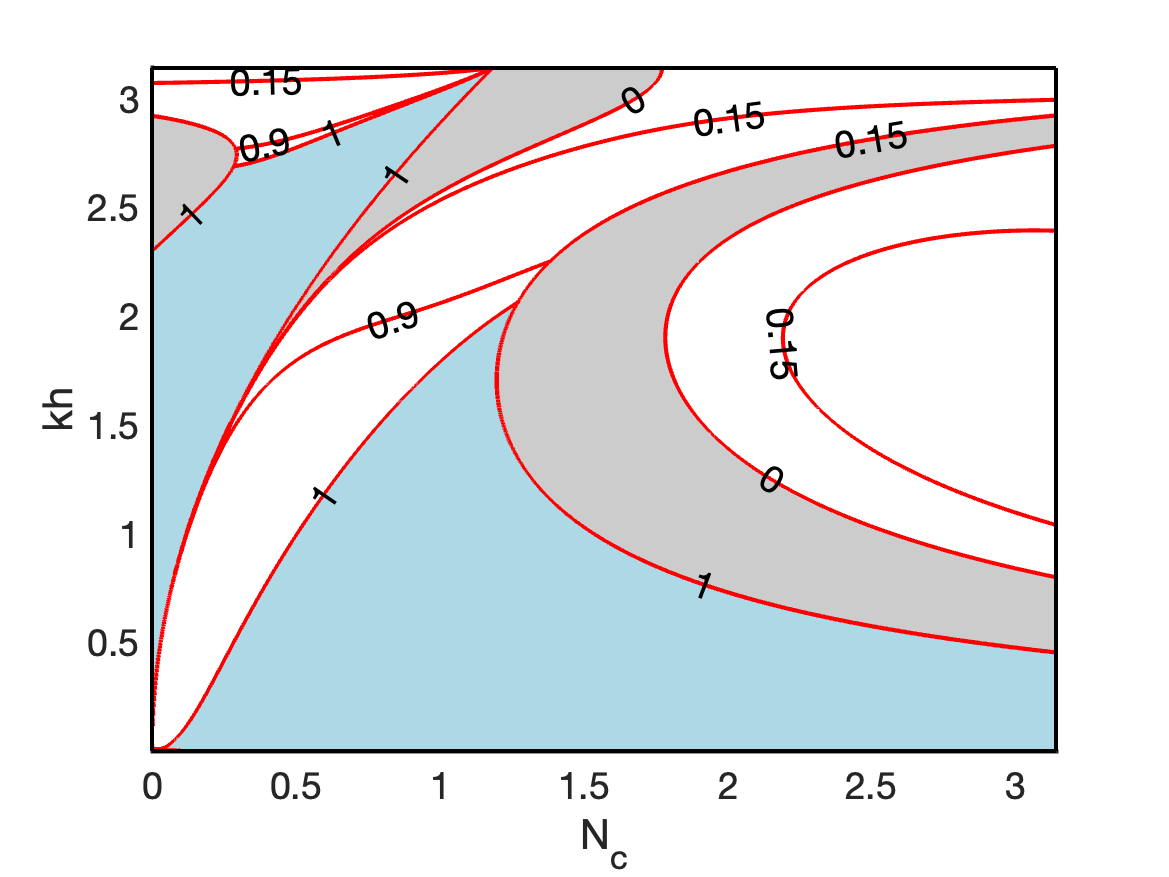}
                \end{minipage}
            }
        \end{minipage}
        &
        \hspace{\ColumnSpace} 
        \begin{minipage}{0.08\textwidth}
            \includegraphics[width=\textwidth]{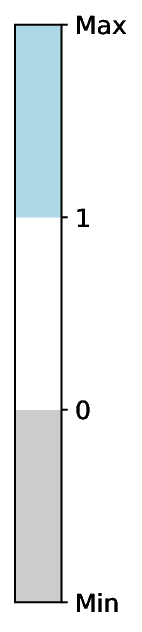}
        \end{minipage}
        \\[1em]

        \begin{minipage}{\mympwidth}
            \fbox{
                \begin{minipage}{\textwidth}
                    \makebox[\linewidth][l]{\RomanLabelSize\RomanLabelFont (iii)} \\[-0.6em]
                    {\centering \TitleSize\TitleFont {$v_{\text{g,num}}/v_{\text{g}}$} \par}
                    {\MinMaxSize\MinMaxFont \makebox[\linewidth]{\hspace{\ValueHSpace}Min = -4.03e+06\hfill Max = 1.16e+06\hspace{\ValueHSpace}}} \\[0.2em]
                    \includegraphics[width=\linewidth, trim=20 2 20 20, clip]{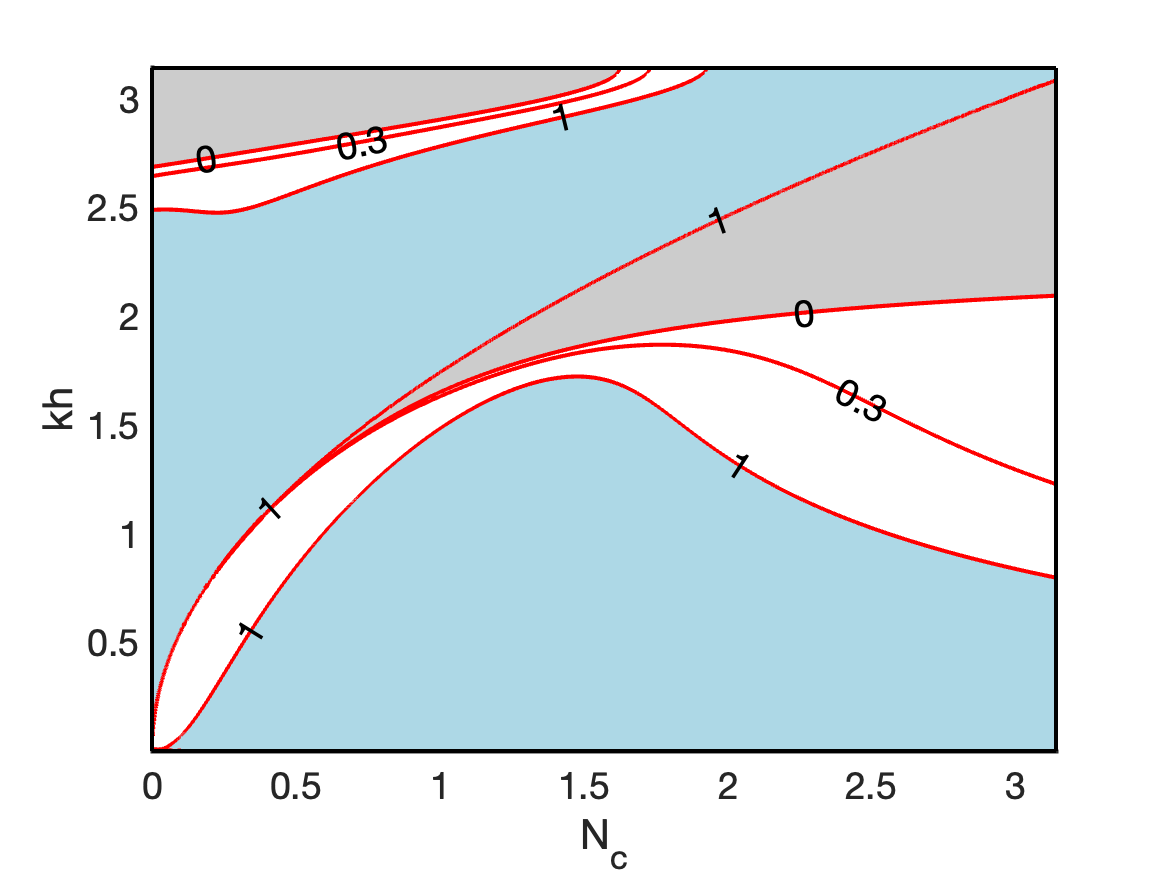}
                \end{minipage}
                }
        \end{minipage}
        &
        \hspace{\ColumnSpace} 
        \begin{minipage}{\mympwidth}
            \fbox{
                \begin{minipage}{\textwidth}
                    \makebox[\linewidth][l]{\RomanLabelSize\RomanLabelFont (vi)} \\[-0.6em]
                    {\centering \TitleSize\TitleFont {$v_{\text{g,num}}/v_{\text{g}}$} \par}
                    {\MinMaxSize\MinMaxFont \makebox[\linewidth]{\hspace{\ValueHSpace}Min = -3.35e+06\hfill Max = 1.36e+06\hspace{\ValueHSpace}}} \\[0.2em]
                    \includegraphics[width=\linewidth, trim=20 2 20 20, clip]{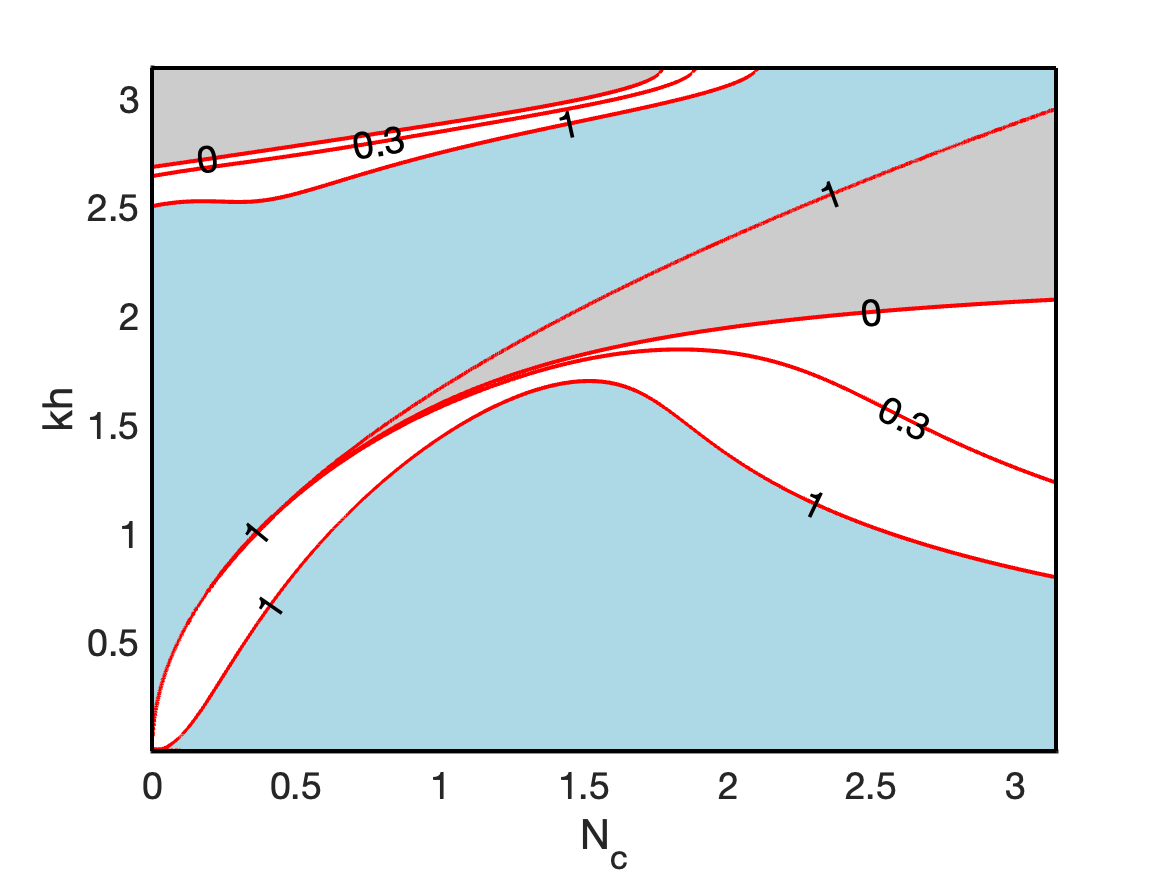}
                \end{minipage}
                }
        \end{minipage}
        &
        \hspace{\ColumnSpace} 
        \begin{minipage}{0.08\textwidth}
            \includegraphics[width=\textwidth]{Figures/3_colorbar.eps}
        \end{minipage}
    \end{tabular}
    \end{minipage}

    \caption{Contour plots of $|G_{\text{num}}|$, $\dfrac{c_{\text{num}}}{c_{\text{ph}}}$, and $\dfrac{v_{\text{g,num}}}{v_{\text{g}}}$ in the $(N_c, kh)$ plane for Eq.~\eqref{eqn:Main}, computed using the SSPRK3-CNCS6 scheme with the indicated values of $D_{\alpha}$.}
    \label{Fig:A1}
\end{figure}

Fig.~\ref{Fig:A1} shows the contours of the  \( |G_{\text{num}}| \), \( \dfrac{c_{\text{num}}}{c_{\text{ph}}} \), and \( \dfrac{v_{\text{g,num}}}{v_{\text{g}}} \) in $(N_{c}, kh)$ plane for the scheme SSPRK3-CNCS6. We have compared these properties for two distinct dispersion numbers $D_{\alpha} = \{0.11, 0.12\}$. In the idealized continuum limit corresponding to $kh \to 0$ and $N_c \to 0$, all numerical properties are expected to converge to their exact analytical counterparts, ensuring consistency. However, deviations from unity in these normalized quantities are expected as the solution moves away from this limit. For the \( |G_{\text{num}}| \), which characterizes numerical stability, the scheme exhibits distinct behaviors for the two dispersion numbers. The light blue shaded regions represent the area where \(|G_{\text{num}}|>1\), which is the unstable region. It is observed from Fig.~\ref{Fig:A1}(i) that for $D_{\alpha} = 0.11$, SSPRK3-CNCS6 remains stable throughout the Nyquist range ($kh = \pi$). However, for $D_{\alpha} = 0.12$, the stability range reduces to $kh \approx 2.51$ for $N_c \leq 1.38$. Consequently, $D_\alpha = 0.12$ is identified as the critical dispersion number, denoted by $D_{\alpha,\mathrm{cr}}$, beyond which the scheme may exhibit instability. Fig.~\ref{Fig:A1}(ii)-\ref{Fig:A1}(vi) shows the contours for  \( \dfrac{c_{\text{num}}}{c_{\text{ph}}} \), and \( \dfrac{v_{\text{g,num}}}{v_{\text{g}}} \) which quantify the phase and dispersion error. Here, the grey-shaded region represents negative contour values, indicating the presence of spurious waves. These spurious phenomena are referred to as $q$-waves~\cite{sengupta2012spurious}, and their appearance is purely a numerical artifact arising from specific combinations of discretization parameters. From Figs.~\ref{Fig:A1}(ii) and \ref{Fig:A1}(v), it is evident that increasing \( D_{\alpha} \) from 0.11 to 0.12 leads to an expansion of both the unstable region and the presence of reversed \( q \)-waves. This effect is particularly pronounced in the arc-shaped region around \( kh \in (2.4, 2.9) \), where the area affected by instability and negative wave propagation noticeably increases. Similar behavior is also observed for the normalized group velocity \( \dfrac{v_{\text{g,num}}}{v_{\text{g}}} \), as shown in Figs.~\ref{Fig:A1}(iii) and \ref{Fig:A1}(iv). As \( D_{\alpha} \) increases from 0.11 to 0.12, the unstable region expands significantly, particularly around \( kh = 0.15 \), indicating a broader range of \( q \)-waves experiencing inaccurate energy transport. This expansion of the negative and unstable contours suggests a degradation in the scheme’s ability to correctly capture group behavior, highlighting the sensitivity of energy propagation to the choice of $D_{\alpha}$.

\subsubsection{SSPRK3-CNCS8}
\begin{figure}[htbp!]
    \centering
    \begin{minipage}{\textwidth}
    \hspace*{0cm} 
    \begin{tabular}{ccc}
        \multicolumn{3}{c}{\makebox[0pt]{\hspace{0em}\textcolor{red}{\textbf{SSPRK3-CNCS8}}}} \\
        \multicolumn{1}{c}{\textbf{$\bm{D}_{\bm{\alpha}} = \bm{0.11}$}} & \textbf{$\bm{D}_{\bm{\alpha,\mathrm{cr}}} = \bm{0.12}$} \\[0.5em]

        \begin{minipage}{\mympwidth}
            \fbox{
               \begin{minipage}{\textwidth}
                    \makebox[\linewidth][l]{\RomanLabelSize\RomanLabelFont (i)} \\[-0.6em]
                    {\centering \TitleSize\TitleFont {$|G_{\text{num}}|$} \par}
                    {\MinMaxSize\MinMaxFont \makebox[\linewidth]{\hspace{\ValueHSpace}Min = 0.94\hfill Max = 25.16\hspace{\ValueHSpace}}} \\[0.2em]
                    \includegraphics[width=\linewidth, trim=20 2 20 20, clip]{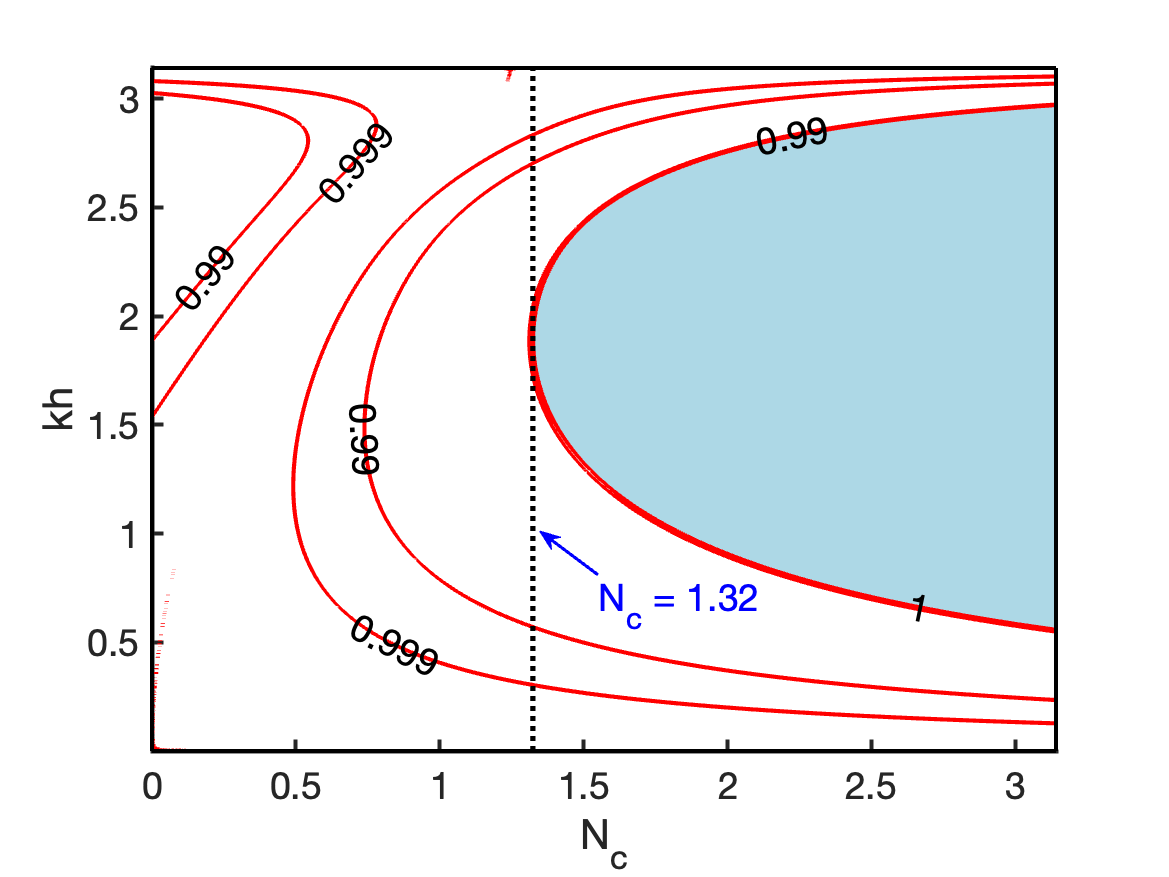}
               \end{minipage}
            }
        \end{minipage}
        &
        \hspace{\ColumnSpace}
        \begin{minipage}{\mympwidth}
            \fbox{
                \begin{minipage}{\textwidth}
                    \makebox[\linewidth][l]{\RomanLabelSize\RomanLabelFont (iv)} \\[-0.6em]
                    {\centering \TitleSize\TitleFont {$|G_{\text{num}}|$} \par}
                    {\MinMaxSize\MinMaxFont \makebox[\linewidth]{\hspace{\ValueHSpace}Min = 0.94\hfill Max = 23.68\hspace{\ValueHSpace}}} \\[0.2em]
                    \includegraphics[width=\linewidth, trim=20 2 20 20, clip]{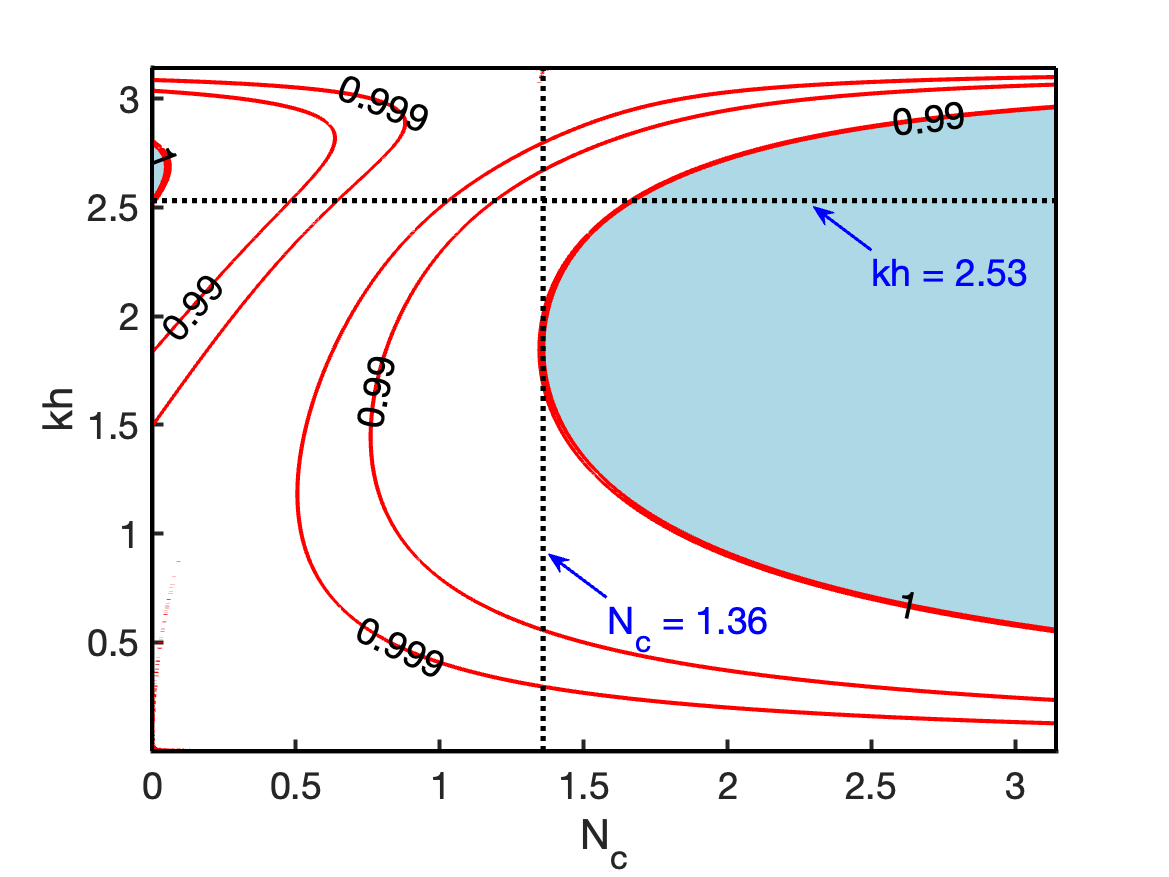}
                \end{minipage}
            }
        \end{minipage}
        &
        \hspace{\ColumnSpace}
        \begin{minipage}{0.08\textwidth}
            \includegraphics[width=\textwidth]{Figures/2_colorbar.eps}
        \end{minipage}
        \\[1em]

        \begin{minipage}{\mympwidth}
            \fbox{
                \begin{minipage}{\textwidth}
                    \makebox[\linewidth][l]{\RomanLabelSize\RomanLabelFont (ii)} \\[-0.6em]
                    {\centering \TitleSize\TitleFont {$c_{\text{num}}/c_{\text{ph}}$} \par}
                    {\MinMaxSize\MinMaxFont \makebox[\linewidth]{\hspace{\ValueHSpace}Min = -4.59e+04\hfill Max = 2.80e+04\hspace{\ValueHSpace}}} \\[0.2em]
                    \includegraphics[width=\linewidth, trim=20 2 20 20, clip]{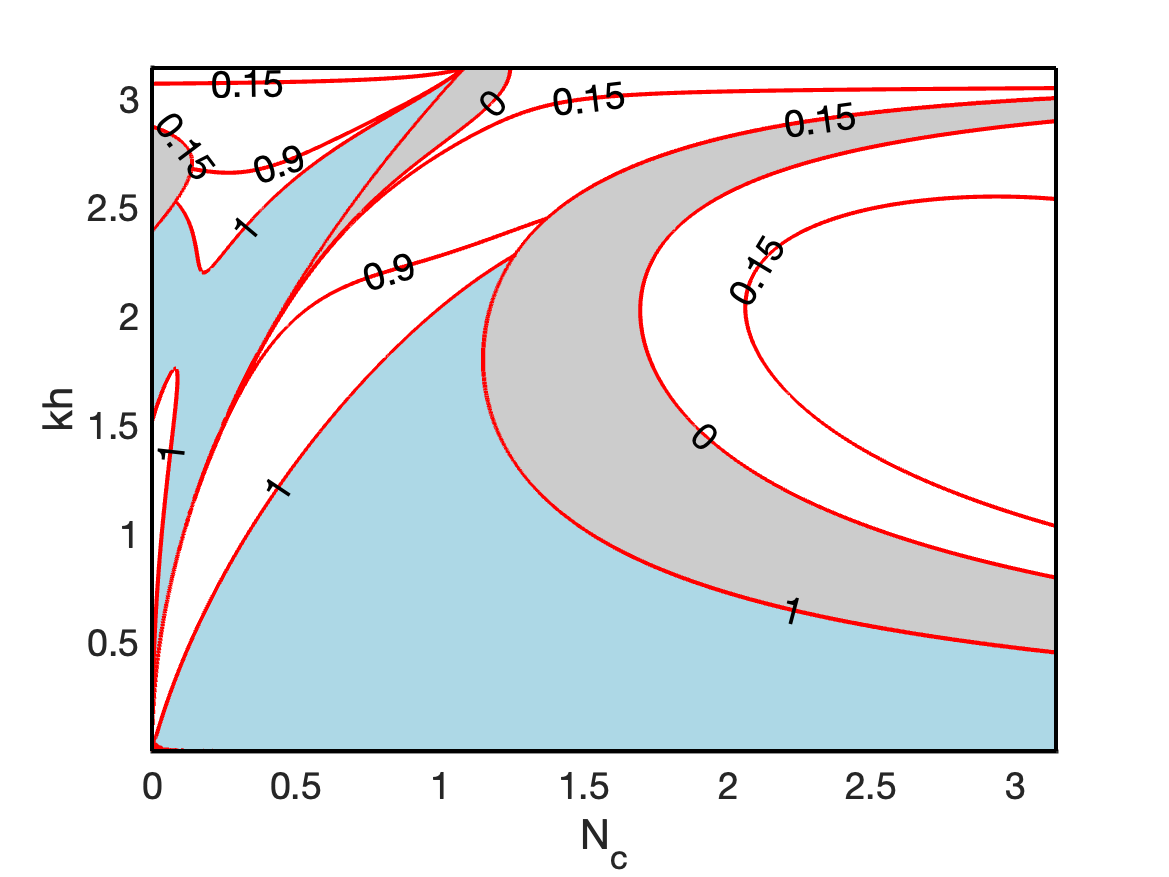}
                \end{minipage}
            }
        \end{minipage}
        &
        \hspace{\ColumnSpace}
        \begin{minipage}{\mympwidth}
            \fbox{
                \begin{minipage}{\textwidth}
                    \makebox[\linewidth][l]{\RomanLabelSize\RomanLabelFont (v)} \\[-0.6em]
                    {\centering \TitleSize\TitleFont {$c_{\text{num}}/c_{\text{ph}}$} \par}
                    {\MinMaxSize\MinMaxFont \makebox[\linewidth]{\hspace{\ValueHSpace}Min = -9.66e+04\hfill Max = 1.65e+04\hspace{\ValueHSpace}}} \\[0.2em]
                    \includegraphics[width=\linewidth, trim=20 2 20 20, clip]{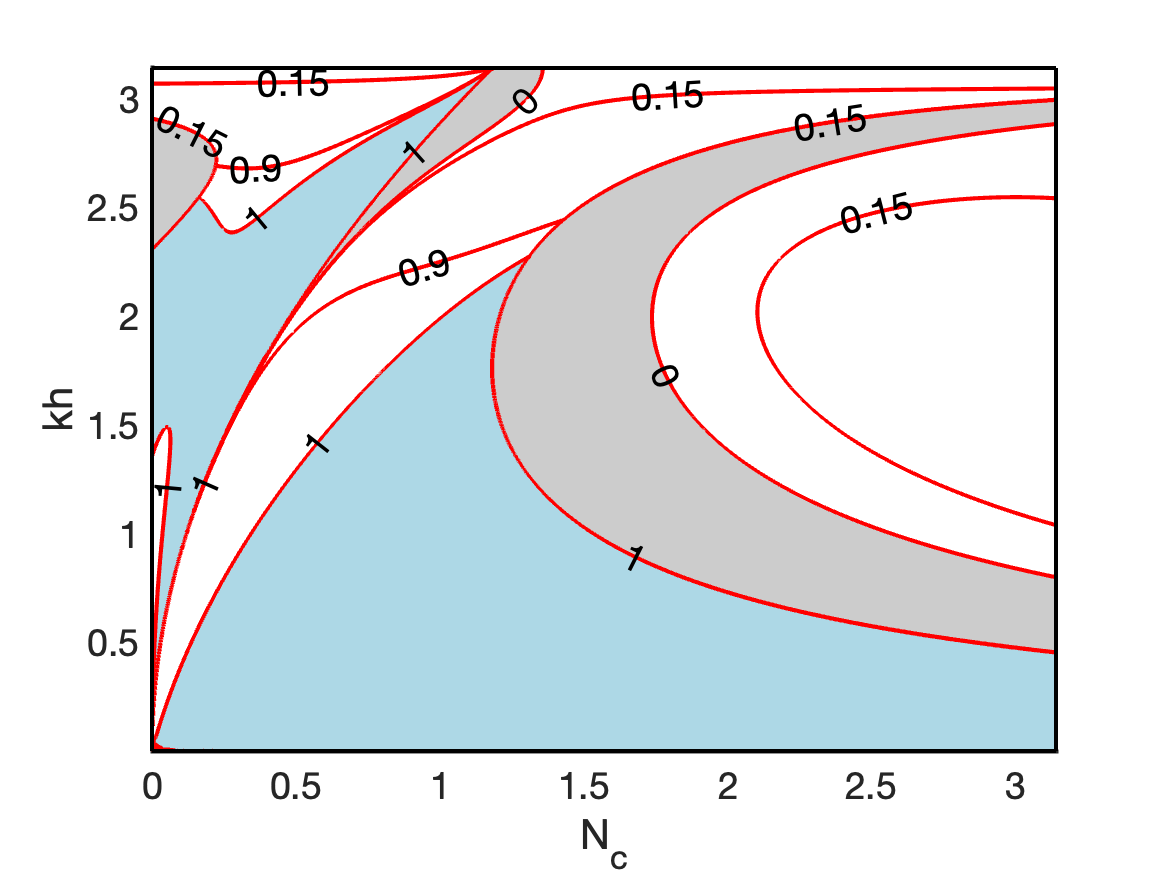}
                \end{minipage}
            }
        \end{minipage}
        &
        \hspace{\ColumnSpace}
        \begin{minipage}{0.08\textwidth}
            \includegraphics[width=\textwidth]{Figures/3_colorbar.eps}
        \end{minipage}
        \\[1em]

        \begin{minipage}{\mympwidth}
            \fbox{
                \begin{minipage}{\textwidth}
                    \makebox[\linewidth][l]{\RomanLabelSize\RomanLabelFont (iii)} \\[-0.6em]
                    {\centering \TitleSize\TitleFont {$v_{\text{g,num}}/v_{\text{g}}$} \par}
                    {\MinMaxSize\MinMaxFont \makebox[\linewidth]{\hspace{\ValueHSpace}Min = -2.51e+06\hfill Max = 1.76e+06\hspace{\ValueHSpace}}} \\[0.2em]
                    \includegraphics[width=\linewidth, trim=20 2 20 20, clip]{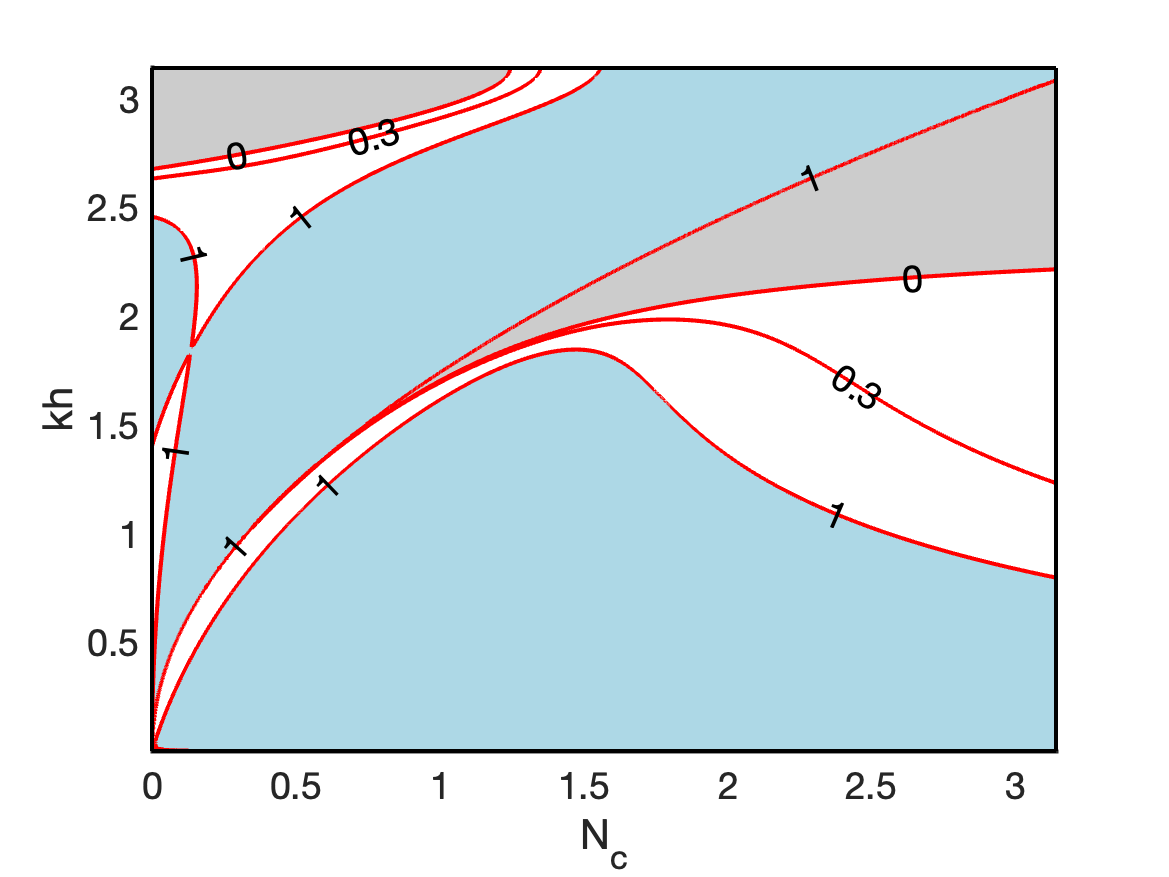}
                \end{minipage}
                }
        \end{minipage}
        &
        \hspace{\ColumnSpace}
        \begin{minipage}{\mympwidth}
            \fbox{
                \begin{minipage}{\textwidth}
                    \makebox[\linewidth][l]{\RomanLabelSize\RomanLabelFont (vi)} \\[-0.6em]
                    {\centering \TitleSize\TitleFont {$v_{\text{g,num}}/v_{\text{g}}$} \par}
                    {\MinMaxSize\MinMaxFont \makebox[\linewidth]{\hspace{\ValueHSpace}Min = -2.76e+06\hfill Max = 5.88e+05\hspace{\ValueHSpace}}} \\[0.2em]
                    \includegraphics[width=\linewidth, trim=20 2 20 20, clip]{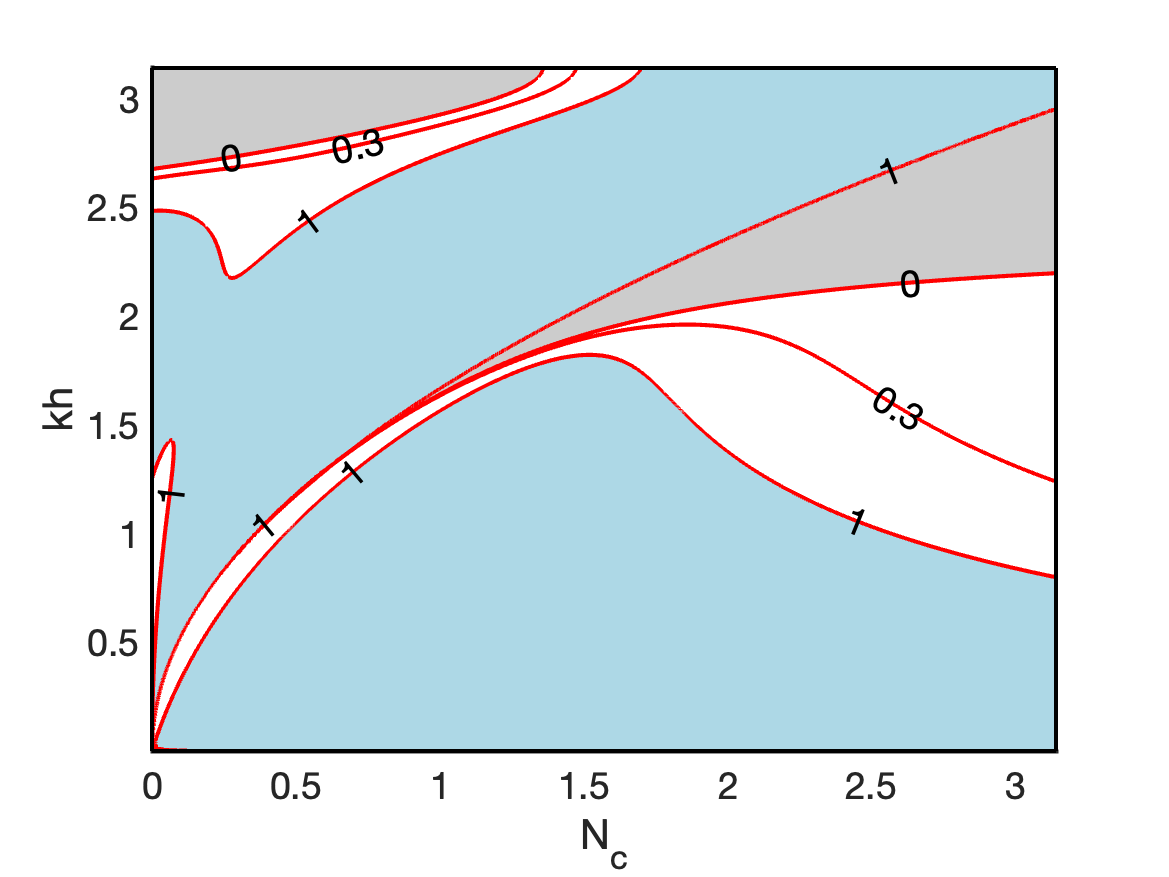}
                \end{minipage}
                }
        \end{minipage}
        &
        \hspace{\ColumnSpace}
        \begin{minipage}{0.08\textwidth}
            \includegraphics[width=\textwidth]{Figures/3_colorbar.eps}
        \end{minipage}
    \end{tabular}
    \end{minipage}
    \caption{Contour plots of $|G_{\text{num}}|$, $\dfrac{c_{\text{num}}}{c_{\text{ph}}}$, and $\dfrac{v_{\text{g,num}}}{v_{\text{g}}}$ in the $(N_c, kh)$ plane for Eq.~\eqref{eqn:Main}, computed using the SSPRK3-CNCS8 scheme with the indicated values of $D_{\alpha}$.}
    \label{Fig:A2}
\end{figure}

The spectral analysis is extended to the eighth-order scheme (SSPRK3-CNCS8) to examine its stability and dispersion characteristics. Fig.~\ref{Fig:A2}(i) presents the computed property charts for the dispersion number $D_{\alpha} = 0.11$, which corresponds to the highest stable value for this scheme. In this case, numerical stability is maintained over a Courant number range $N_c \leq 1.32$, with the scheme capable of resolving the entire admissible wavenumber domain up to $kh = \pi$. A slight increment in the dispersion number to $D_{\alpha,\mathrm{cr}} = 0.12$  results in the property charts shown in Fig.~\ref{Fig:A2}(iv), where stability is preserved up to $N_c = 1.36$. However, the resolved spectrum is restricted to wavenumbers $kh \leq 2.53$. Any further increase in $D_{\alpha}$ beyond this threshold leads to instability, as indicated by amplification factors exceeding unity.\par
Figs.~\ref{Fig:A2}(ii),~\ref{Fig:A2}(v),~\ref{Fig:A2}(iii), and~\ref{Fig:A2}(vi) display the corresponding contours for the normalized phase speed and group velocity. The qualitative trends observed in these plots are consistent with those obtained for the sixth-order scheme, revealing similar regions of phase reversal and non-physical energy propagation. These results show that spurious numerical behavior, such as \( q \)-waves, can persist even at higher spatial accuracy unless the numerical parameters are carefully chosen.

\subsubsection{SSPRK3-CCS8}
\begin{figure}[htbp!]
    \centering
    \begin{minipage}{\textwidth}
    \hspace*{0cm}
    \begin{tabular}{ccc}
        \multicolumn{3}{c}{\makebox[0pt]{\hspace{0em}\textcolor{red}{\textbf{SSPRK3-CCS8}}}} \\
        \multicolumn{1}{c}{\textbf{$\bm{D}_{\bm{\alpha}} = \bm{0.011}$}} & \textbf{$\bm{D}_{\bm{\alpha,\mathrm{cr}}} = \bm{0.012}$} & \\[0.5em]

        \begin{minipage}{\mympwidth}
            \fbox{
               \begin{minipage}{\textwidth}
                    \makebox[\linewidth][l]{\RomanLabelSize\RomanLabelFont (i)} \\[-0.6em]
                    {\centering \TitleSize\TitleFont {$|G_{\text{num}}|$} \par}
                    {\MinMaxSize\MinMaxFont \makebox[\linewidth]{\hspace{\ValueHSpace}Min = 0.94\hfill Max = 2.60e+02\hspace{\ValueHSpace}}} \\[0.2em]
                    \includegraphics[width=\linewidth, trim=20 2 20 20, clip]{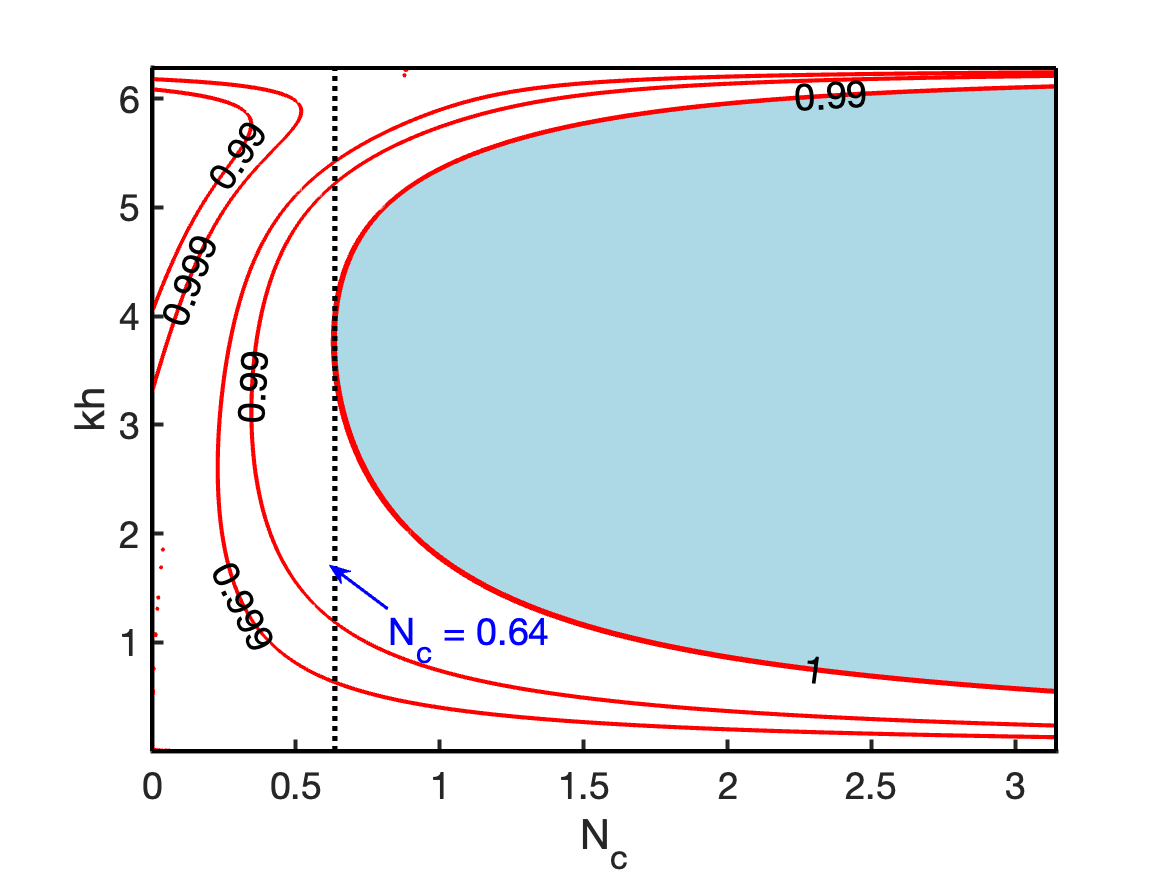}
               \end{minipage}
            }
        \end{minipage}
        &
        \hspace{\ColumnSpace}
        \begin{minipage}{\mympwidth}
            \fbox{
                \begin{minipage}{\textwidth}
                    \makebox[\linewidth][l]{\RomanLabelSize\RomanLabelFont (iv)} \\[-0.6em]
                    {\centering \TitleSize\TitleFont {$|G_{\text{num}}|$} \par}
                    {\MinMaxSize\MinMaxFont \makebox[\linewidth]{\hspace{\ValueHSpace}Min = 0.94\hfill Max = 2.54e+02\hspace{\ValueHSpace}}} \\[0.2em]
                    \includegraphics[width=\linewidth, trim=20 2 20 20, clip]{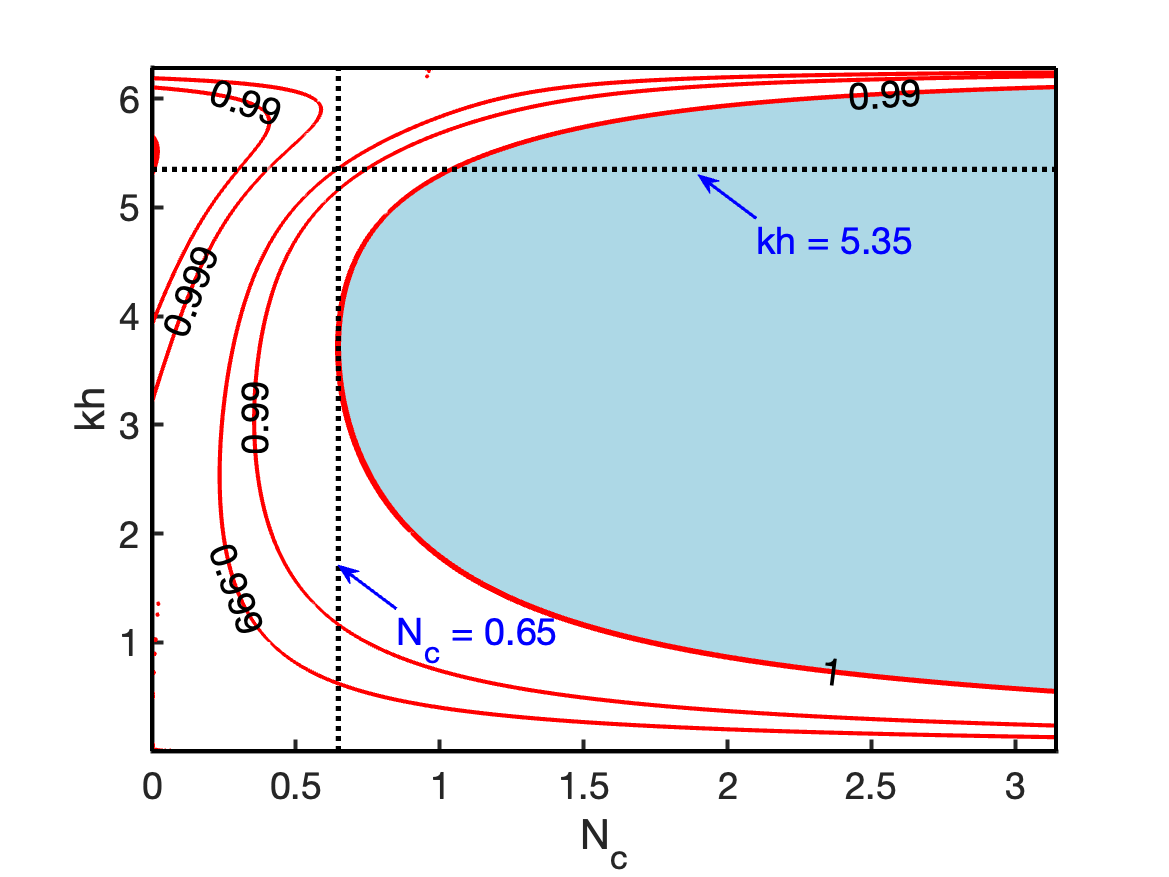}
                \end{minipage}
            }
        \end{minipage}
        &
        \hspace{\ColumnSpace}
        \begin{minipage}{0.08\textwidth}
            \includegraphics[width=\textwidth]{Figures/2_colorbar.eps}
        \end{minipage}
        \\[1em]

        \begin{minipage}{\mympwidth}
            \fbox{
                \begin{minipage}{\textwidth}
                    \makebox[\linewidth][l]{\RomanLabelSize\RomanLabelFont (ii)} \\[-0.6em]
                    {\centering \TitleSize\TitleFont {$c_{\text{num}}/c_{\text{ph}}$} \par}
                    {\MinMaxSize\MinMaxFont \makebox[\linewidth]{\hspace{\ValueHSpace}Min = -1.81e+04 \hfill Max = 2.47e+05\hspace{\ValueHSpace}}} \\[0.2em]
                    \includegraphics[width=\linewidth, trim=20 2 20 20, clip]{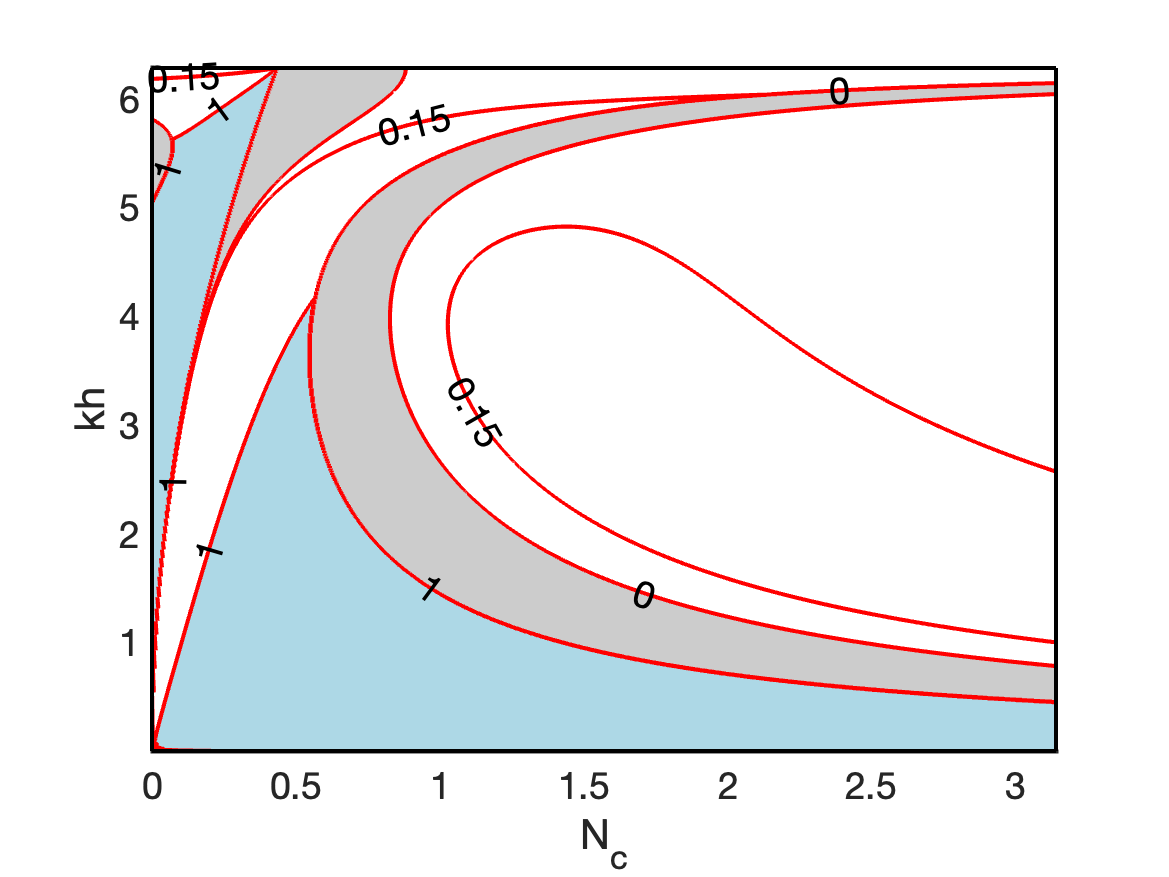}
                \end{minipage}
            }
        \end{minipage}
        &
        \hspace{\ColumnSpace}
        \begin{minipage}{\mympwidth}
            \fbox{
                \begin{minipage}{\textwidth}
                    \makebox[\linewidth][l]{\RomanLabelSize\RomanLabelFont (v)} \\[-0.6em]
                    {\centering \TitleSize\TitleFont {$c_{\text{num}}/c_{\text{ph}}$} \par}
                    {\MinMaxSize\MinMaxFont \makebox[\linewidth]{\hspace{\ValueHSpace}Min = -3.51e+03\hfill Max = 1.16e+04\hspace{\ValueHSpace}}} \\[0.2em]
                    \includegraphics[width=\linewidth, trim=20 2 20 20, clip]{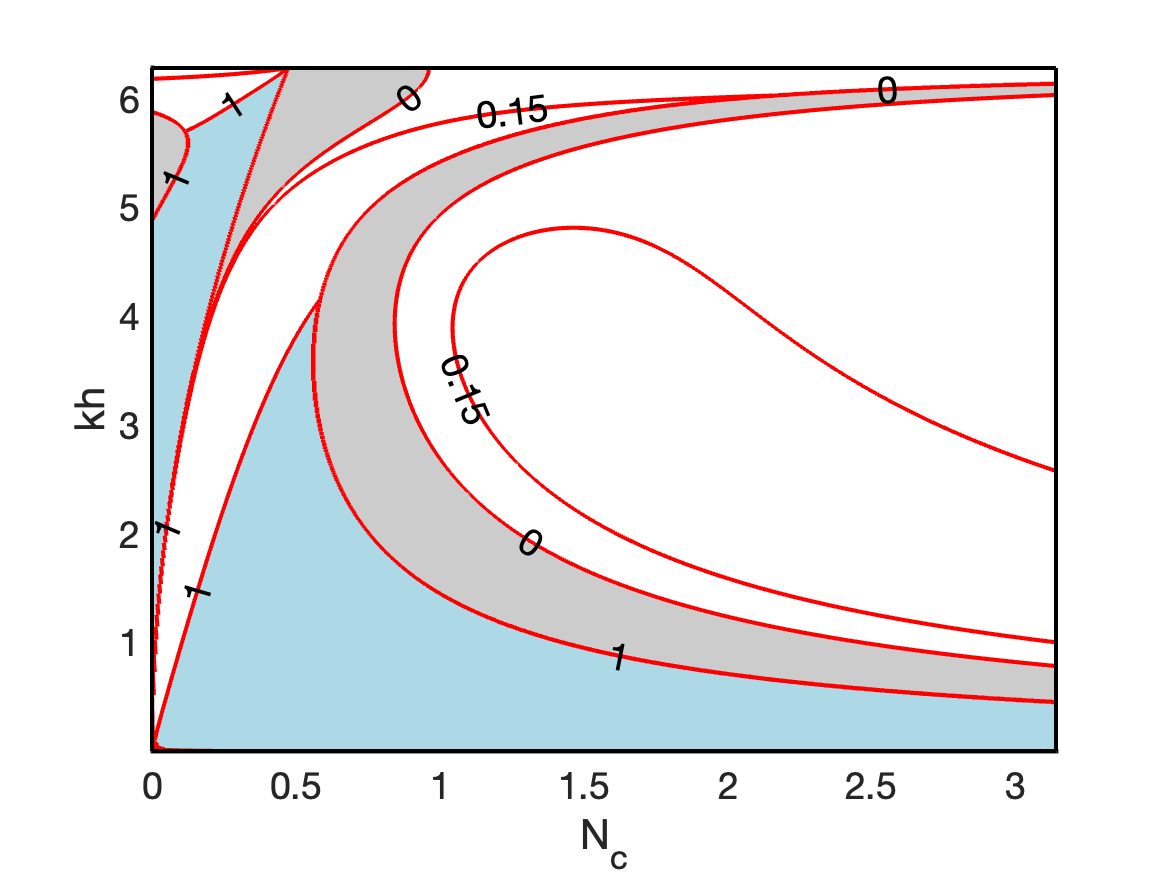}
                \end{minipage}
            }
        \end{minipage}
        &
        \hspace{\ColumnSpace}
        \begin{minipage}{0.08\textwidth}
            \includegraphics[width=\textwidth]{Figures/3_colorbar.eps}
        \end{minipage}
        \\[1em]

        \begin{minipage}{\mympwidth}
            \fbox{
                \begin{minipage}{\textwidth}
                    \makebox[\linewidth][l]{\RomanLabelSize\RomanLabelFont (iii)} \\[-0.6em]
                    {\centering \TitleSize\TitleFont {$v_{\text{g,num}}/v_{\text{g}}$} \par}
                    {\MinMaxSize\MinMaxFont \makebox[\linewidth]{\hspace{\ValueHSpace}Min = -8.33e+06\hfill Max = 3.19e+06\hspace{\ValueHSpace}}} \\[0.2em]
                    \includegraphics[width=\linewidth, trim=20 2 20 20, clip]{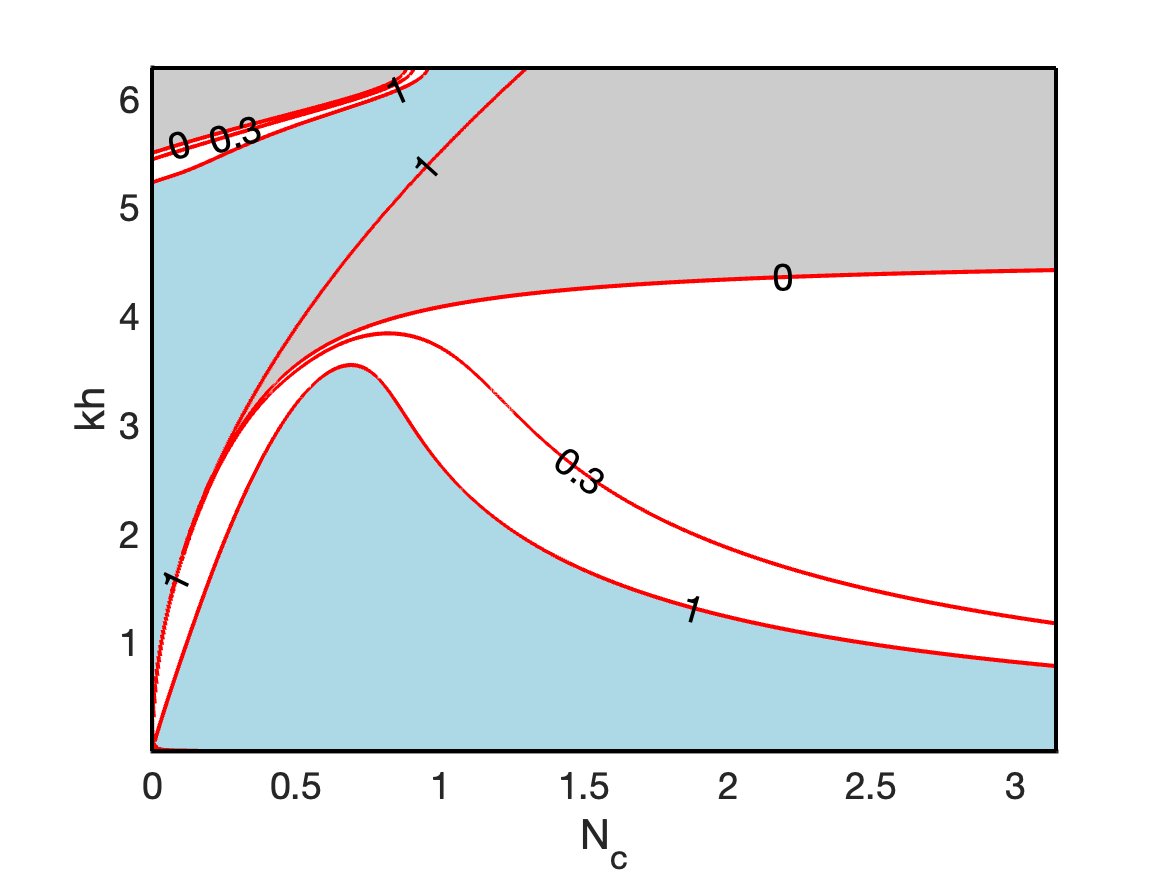}
                \end{minipage}
            }
        \end{minipage}
        &
        \hspace{\ColumnSpace}
        \begin{minipage}{\mympwidth}
            \fbox{
                \begin{minipage}{\textwidth}
                    \makebox[\linewidth][l]{\RomanLabelSize\RomanLabelFont (vi)} \\[-0.6em]
                    {\centering \TitleSize\TitleFont {$v_{\text{g,num}}/v_{\text{g}}$} \par}
                    {\MinMaxSize\MinMaxFont \makebox[\linewidth]{\hspace{\ValueHSpace}Min = -3.46e+05\hfill Max = 1.01e+06\hspace{\ValueHSpace}}} \\[0.2em]
                    \includegraphics[width=\linewidth, trim=20 2 20 20, clip]{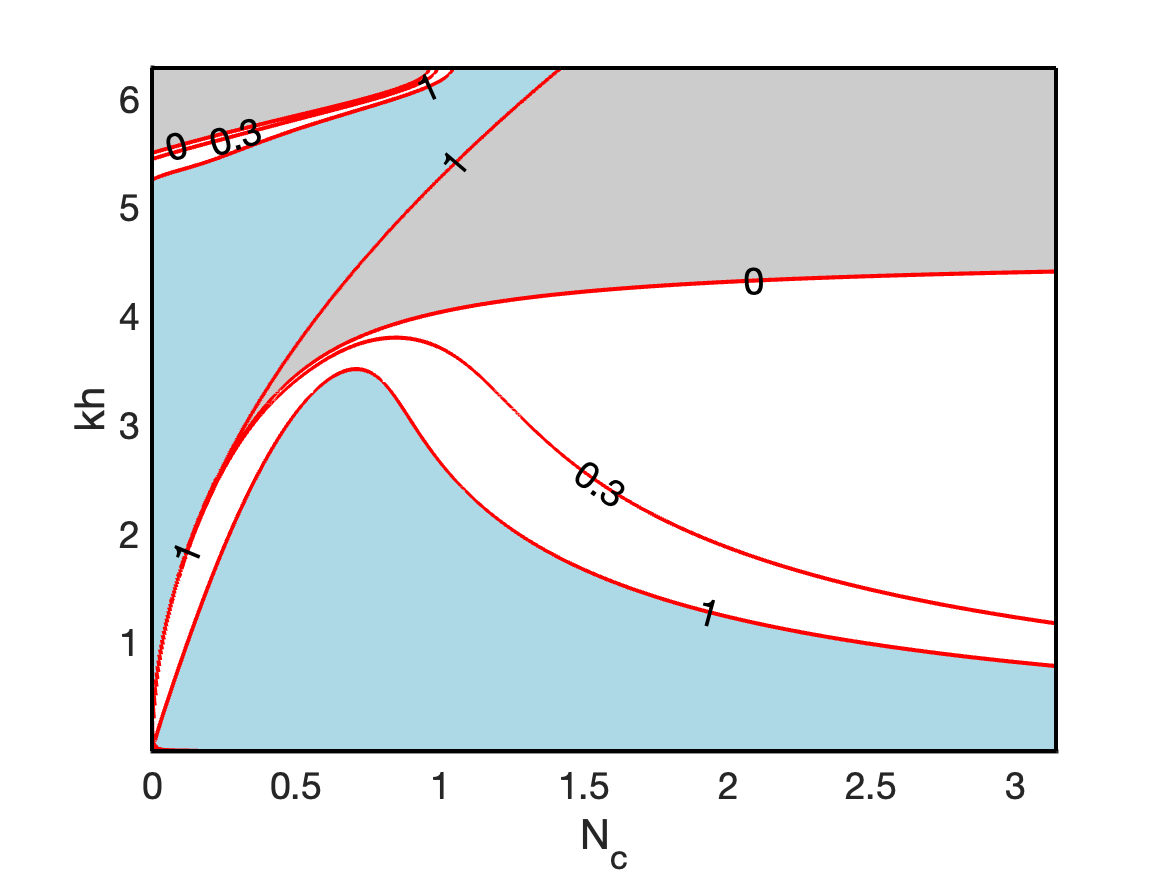}
                \end{minipage}
            }
        \end{minipage}
        &
        \hspace{\ColumnSpace}
        \begin{minipage}{0.08\textwidth}
            \includegraphics[width=\textwidth]{Figures/3_colorbar.eps}
        \end{minipage}
    \end{tabular}
    \end{minipage}
    \caption{Contour plots of $|G_{\text{num}}|$, $\dfrac{c_{\text{num}}}{c_{\text{ph}}}$, and $\dfrac{v_{\text{g,num}}}{v_{\text{g}}}$ in the $(N_c, kh)$ plane for Eq.~\eqref{eqn:Main}, computed using the SSPRK3-CCS8 scheme with the indicated values of $D_{\alpha}$.}
    \label{Fig:A3}
\end{figure}

Next, we analyzed the eighth-order cell-centered scheme (SSPRK3-CCS8). Fig.~\ref{Fig:A3}(i) illustrates the numerical properties corresponding to the maximum stable dispersion number, \( D_{\alpha} = 0.011 \), which is approximately an order of magnitude smaller than that employed for the SSPRK3-CNCS6. At this value of \( D_{\alpha} \), the CCS8 remains stable for Courant numbers up to \( N_c \leq 0.64 \) and is capable of resolving the full spectral range up to the Nyquist limit \( kh = 2\pi \). The extended spectral range arises from the cell-centered discretization, where the degrees of freedom are located at cell centers rather than grid nodes, allowing the scheme to represent Fourier modes over a doubled wavenumber interval compared to node-centered schemes. A modest increase in the dispersion number to $D_{\alpha,\mathrm{cr}}= 0.012$   yields the property charts shown in Fig.~\ref{Fig:A3}(iv), where stability is preserved up to $N_c = 0.65$, and the resolvable spectrum extends to $kh \approx 5.35$. Beyond this critical value, instability arises, as evidenced by amplification factors exceeding unity.

\par
The results indicate that the cell-centered scheme provides superior resolution of higher wavenumbers, even at relatively lower Courant numbers. In contrast, the node-centered compact scheme maintains stability at higher values of $N_c$ but supports a narrower wavenumber range. The normalized phase speed and group velocity contours for $D_{\alpha} = 0.011$ and $D_{\alpha,\mathrm{cr}} = 0.012$ are depicted in Figs.~\ref{Fig:A3}(ii)--\ref{Fig:A3}(iii) and Figs.~\ref{Fig:A3}(v)--\ref{Fig:A3}(vi), respectively. While these contours reflect qualitatively similar behavior to those of the CNCS8, they extend over a considerably wider wavenumber range, highlighting the enhanced spectral resolution capabilities of the CCS8 formulation.
\subsection{Spectral analysis of linear 2-D convection dispersion equation}\label{Sec:2D_GSA}
In this section, we extend the analysis of the 1D convection-dispersion equation to the 2D case. The governing equation for the 2D linear convection-dispersion model with periodic boundary conditions is given by:
\begin{align}\label{eqn:Main2D}
    &u_t + c_x u_x + c_y u_y + \alpha \left( u_{xxx} + u_{yyy}  \right) = 0, \quad (x, y, t) \in \mathbb{R} \times \mathbb{R} \times \mathbb{R}^+, \nonumber \\
    &u(x, y, 0) = u_{0},
\end{align}
\noindent where $c_x$ and $c_y$ represent the convection speeds along the $x$- and $y$-directions, respectively, while $\alpha$ denotes the dispersion coefficient, which is assumed to be isotropic for this study. Applying the GSA approach as in the one-dimensional case, we derive the fundamental expressions for the two-dimensional case. The physical dispersion relation is expressed as:
\begin{equation}
    \omega = c_x k_x + c_y k_y - \alpha(k_x^3 + k_y^3).
\end{equation}
\noindent The physical phase speed is derived from the dispersion relation as:
\begin{equation}
    c_{\text{\text{ph}}} = \dfrac{\omega}{\sqrt{k_x^2 + k_y^2}} = \dfrac{c_x k_x + c_y k_y - \alpha(k_x^3 + k_y^3)}{\sqrt{k_x^2 + k_y^2}}.
\end{equation}
\noindent The components of the physical group velocity are given by:
\begin{equation}
    v_{\text{g},x} = \dfrac{\partial \omega}{\partial k_x} = c_x - 3\alpha k_x^2,
\end{equation}
\begin{equation}
    v_{\text{g},y} = \dfrac{\partial \omega}{\partial k_y} = c_y - 3\alpha k_y^2.
\end{equation}

\noindent The physical amplification factor, expressed in terms of $N_{cz}$, and $D_{\alpha z}$, where $z = \{x, y\}$ denotes the spatial directions, is given by

\begin{equation}
    G_{\text{\text{ph}}} = e^{\iota [D_{\alpha x}(k_x h_x)^3 + D_{\alpha y}(k_y h_y)^3]} e^{-\iota [N_{cx} k_x h_x + N_{cy} k_y h_y]}. 
\end{equation}

\noindent Here, $N_{cz} = c_z \Delta t / \Delta z$ denotes the Courant number, with $c_z$ representing the wave speed in the $z$-direction, while $D_{\alpha z} = \alpha_z \Delta t / \Delta z^3$ denotes the dimensionless dispersion parameter, where $\alpha_z$ is the dispersion coefficient in the $z$-direction. For this study, it is assumed that $\alpha = \alpha_x = \alpha_y$. The numerical dispersion relation is given as:

\begin{equation}
    \omega_{\text{num}} = c_{\text{num}} \left(\sqrt{k_x^2 + k_y^2}\right)  - \alpha_{\text{num}}(k_x^3 + k_y^3),
\end{equation}

\noindent where $c_{\text{num}}$ and $\alpha_{\text{num}}$ vary depending on the chosen numerical discretization scheme. The numerical amplification factor for the two-dimensional convection–dispersion equation, discretized using a third-order SSP Runge–Kutta method, is expressed as
\begin{equation}
G_{\text{num}} = 1 - A + \frac{A^2}{2} - \frac{A^3}{6},
\end{equation}
where the operator \( A \) contains both convective and dispersive operator similar to 1D and for 2D it is defined as
\begin{equation}
A = N_{cx} \left( \iota \, k^{[1]}_{\text{eq},x} h_x \right) + N_{cy}\left(\iota \, k^{[1]}_{\text{eq},y} h_y \right)
+ D_{\alpha x} \left( -\iota (k^{[2]}_{\text{eq},x})^3 h_x^3 \right)+ D_{\alpha y} \left( -\iota (k^{[2]}_{\text{eq},y})^3 h_y^3 \right).
\end{equation}
Here, \( h_x \) and \( h_y \) are the spatial grid spacings, and \( k^{[1]}_{\text{eq},x} \), \( k^{[1]}_{\text{eq},y} \), \( k^{[2]}_{\text{eq},x} \), and \( k^{[2]}_{\text{eq},y} \) are the equivalent wavenumbers corresponding to first- and third-derivative approximations in the $x$- and $y$-directions, respectively. 

\noindent The numerical phase shift per time step is:

\begin{equation}
    \tan(\beta_{\text{num}}) = -\Bigl(\dfrac{(G_{\text{num}})_{\text{Imag}}}{(G_{\text{num}})_{\text{Real}}} \Bigr),
\end{equation}

\noindent where $\beta_{\text{num}} = \Biggl( c_{\text{num}} \Bigl(\sqrt{k_x^2 + k_y^2}\Bigr) - \alpha_{\text{num}}(k_x^3 + k_y^3)\Biggr) \Delta t$.\\

\noindent The dimensionless numerical phase speed is:

\begin{equation}
    \dfrac{c_{\text{num}}}{c_{\text{ph}}} = \dfrac{\beta_{\text{num}}}{\Bigl( N_{cx} (k_x h_x) + N_{cy} (k_y h_y)\Bigr) - \Bigl( D_{\alpha x} (k_x h_x)^3 + D_{\alpha y} (k_y h_y)^3\Bigr)}.
\end{equation}

\noindent The numerical group velocity can be derived from the numerical dispersion relation using $(v_{\text{g},z})_{\text{num}} = \dfrac{\partial}{\partial k_z} (\omega_{\text{num}})$, yielding:

\begin{equation}
    \dfrac{(v_{\text{g},x})_{\text{num}}}{v_{\text{g},x}} = \dfrac{1}{N_{cx} - 3 D_{\alpha x}(k_x h_x)^2} \dfrac{d \beta_{\text{num}}}{d (k_x h_x)},
\end{equation}

\begin{equation}
    \dfrac{(v_{\text{g},y})_{\text{num}}}{v_{\text{g},y}} = \dfrac{1}{N_{cy} - 3 D_{\alpha y}(k_y h_y)^2} \dfrac{d \beta_{\text{num}}}{d (k_y h_y)}.
\end{equation}

\noindent For accurate solutions, the quantities $ \dfrac{c_{\text{num}}}{c_{\text{ph}}} $, $ \dfrac{(v_{\text{g},x})_{\text{num}}}{v_{\text{g},x}} $, and $ \dfrac{(v_{\text{g},y})_{\text{num}}}{v_{\text{g},y}} $ should be as close to unity as possible. Following the 1D case, we extend the GSA to the 2D setting for three schemes: SSPRK3-CNCS6, SSPRK3-CNCS8, and SSPRK3-CCS8.
Here, this 2D spatial discretization introduces two key parameters: the spatial aspect ratio of the cell, $ AR = \dfrac{h_y}{h_x} $, and the imposed wave propagation angle~\cite{sengupta2011analysis}, $ \theta = \tan^{-1} \left( \dfrac{c_y}{c_x} \right) $. 
The aspect ratio governs the geometric anisotropy \cite{sengupta2011analysis} of the mesh and affects directional accuracy and numerical diffusion, particularly for non-square cells. The wave propagation angle \( \theta \) defines the orientation of the imposed wave with respect to the grid axes, and is crucial for analyzing the scheme’s directional dispersion and anisotropic behavior.
For all analyses, we have set the aspect ratio of the cell to $ AR = 1 $ and the wave propagation angle to $ \theta = 45^\circ $. To assess the performance of the two numerical methods, the unit domain is discretized using a uniform grid with 2000 equally spaced points.

\subsubsection{SSPRK3-CNCS6}
\begin{figure}[htbp!]
    \centering
    \begin{minipage}{\textwidth}
    \hspace*{0cm} 
    \begin{tabular}{ccc}
        \multicolumn{3}{c}{\makebox[0pt]{\hspace{0em}\textcolor{red}{\textbf{SSPRK3-CNCS6}}}} \\
        \multicolumn{1}{c}{\textbf{$\bm{D}_{\bm{\alpha}} = \bm{0.11}$}} & \textbf{$\bm{D}_{\bm{\alpha,\mathrm{cr}}} = \bm{0.12}$} \\[0.5em]

        \begin{minipage}{\mympwidth}
            \setlength{\fboxrule}{0.4pt} 
            \fcolorbox{black}{white}{
               \begin{minipage}{\textwidth}
                    \makebox[\linewidth][l]{\RomanLabelSize\RomanLabelFont (i)} \\[-0.6em]
                    {\centering \TitleSize\TitleFont {$|G_{\text{num}}|$} \par}
                    {\MinMaxSize\MinMaxFont \makebox[\linewidth]{\hspace{\ValueHSpace}Min = 0.94\hfill Max = 1.00 \hspace{\ValueHSpace}}} \\[0.2em]
                    \includegraphics[width=\linewidth, trim=20 2 20 20, clip]{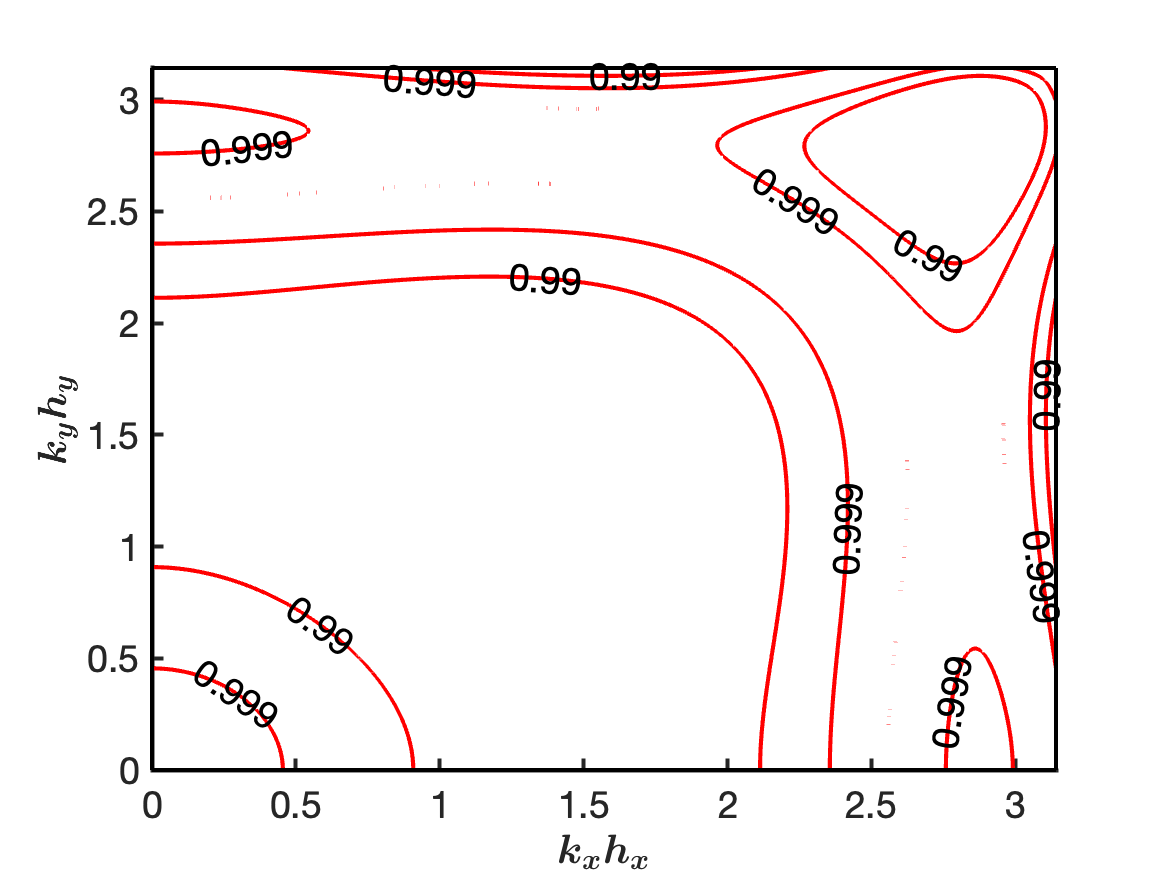}
               \end{minipage}
            }
        \end{minipage}
        &
        \hspace{\ColumnSpace} 
        \begin{minipage}{\mympwidth}
            \fbox{
                \begin{minipage}{\textwidth}
                    \makebox[\linewidth][l]{\RomanLabelSize\RomanLabelFont (iv)} \\[-0.6em]
                    {\centering \TitleSize\TitleFont {$|G_{\text{num}}|$} \par}
                    {\MinMaxSize\MinMaxFont \makebox[\linewidth]{\hspace{\ValueHSpace}Min = 0.94\hfill Max = 1.10\hspace{\ValueHSpace}}} \\[0.2em]
                    \includegraphics[width=\linewidth, trim=20 2 20 20, clip]{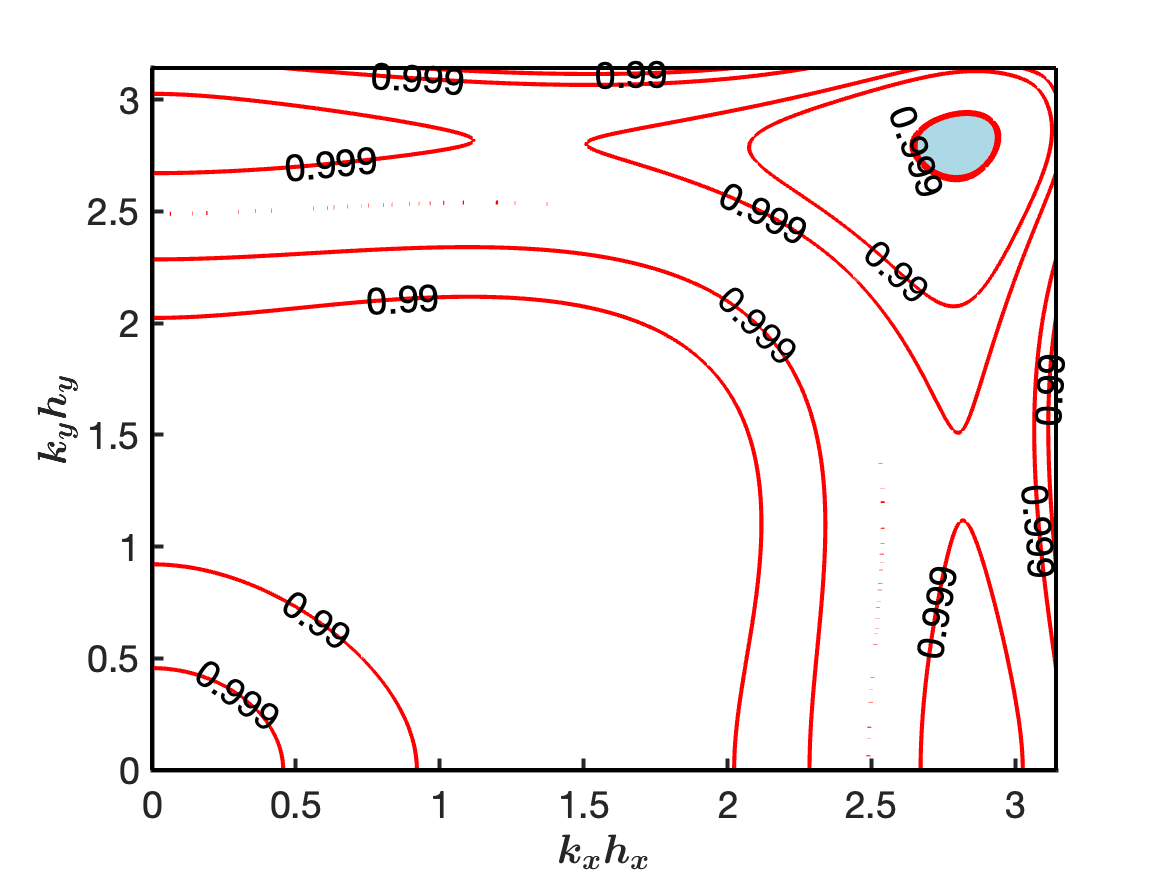}
                \end{minipage}
            }
        \end{minipage}
        &
        \hspace{\ColumnSpace} 
        \begin{minipage}{0.08\textwidth}
            \includegraphics[width=\textwidth]{Figures/2_colorbar.eps}
        \end{minipage}
        \\[1em]

        \begin{minipage}{\mympwidth}
            \fbox{
                \begin{minipage}{\textwidth}
                    \makebox[\linewidth][l]{\RomanLabelSize\RomanLabelFont (ii)} \\[-0.6em]
                    {\centering \TitleSize\TitleFont {$c_{\text{num}}/c_{\text{ph}}$} \par}
                    {\MinMaxSize\MinMaxFont \makebox[\linewidth]{\hspace{\ValueHSpace}Min = -4.40e+05\hfill Max = 8.74e+05\hspace{\ValueHSpace}}} \\[0.2em]
                    \includegraphics[width=\linewidth, trim=20 2 20 20, clip]{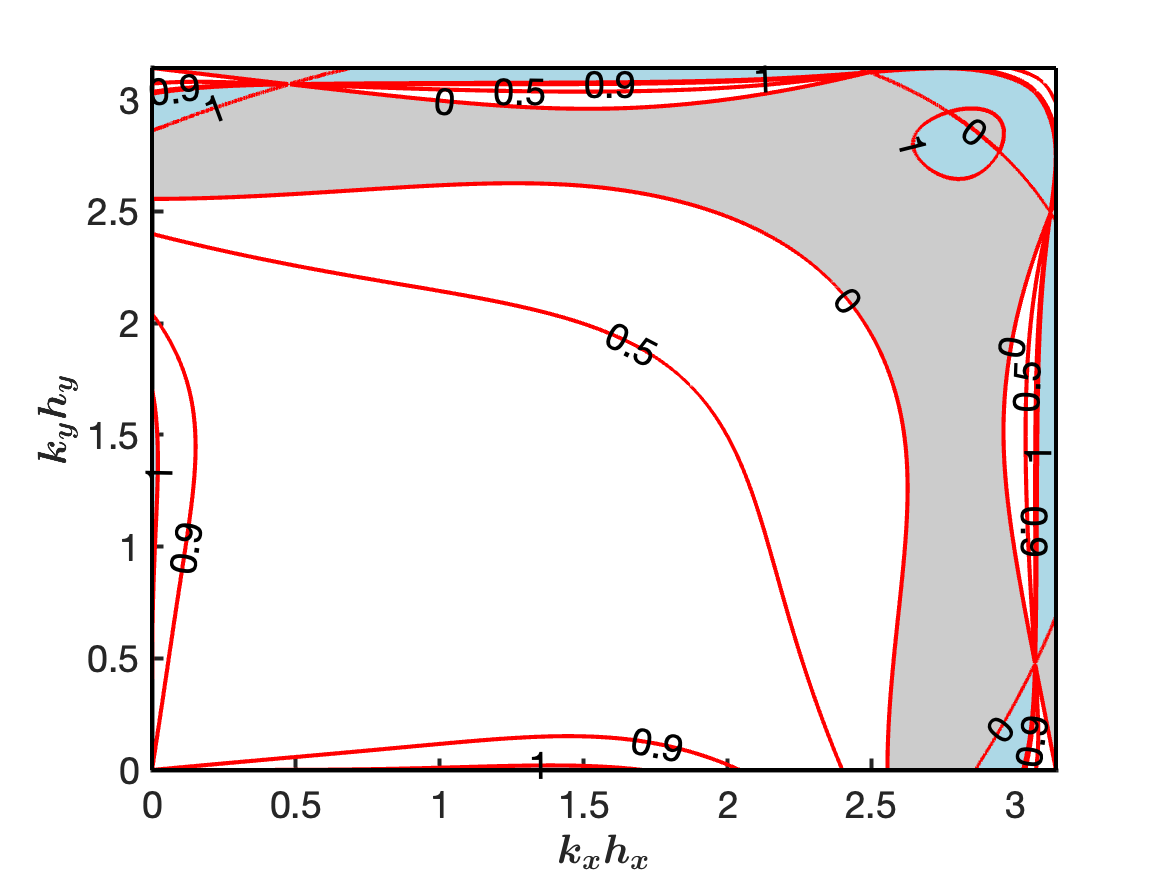}
                \end{minipage}
            }
        \end{minipage}
        &
        \hspace{\ColumnSpace} 
        \begin{minipage}{\mympwidth}
            \fbox{
                \begin{minipage}{\textwidth}
                    \makebox[\linewidth][l]{\RomanLabelSize\RomanLabelFont (v)} \\[-0.6em]
                    {\centering \TitleSize\TitleFont {$c_{\text{num}}/c_{\text{ph}}$} \par}
                    {\MinMaxSize\MinMaxFont \makebox[\linewidth]{\hspace{\ValueHSpace}Min = -1.47e+05\hfill Max = 1.24e+05\hspace{\ValueHSpace}}} \\[0.2em]
                    \includegraphics[width=\linewidth, trim=20 2 20 20, clip]{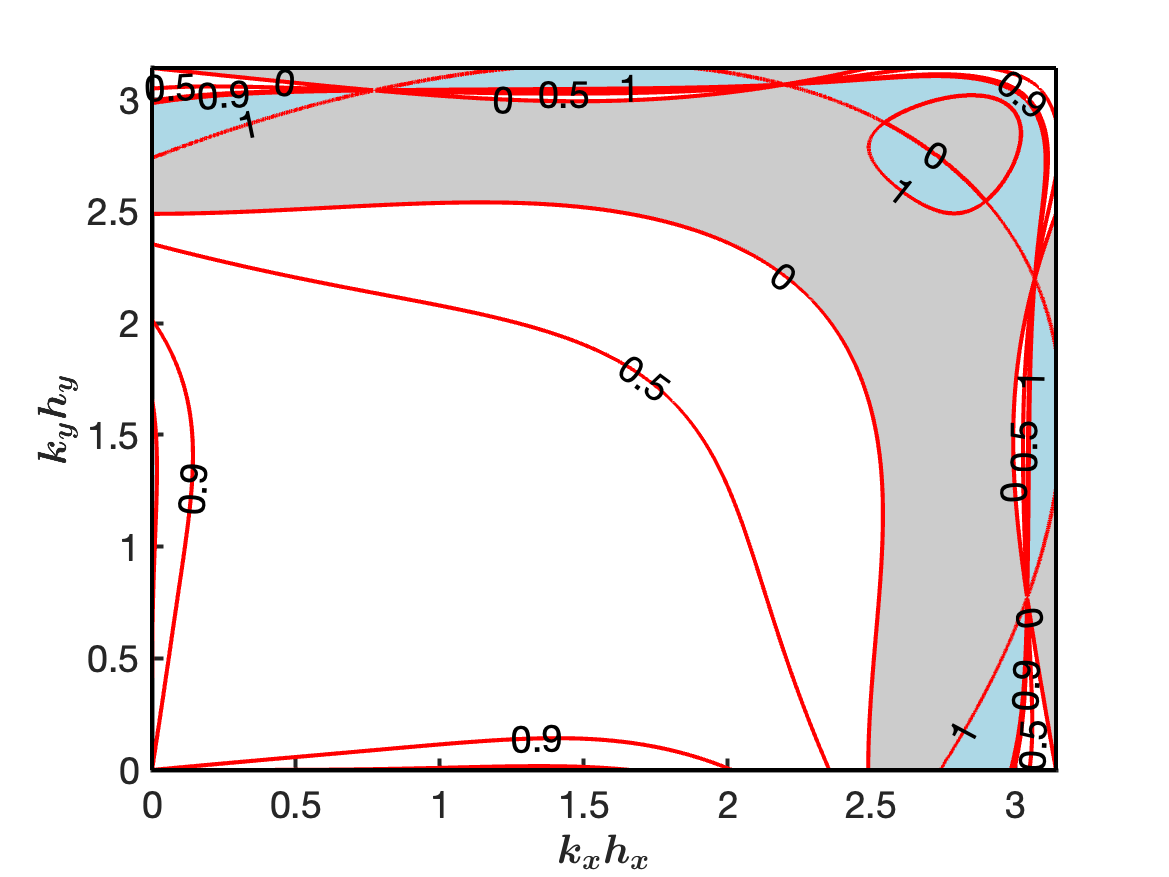}
                \end{minipage}
            }
        \end{minipage}
        &
        \hspace{\ColumnSpace} 
        \begin{minipage}{0.08\textwidth}
            \includegraphics[width=\textwidth]{Figures/3_colorbar.eps}
        \end{minipage}
        \\[1em]

        \begin{minipage}{\mympwidth}
            \fbox{
                \begin{minipage}{\textwidth}
                    \makebox[\linewidth][l]{\RomanLabelSize\RomanLabelFont (iii)} \\[-0.6em]
                    {\centering \TitleSize\TitleFont {$(v_{\text{g},x})_{\text{num}}/{v_{\text{g},x}}$} \par}
                    {\MinMaxSize\MinMaxFont \makebox[\linewidth]{\hspace{\ValueHSpace}Min = -2.90e+02\hfill Max = 1.30e+03\hspace{\ValueHSpace}}} \\[0.2em]
                    \includegraphics[width=\linewidth, trim=20 2 20 20, clip]{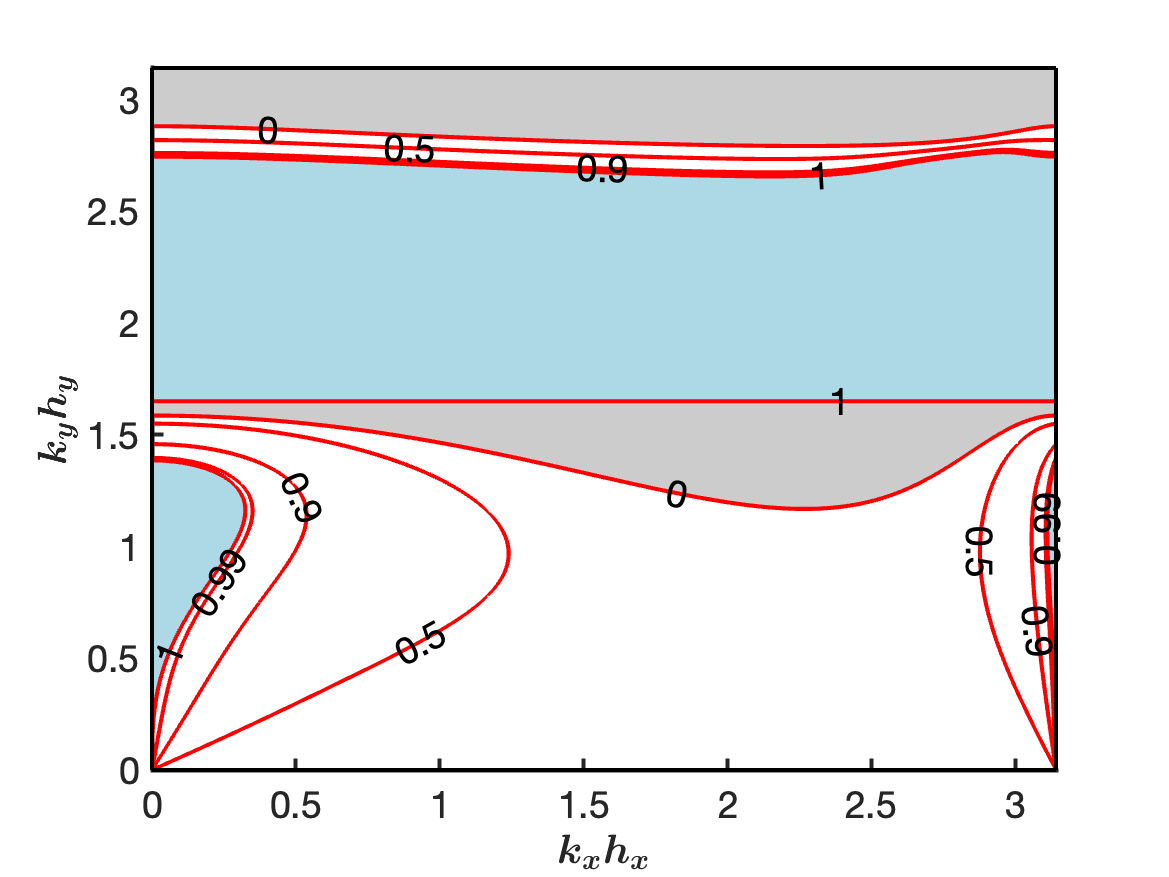}
                \end{minipage}
                }
        \end{minipage}
        &
        \hspace{\ColumnSpace} 
        \begin{minipage}{\mympwidth}
            \fbox{
                \begin{minipage}{\textwidth}
                    \makebox[\linewidth][l]{\RomanLabelSize\RomanLabelFont (vi)} \\[-0.6em]
                    {\centering \TitleSize\TitleFont {$(v_{\text{g},x})_{\text{num}}/{v_{\text{g},x}}$} \par}
                    {\MinMaxSize\MinMaxFont \makebox[\linewidth]{\hspace{\ValueHSpace}Min = -2.71e+03\hfill Max = 2.39e+02\hspace{\ValueHSpace}}} \\[0.2em]
                    \includegraphics[width=\linewidth, trim=20 2 20 20, clip]{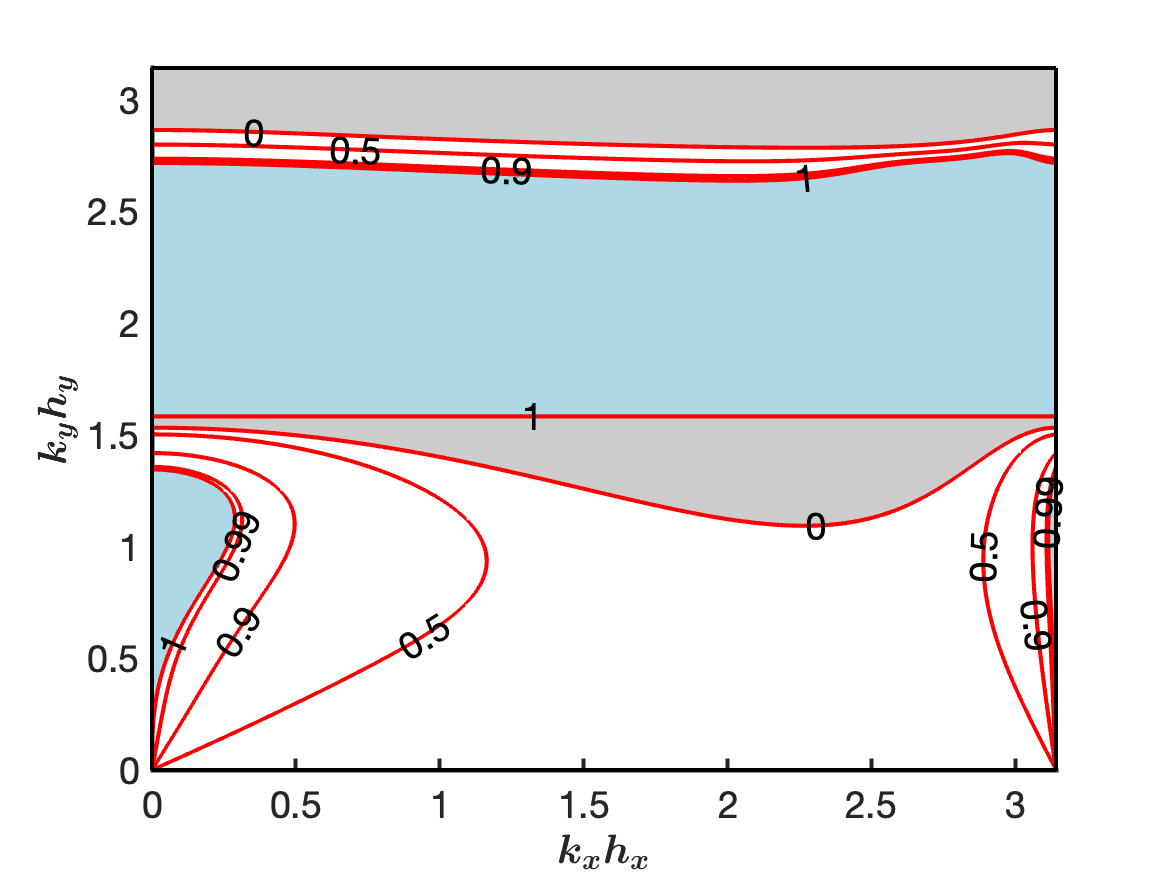}
                \end{minipage}
                }
        \end{minipage}
        &
        \hspace{\ColumnSpace} 
        \begin{minipage}{0.08\textwidth}
            \includegraphics[width=\textwidth]{Figures/3_colorbar.eps}
        \end{minipage}
    \end{tabular}
    \end{minipage}

    \caption{Contour plots of $|G_{\text{num}}|$, $\dfrac{c_{\text{num}}}{c_{\text{ph}}}$, $\dfrac{(v_{\text{g},x})_{\text{num}}}{v_{\text{g},x}}$ in the $(k_xh_x, k_yh_y)-$ plane for Eq.~\eqref{eqn:Main2D} with $N_c = 0.9$, computed using the SSPRK3-CNCS6 scheme with the indicated values of $D_{\alpha}$.}
    \label{Fig:B1}
\end{figure}

Fig.~\ref{Fig:B1} shows the numerical property charts of $|G_{\text{num}}|$, $\dfrac{c_{\text{num}}}{c_{\text{ph}}}$, and $\dfrac{(v_{\text{g},x})_{\text{num}}}{v_{\text{g},x}}$ for the SSPRK3-CNCS6, corresponding to two distinct values of the dispersion parameter $D_\alpha$. For the specific case with $N_c = 0.9$ and $D_\alpha = 0.11$, the variation of these quantities is examined over the $(k_xh_x, k_yh_y)$-plane, covering the full Nyquist range in both spatial directions. As observed in Fig.~\ref{Fig:B1}(i), the amplification factor $|G_{\text{num}}|$ stays within the stability limit (i.e., $|G_{\text{num}}|\leq 1$) throughout the range of wavenumber. However, in the case illustrated in Fig.~\ref{Fig:B1}(iv), corresponding to a slightly increased value $D_\alpha = 0.12$, the maximum value of $|G_{\text{num}}|$ marginally exceeds one. This indicates the onset of instability in the numerical scheme for this parameter setting. Accordingly, $D_{\alpha,\mathrm{cr}} = 0.12$ is identified as the critical dispersion parameter beyond which the scheme may become unstable.

\par
Additional insights are revealed by analysing phase and group velocity ratios, where regions of non-physical behavior are detected. In particular, the contour plots of $\dfrac{c_{\text{num}}}{c_{\text{ph}}}$ shown in Fig.~\ref{Fig:B1}(ii) exhibit zones beyond $(k_xh_x, k_yh_y) \approx (2.5,2.5)$ where the numerical phase velocity becomes negative. This implies a spurious wave phase propagation. Likewise, the group velocity plots in Fig.~\ref{Fig:B1}(iii) highlight two distinct regions, one near $k_xh_x = 0$ and another for $k_yh_y > 1.5$, where the numerical group velocity ratio $\dfrac{(v_{\text{g},x})_{\text{num}}}{v_{\text{g},x}}$ exceeds unity, suggesting artificially accelerated energy transport. Furthermore, two additional regions are observed, one in the range $1.1 < k_yh_y < 1.6$ and another beyond $k_yh_y \approx 2.7$, where the group velocity ratio becomes negative \bigg($\dfrac{(v_{\text{g},x})_{\text{num}}}{v_{\text{g},x}} < 0$\bigg), indicating energy propagation in the opposite direction, which is physically inconsistent.

\subsubsection{SSPRK3-CNCS8}
\begin{figure}[htbp!]
    \centering
    \begin{minipage}{\textwidth}
    \hspace*{0cm} 
    \begin{tabular}{ccc}
        \multicolumn{3}{c}{\makebox[0pt]{\hspace{0em}\textcolor{red}{\textbf{SSPRK3-CNCS8}}}} \\
        \multicolumn{1}{c}{\textbf{$\bm{D}_{\bm{\alpha}} = \bm{0.11}$}} & \textbf{$\bm{D}_{\bm{\alpha,\mathrm{cr}}} = \bm{0.12}$} \\[0.5em]

        \begin{minipage}{\mympwidth}
            \fbox{
               \begin{minipage}{\textwidth}
                    \makebox[\linewidth][l]{\RomanLabelSize\RomanLabelFont (i)} \\[-0.6em]
                    {\centering \TitleSize\TitleFont {$|G_{\text{num}}|$} \par}
                    {\MinMaxSize\MinMaxFont \makebox[\linewidth]{\hspace{\ValueHSpace}Min = 0.94\hfill Max = 1.00\hspace{\ValueHSpace}}} \\[0.2em]
                    \includegraphics[width=\linewidth, trim=20 2 20 20, clip]{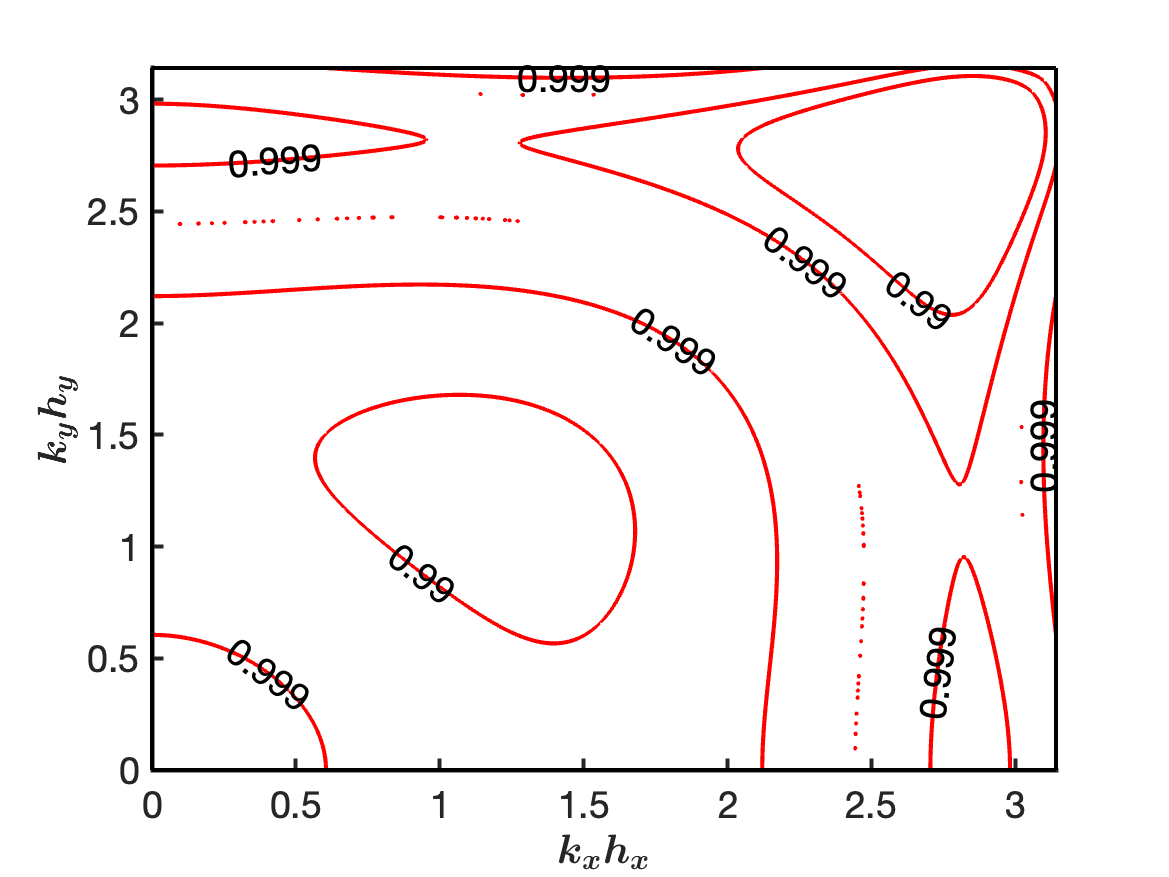}
               \end{minipage}
            }
        \end{minipage}
        &
        \hspace{\ColumnSpace}
        \begin{minipage}{\mympwidth}
            \fbox{
                \begin{minipage}{\textwidth}
                    \makebox[\linewidth][l]{\RomanLabelSize\RomanLabelFont (iv)} \\[-0.6em]
                    {\centering \TitleSize\TitleFont {$|G_{\text{num}}|$} \par}
                    {\MinMaxSize\MinMaxFont \makebox[\linewidth]{\hspace{\ValueHSpace}Min = 0.94\hfill Max = 1.11\hspace{\ValueHSpace}}} \\[0.2em]
                    \includegraphics[width=\linewidth, trim=20 2 20 20, clip]{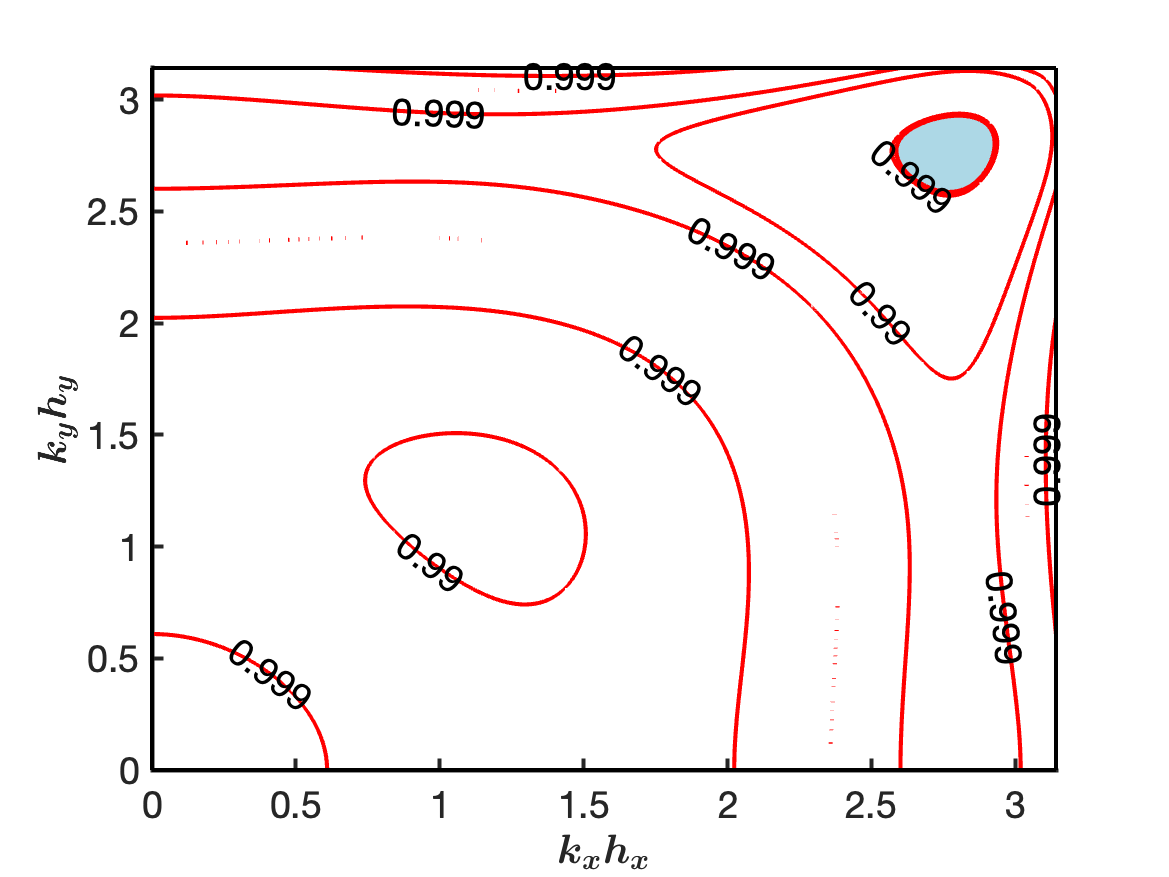}
                \end{minipage}
            }
        \end{minipage}
        &
        \hspace{\ColumnSpace}
        \begin{minipage}{0.08\textwidth}
            \includegraphics[width=\textwidth]{Figures/2_colorbar.eps}
        \end{minipage}
        \\[1em]

        \begin{minipage}{\mympwidth}
            \fbox{
                \begin{minipage}{\textwidth}
                    \makebox[\linewidth][l]{\RomanLabelSize\RomanLabelFont (ii)} \\[-0.6em]
                    {\centering \TitleSize\TitleFont {$c_{\text{num}}/c_{\text{ph}}$} \par}
                    {\MinMaxSize\MinMaxFont \makebox[\linewidth]{\hspace{\ValueHSpace}Min = -8.01e+04\hfill Max = 2.56e+05\hspace{\ValueHSpace}}} \\[0.2em]
                    \includegraphics[width=\linewidth, trim=20 2 20 20, clip]{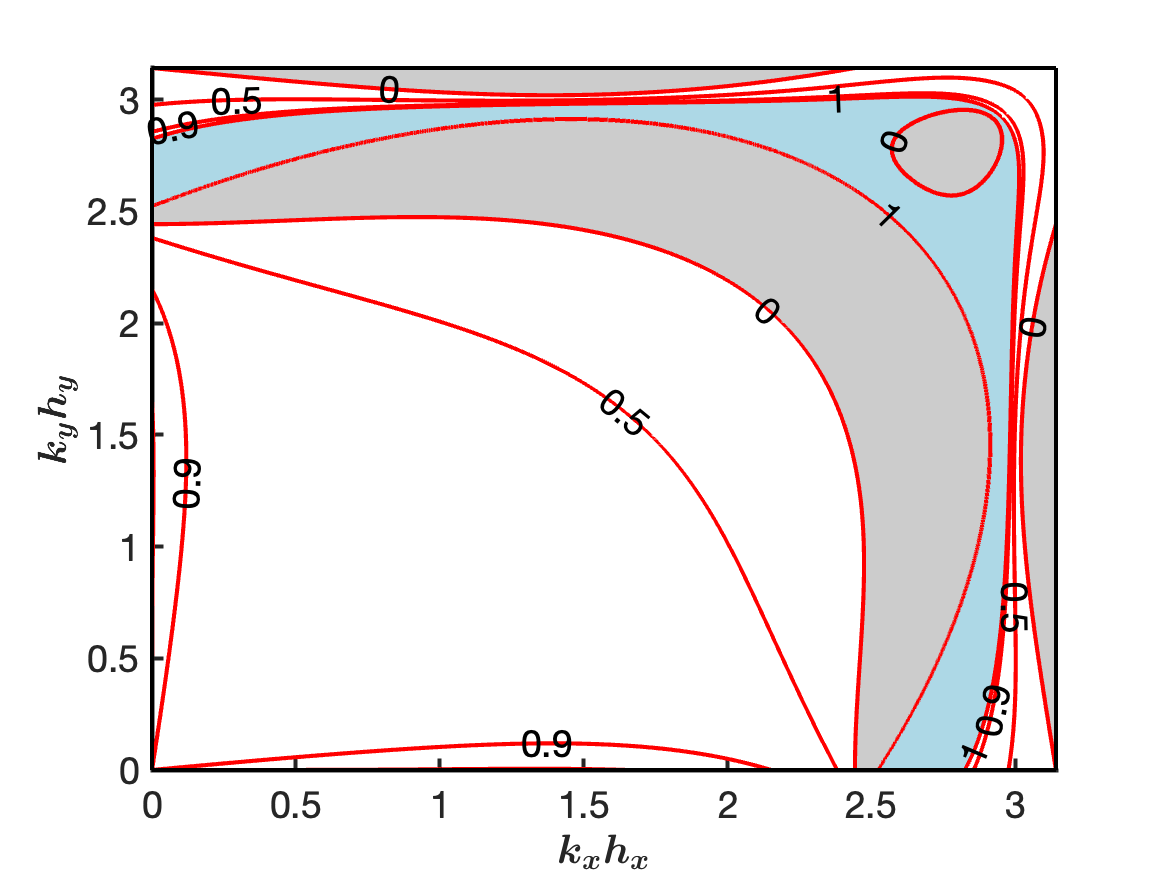}
                \end{minipage}
            }
        \end{minipage}
        &
        \hspace{\ColumnSpace}
        \begin{minipage}{\mympwidth}
            \fbox{
                \begin{minipage}{\textwidth}
                    \makebox[\linewidth][l]{\RomanLabelSize\RomanLabelFont (v)} \\[-0.6em]
                    {\centering \TitleSize\TitleFont {$c_{\text{num}}/c_{\text{ph}}$} \par}
                    {\MinMaxSize\MinMaxFont \makebox[\linewidth]{\hspace{\ValueHSpace}Min = -9.42e+04\hfill Max = 7.24e+04\hspace{\ValueHSpace}}} \\[0.2em]
                    \includegraphics[width=\linewidth, trim=20 2 20 20, clip]{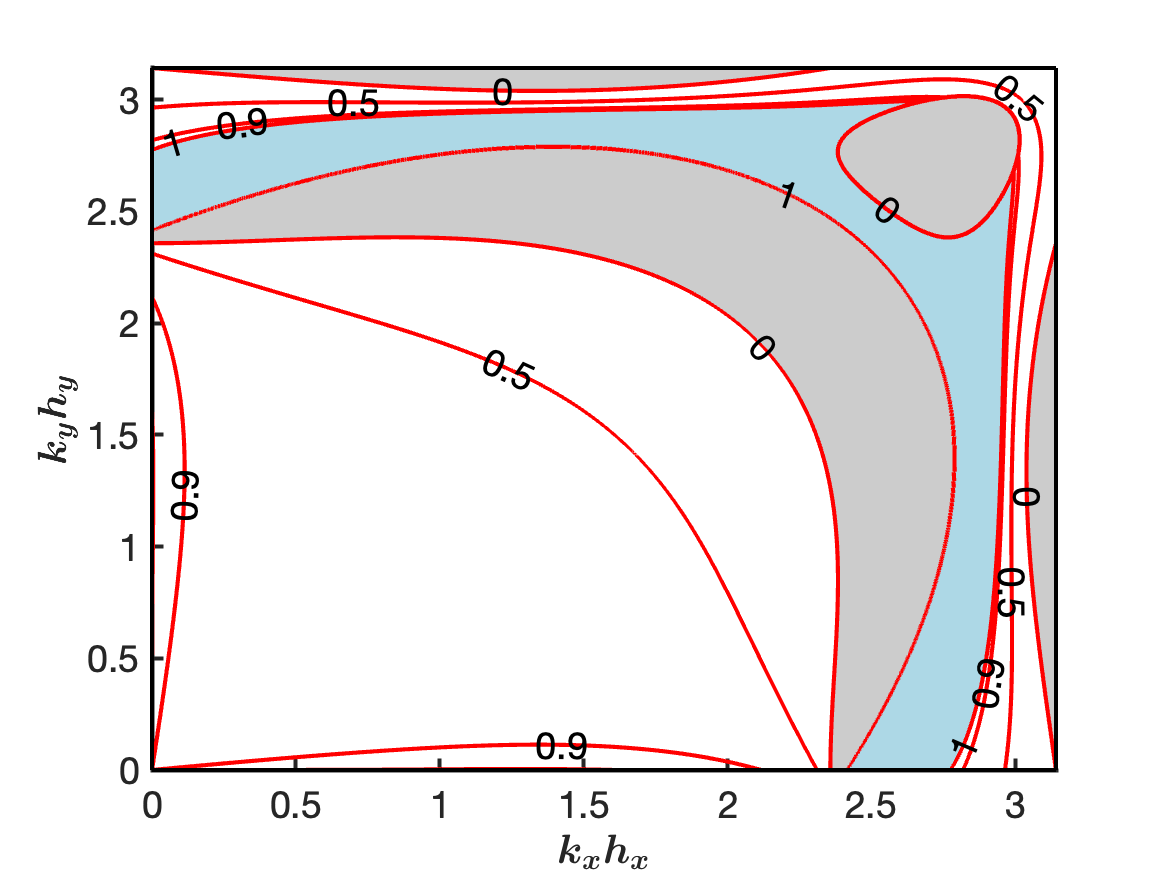}
                \end{minipage}
            }
        \end{minipage}
        &
        \hspace{\ColumnSpace}
        \begin{minipage}{0.08\textwidth}
            \includegraphics[width=\textwidth]{Figures/3_colorbar.eps}
        \end{minipage}
        \\[1em]

        \begin{minipage}{\mympwidth}
            \fbox{
                \begin{minipage}{\textwidth}
                    \makebox[\linewidth][l]{\RomanLabelSize\RomanLabelFont (iii)} \\[-0.6em]
                    {\centering \TitleSize\TitleFont {$(v_{\text{g},x})_{\text{num}}/{v_{\text{g},x}}$} \par}
                    {\MinMaxSize\MinMaxFont \makebox[\linewidth]{\hspace{\ValueHSpace}Min = -2.78e+02\hfill Max = 7.70e+02\hspace{\ValueHSpace}}} \\[0.2em]
                    \includegraphics[width=\linewidth, trim=20 2 20 20, clip]{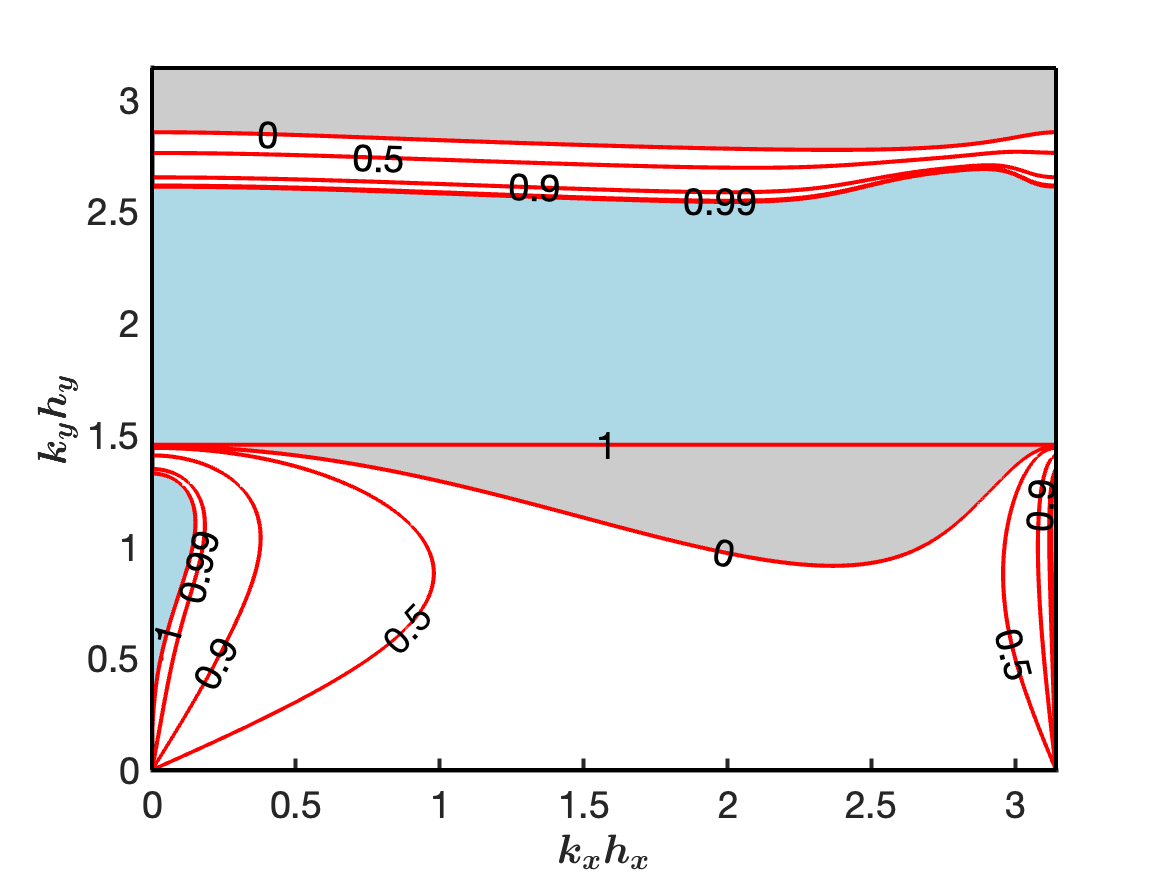}
                \end{minipage}
                }
        \end{minipage}
        &
        \hspace{\ColumnSpace}
        \begin{minipage}{\mympwidth}
            \fbox{
                \begin{minipage}{\textwidth}
                    \makebox[\linewidth][l]{\RomanLabelSize\RomanLabelFont (vi)} \\[-0.6em]
                    {\centering \TitleSize\TitleFont {$(v_{\text{g},x})_{\text{num}}/{v_{\text{g},x}}$} \par}
                    {\MinMaxSize\MinMaxFont \makebox[\linewidth]{\hspace{\ValueHSpace}Min = -7.23e+02\hfill Max = 2.82e+02\hspace{\ValueHSpace}}} \\[0.2em]
                    \includegraphics[width=\linewidth, trim=20 2 20 20, clip]{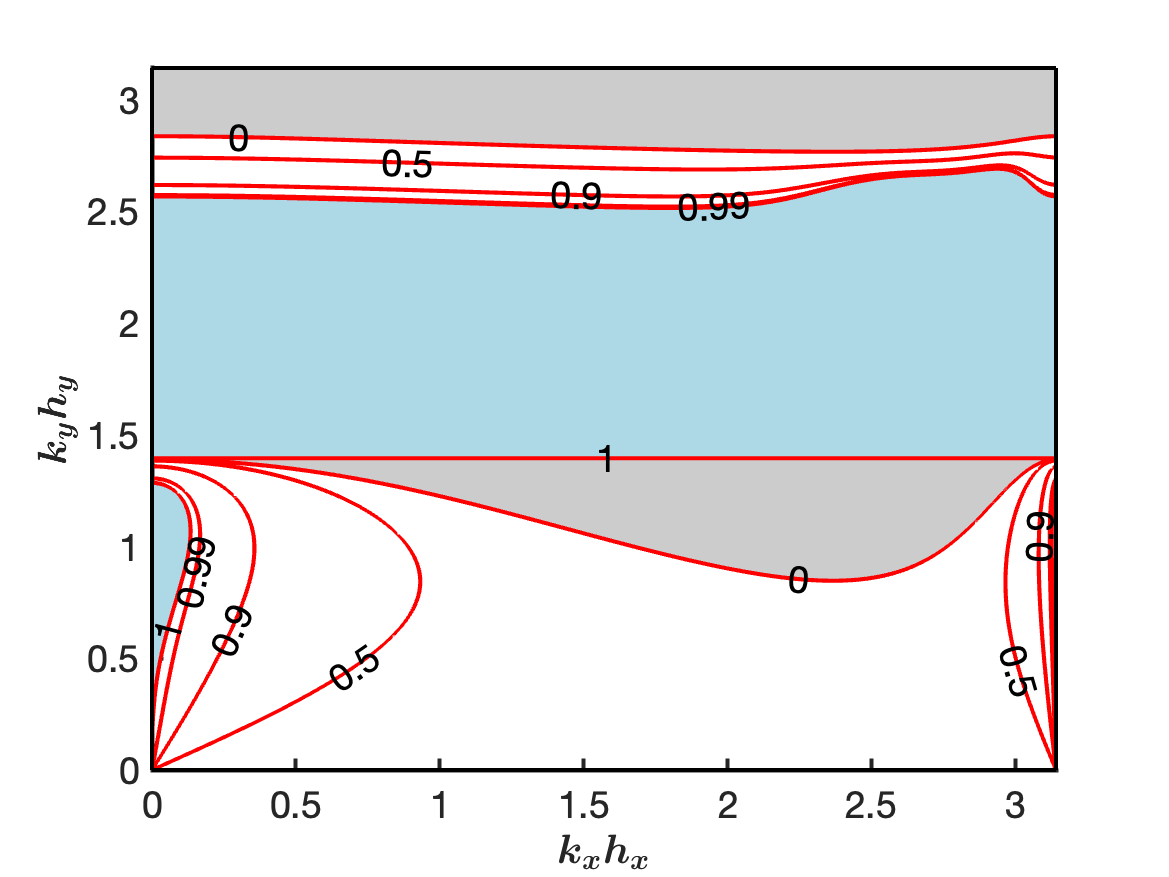}
                \end{minipage}
                }
        \end{minipage}
        &
        \hspace{\ColumnSpace}
        \begin{minipage}{0.08\textwidth}
            \includegraphics[width=\textwidth]{Figures/3_colorbar.eps}
        \end{minipage}
    \end{tabular}
    \end{minipage}
    \caption{Contour plots of $|G_{\text{num}}|$, $\dfrac{c_{\text{num}}}{c_{\text{ph}}}$, $\dfrac{(v_{\text{g},x})_{\text{num}}}{v_{\text{g},x}}$ in the $(k_xh_x, k_yh_y)-$ plane for Eq.~\eqref{eqn:Main2D} with $N_c = 0.7$, computed using the SSPRK3-CNCS8 scheme with the indicated values of $D_{\alpha}$.}
    \label{Fig:B2}
\end{figure}

Next, we have performed a similar analysis using the eighth-order SSPRK3-CNCS8. Fig.~\ref{Fig:B2} presents the behavior of $|G_{\text{num}}|$, $\dfrac{c_{\text{num}}}{c_{\text{ph}}}$, and $\dfrac{(v_{\text{g},x})_{\text{num}}}{v_{\text{g},x}}$ for two different values of $D_\alpha$. For $N_c = 0.7$ and $D_\alpha = 0.11$, the plots across the $(k_xh_x, k_yh_y)$ domain confirm that the amplification factor remains within stable bounds. However, when $D_\alpha$ is slightly increased to $0.12$ (see Fig.~\ref{Fig:B2}(iv)), $|G_{\text{num}}|$ slightly exceeds unity, hinting at the onset of instability. Therefore, $D_{\alpha,\mathrm{cr}} = 0.12$ is identified as the threshold value for this setup.\par
Further examination of the velocity ratios reveals artifacts that indicate unphysical wave behavior. For example, negative phase speeds appear beyond $kh \approx 2.5$ in Fig.~\ref{Fig:B2}(ii), and the group velocity plots (Fig.~\ref{Fig:B2}(iii)) show regions of anomalous energy transport, including zones where energy travels faster than expected or in the wrong direction. These artifacts, also observed in the sixth-order CNCS6, highlight that higher-order accuracy alone does not entirely eliminate spurious numerical effects such as reversed \( q \)-waves unless the dispersion parameters are carefully tuned. However, in the case of CNCS8, a notable reduction in the extent of the negative group velocity region is observed. Specifically, the area of the arc-shaped zone with \( \dfrac{(v_{\text{g},x})_{\text{num}}}{v_{\text{g},x}} < 0 \) between \( kh \in [0.5, 1.5] \) is significantly smaller compared to CNCS6, as evident in the group velocity contours. This reduction indicates improved directional accuracy and energy transport behavior in mid-frequency ranges, suggesting that the eighth-order discretization offers better suppression of nonphysical wave reversal, particularly in the moderately dispersive regime. Nevertheless, the persistence of localized instability and anomalous transport at higher wavenumbers reinforces the need for precise calibration of scheme parameters to ensure both stability and physical fidelity in multidimensional simulations.
 
\subsubsection{SSPRK3-CCS8}
\begin{figure}[htbp!]
    \centering
    \begin{minipage}{\textwidth}
    \hspace*{0cm} 
    \begin{tabular}{ccc}
        \multicolumn{3}{c}{\makebox[0pt]{\hspace{0em}\textcolor{red}{\textbf{SSPRK3-CCS8}}}} \\
        \multicolumn{1}{c}{\textbf{$\bm{D}_{\bm{\alpha}} = \bm{0.011}$}} & \textbf{$\bm{D}_{\bm{\alpha,\mathrm{cr}}} = \bm{0.012}$} \\[0.5em]

        \begin{minipage}{\mympwidth}
            \fbox{
               \begin{minipage}{\textwidth}
                    \makebox[\linewidth][l]{\RomanLabelSize\RomanLabelFont (i)} \\[-0.6em]
                    {\centering \TitleSize\TitleFont {$|G_{\text{num}}|$} \par}
                    {\MinMaxSize\MinMaxFont \makebox[\linewidth]{\hspace{\ValueHSpace}Min = 0.94\hfill Max = 1.00\hspace{\ValueHSpace}}} \\[0.2em]
                    \includegraphics[width=\linewidth, trim=20 2 20 20, clip]{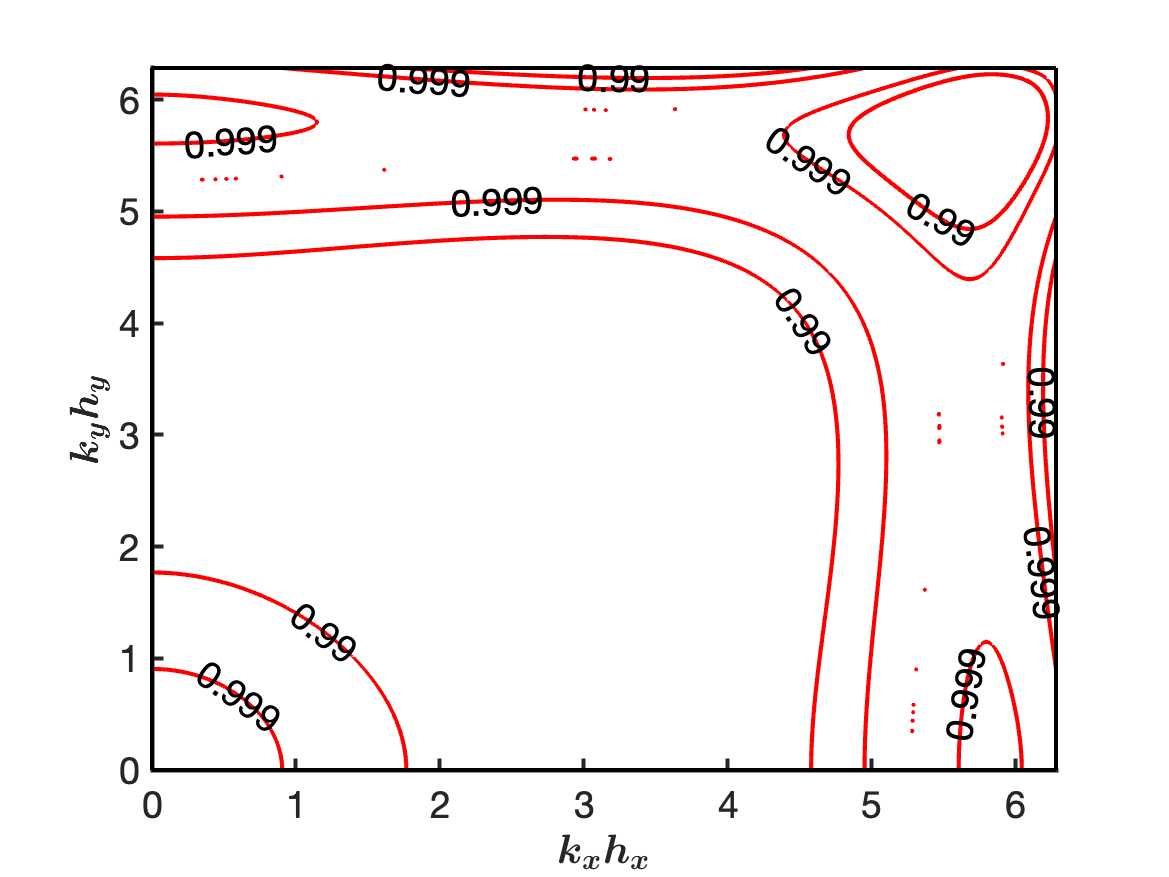}
               \end{minipage}
            }
        \end{minipage}
        &
        \hspace{\ColumnSpace}
        \begin{minipage}{\mympwidth}
            \fbox{
                \begin{minipage}{\textwidth}
                    \makebox[\linewidth][l]{\RomanLabelSize\RomanLabelFont (iv)} \\[-0.6em]
                    {\centering \TitleSize\TitleFont {$|G_{\text{num}}|$} \par}
                    {\MinMaxSize\MinMaxFont \makebox[\linewidth]{\hspace{\ValueHSpace}Min = 0.94\hfill Max = 1.07\hspace{\ValueHSpace}}} \\[0.2em]
                    \includegraphics[width=\linewidth, trim=20 2 20 20, clip]{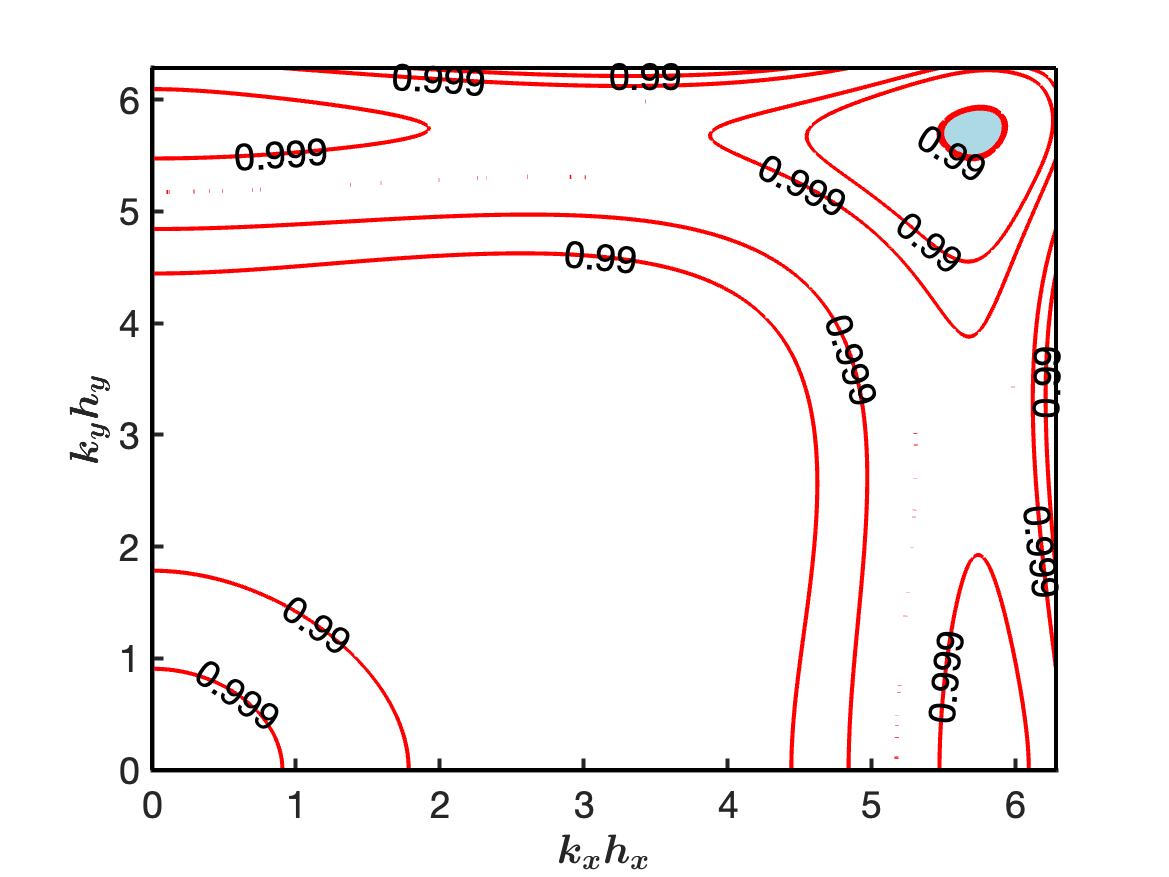}
                \end{minipage}
            }
        \end{minipage}
        &
        \hspace{\ColumnSpace}
        \begin{minipage}{0.08\textwidth}
            \includegraphics[width=\textwidth]{Figures/2_colorbar.eps}
        \end{minipage}
        \\[1em]

        \begin{minipage}{\mympwidth}
            \fbox{
                \begin{minipage}{\textwidth}
                    \makebox[\linewidth][l]{\RomanLabelSize\RomanLabelFont (ii)} \\[-0.6em]
                    {\centering \TitleSize\TitleFont {$c_{\text{num}}/c_{\text{ph}}$} \par}
                    {\MinMaxSize\MinMaxFont \makebox[\linewidth]{\hspace{\ValueHSpace}Min = -2.07\hfill Max = 2.06\hspace{\ValueHSpace}}} \\[0.2em]
                    \includegraphics[width=\linewidth, trim=20 2 20 20, clip]{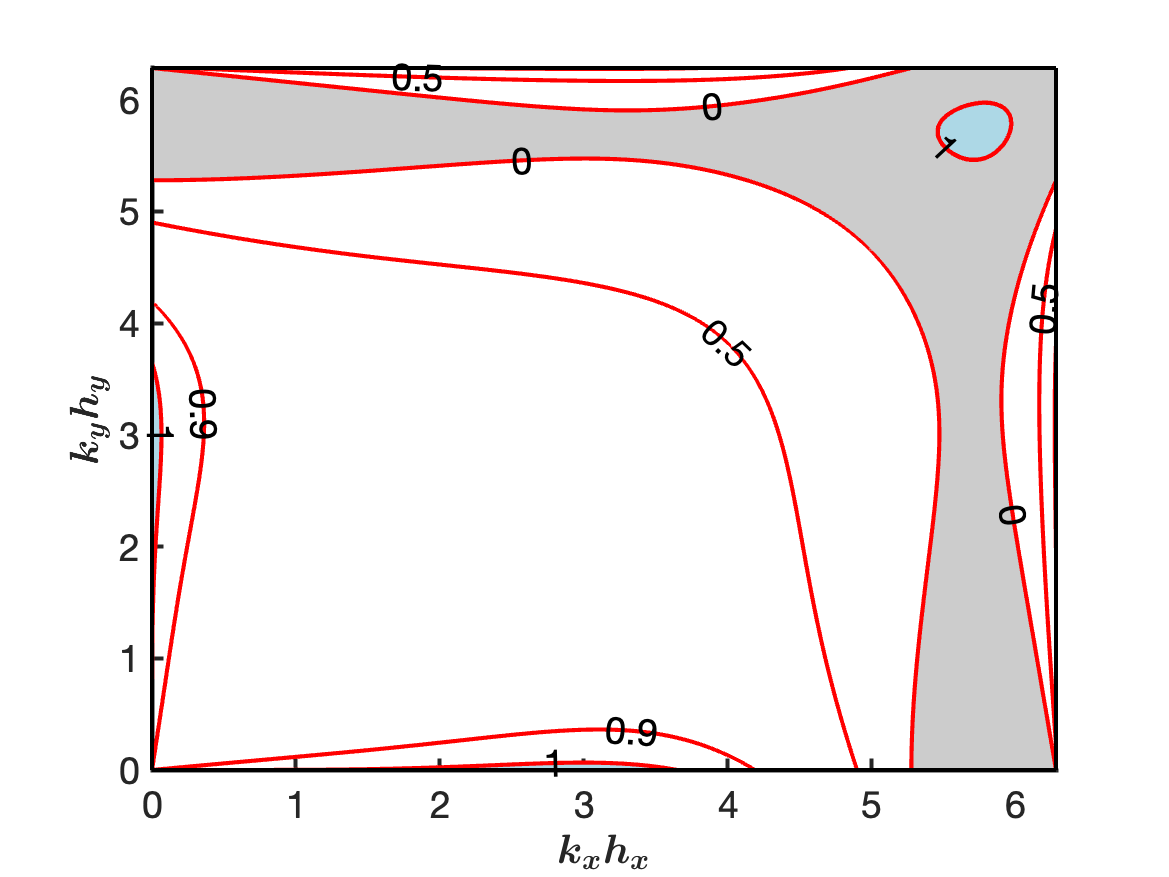}
                \end{minipage}
            }
        \end{minipage}
        &
        \hspace{\ColumnSpace}
        \begin{minipage}{\mympwidth}
            \fbox{
                \begin{minipage}{\textwidth}
                    \makebox[\linewidth][l]{\RomanLabelSize\RomanLabelFont (v)} \\[-0.6em]
                    {\centering \TitleSize\TitleFont {$c_{\text{num}}/c_{\text{ph}}$} \par}
                    {\MinMaxSize\MinMaxFont \makebox[\linewidth]{\hspace{\ValueHSpace}Min = -4.16e+04\hfill Max = 1.42e+04\hspace{\ValueHSpace}}} \\[0.2em]
                    \includegraphics[width=\linewidth, trim=20 2 20 20, clip]{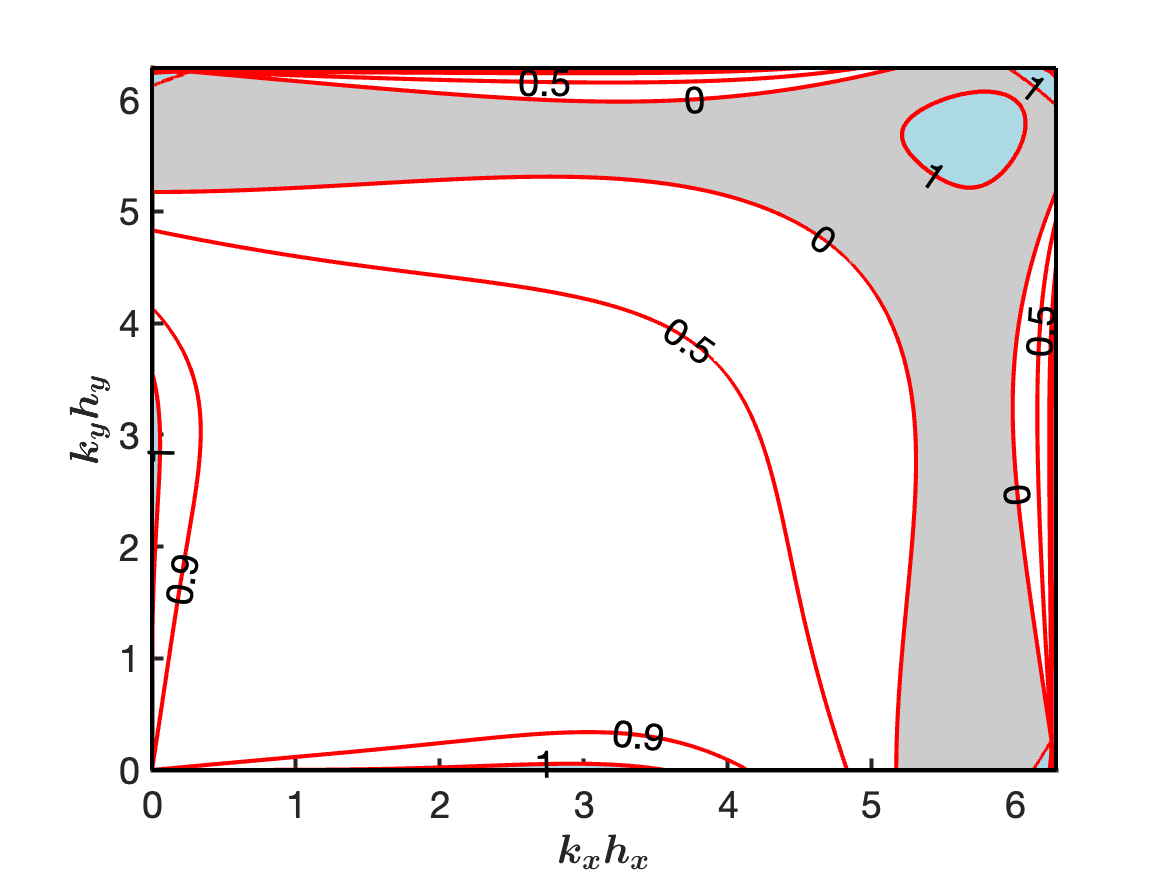}
                \end{minipage}
            }
        \end{minipage}
        &
        \hspace{\ColumnSpace}
        \begin{minipage}{0.08\textwidth}
            \includegraphics[width=\textwidth]{Figures/3_colorbar.eps}
        \end{minipage}
        \\[1em]

        \begin{minipage}{\mympwidth}
            \fbox{
                \begin{minipage}{\textwidth}
                    \makebox[\linewidth][l]{\RomanLabelSize\RomanLabelFont (iii)} \\[-0.6em]
                    {\centering \TitleSize\TitleFont {$(v_{\text{g,x}})_{\text{num}}/{v_{\text{g,x}}}$} \par}
                    {\MinMaxSize\MinMaxFont \makebox[\linewidth]{\hspace{\ValueHSpace}Min = -4.03e+02\hfill Max = 2.28e+03\hspace{\ValueHSpace}}} \\[0.2em]
                    \includegraphics[width=\linewidth, trim=20 2 20 20, clip]{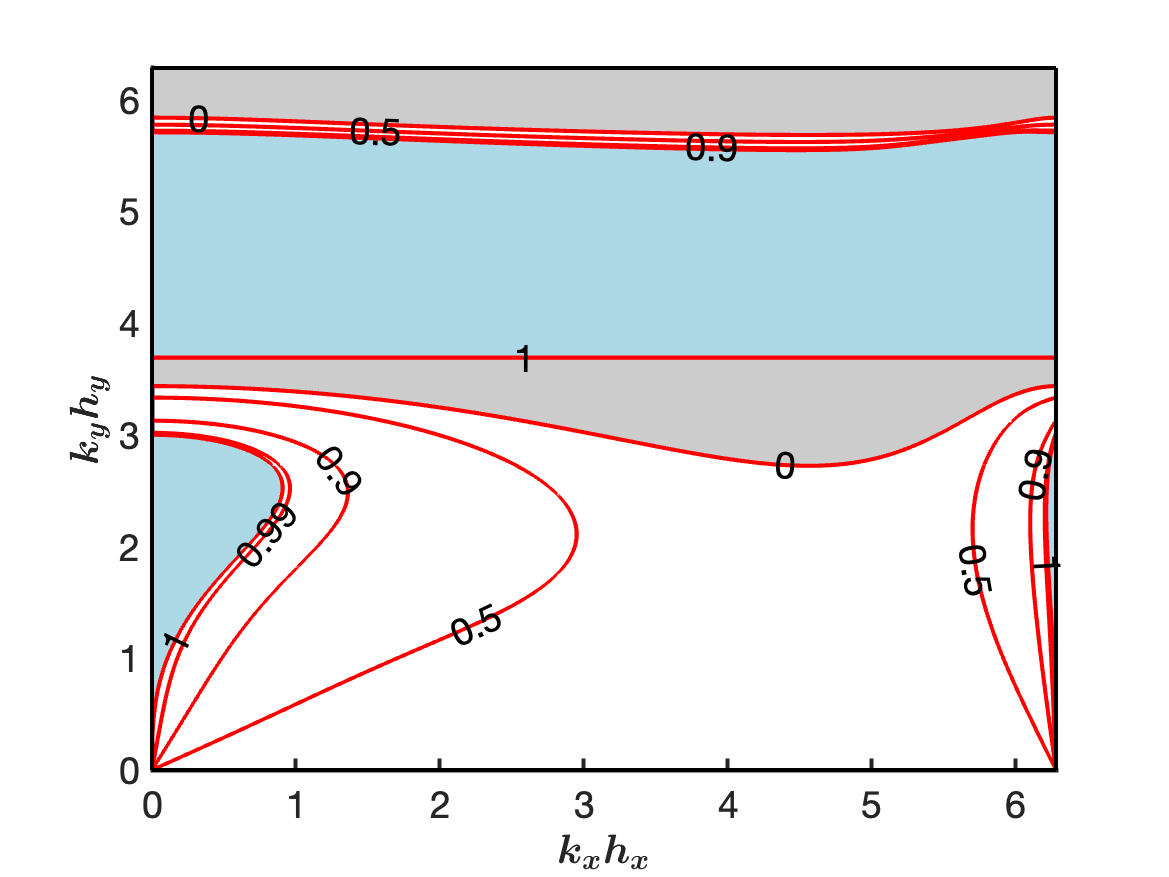}
                \end{minipage}
                }
        \end{minipage}
        &
        \hspace{\ColumnSpace}
        \begin{minipage}{\mympwidth}
            \fbox{
                \begin{minipage}{\textwidth}
                    \makebox[\linewidth][l]{\RomanLabelSize\RomanLabelFont (vi)} \\[-0.6em]
                    {\centering \TitleSize\TitleFont {$(v_{\text{g,x}})_{\text{num}}/{v_{\text{g,x}}}$} \par}
                    {\MinMaxSize\MinMaxFont \makebox[\linewidth]{\hspace{\ValueHSpace}Min = -3.51e+02\hfill Max = 1.75e+03\hspace{\ValueHSpace}}} \\[0.2em]
                    \includegraphics[width=\linewidth, trim=20 2 20 20, clip]{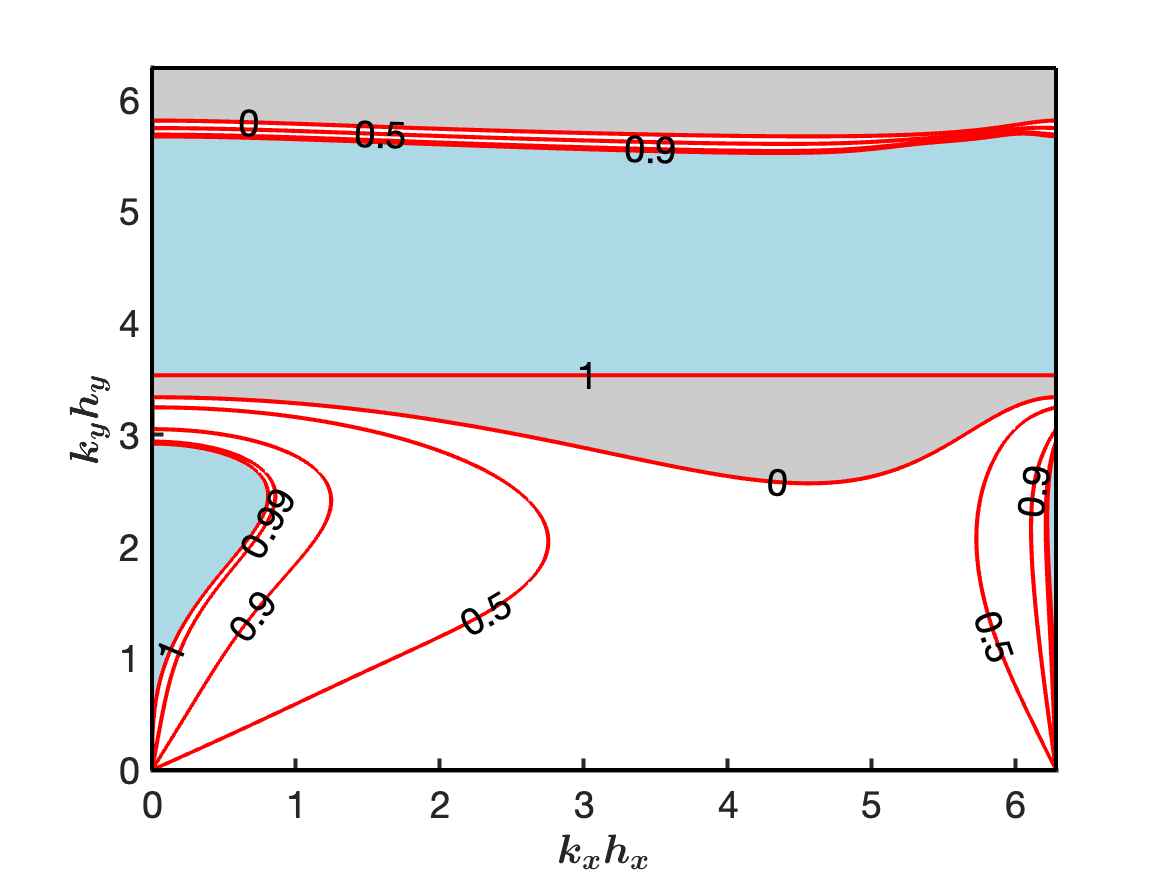}
                \end{minipage}
                }
        \end{minipage}
        &
        \hspace{\ColumnSpace}
        \begin{minipage}{0.08\textwidth}
            \includegraphics[width=\textwidth]{Figures/3_colorbar.eps}
        \end{minipage}
    \end{tabular}
    \end{minipage}
    \caption{Contour plots of $|G_{\text{num}}|$, $\dfrac{c_{\text{num}}}{c_{\text{ph}}}$, $\dfrac{(v_{\text{g},x})_{\text{num}}}{v_{\text{g},x}}$ in the $(k_xh_x, k_yh_y)-$ plane for Eq.~\eqref{eqn:Main2D} with $N_c = 0.45$, computed using the SSPRK3-CCS8 scheme with the indicated values of $D_{\alpha}$.}
    \label{Fig:B3}
\end{figure}

Fig.~\ref{Fig:B3} presents the contour plots of the numerical amplification factor \( |G_{\text{num}}| \), the normalized phase speed \( \dfrac{c_{\text{num}}}{c_{\text{ph}}} \), and the normalized group velocity in the \( x \)-direction \( \dfrac{(v_{\text{g},x})_{\text{num}}}{v_{\text{g},x}} \) for the SSPRK3-CCS8. As in previous cases, these quantities are evaluated for two values of the dispersion parameter, \( D_\alpha = 0.011 \) and \( D_\alpha = 0.012 \), at a fixed Courant number \( N_c = 0.45 \). The spectral analysis is performed over the \( (k_x h_x, k_y h_y) \)-plane, with wavenumbers extending up to \( kh = 2\pi \). This extended spectral domain is appropriate for the CCS8 scheme because it is a cell-centered compact finite difference method, which allows resolution of the full \( kh \in [0, 2\pi] \) range-twice that of node-centered schemes-due to the staggered placement of degrees of freedom. For \( D_\alpha = 0.011 \), the amplification factor remains bounded by unity throughout the entire spectral domain, as shown in Fig.~\ref{Fig:B3}(i), ensuring numerical stability. However, when the dispersion parameter is increased to \( D_\alpha = 0.012 \) (Fig.~\ref{Fig:B3}(iv)), the amplification factor slightly exceeds unity,  indicating the emergence of an unstable region. This behavior is consistent with the instability observed in previously reported schemes for the CCS8 configuration under similar conditions.

\par
Fig.~\ref{Fig:B3}(ii) shows the behavior of the numerical phase velocity. Regions near $ (k_x h_x, k_y h_y) \approx (5.2, 5.2) $ exhibit a sign reversal, where $ \dfrac{c_{\text{num}}}{c_{\text{ph}}} < 0 $, indicating backward wave propagation. The minimum and maximum values of this ratio are $ -2.07 $ and $ 2.06 $, respectively, substantially smaller in magnitude compared to those produced by the SSPRK3-CNCS8, where the extremes reach $ -8.01 \times 10^4 $ and $ 2.56 \times 10^5 $. These differences highlight the reduced tendency of SSPRK3-CCS8 to generate non-physical oscillations and \( q \)-waves.

\par
In Fig.~\ref{Fig:B3}(iii), the numerical group velocity ratio in the $ x $-direction is examined. Values exceeding unity occur near $ k_x h_x \approx 0 $ and for $ k_y h_y > 3.65 $, indicating artificially enhanced energy transport. Additionally, negative group velocity ratios are observed in the bands $ 2.6 < k_y h_y < 3.65 $ and beyond $ k_y h_y \approx 5.6 $, pointing to energy transport in the reverse direction. Although these features share similarities with those observed in SSPRK3-CNCS8, they extend over a larger portion of the spectral domain in the case of SSPRK3-CCS8, suggesting improved spectral accuracy and dispersion control.

\section{Numerical Experiments}\label{Sec:NE}
This section provides a detailed presentation of numerical results for a variety of test problems in both one and two space dimensions. The simulations are performed using the high-order compact numerical schemes CNCS6, CNCS8, and CCS8 for spatial discretization, coupled with the SSPRK3 method for temporal integration. All computations are performed with spatially periodic boundary conditions to ensure consistency with numerical analysis requirements. The time step $\Delta t$ is selected based on a CFL-type condition that ensures numerical stability. For one-dimensional equations of the form $u_t + g(u)_x + f(u)_{xxx} = 0,$ the time step is computed as:
\begin{equation} \label{eqn:timestep1D}
    \Delta t = \dfrac{\text{CFL}}{ \left( \dfrac{\max |g'(u)|}{h_x} + \dfrac{\max |f'(u)|}{h_x^3} \right) },
\end{equation}
where $\text{CFL} = \min (N_c , D_\alpha)$. For the two-dimensional case involving a convection–dispersion equation of the form
$u_t + g_1(u)_x + g_2(u)_y + f_1(u)_{xxx} + f_2(u)_{yyy} = 0,$
the time step is chosen as:
\begin{equation} \label{eqn:timestep2D}
    \Delta t = \dfrac{\text{CFL}}{ \left( \dfrac{\max |g_1'(u)|}{h_x} + \dfrac{\max |g_2'(u)|}{h_y} + \dfrac{\max |f_1'(u)|}{h_x^3} + \dfrac{\max |f_2'(u)|}{h_y^3} \right) }.
\end{equation}
These formulations ensure that both the convection and dispersive terms are properly accounted for in determining a stable time step, thereby enabling accurate resolution of wave propagation and dispersive effects throughout the computational domain. In the numerical experiments presented below, the CFL number is set to $\text{CFL} = 0.11$ for the CNCS  and $\text{CFL} = 0.011$ for the CCS. These values are chosen based on stability insights from the GSA analysis and are further validated through examples~\ref{example:1} and ~\ref{2D_example:1}.
\begin{example} \label{example:1}
\normalfont
Consider the one-dimensional linear homogeneous convection-dispersive equation~\cite{wazwaz2003analytic},
\begin{equation} \label{Eqn:Eg1}
    u_t + 2 u_x + \dfrac{1}{k^2} u_{xxx} = 0, \quad (x,t)\in [0, 2\pi]\times [0, T],
\end{equation}
subject to the initial condition $ u(x,0) = \sin(kx) $. The exact solution is given by $ u(x,t) = \sin(k(x - t)) $. Fig.~\ref{Fig:Eg1a} illustrates the temporal evolution of the maximum absolute error, defined as the difference between the exact and numerical solutions, up to $ t = 0.5 $ with spatial resolution $ N = 100 $. Specifically, Figs.~\ref{Fig:Eg1a}(a) and~\ref{Fig:Eg1a}(b) display the results obtained using the CNCS  for the maximum stable dispersion number $D_{\alpha} = 0.11$ and the critical dispersion number $D_{\alpha,\mathrm{cr}} = 0.12$, respectively, for both sixth- and eighth-order accurate schemes. Fig.~\ref{Fig:Eg1a}(c) shows the corresponding results using the eighth-order CCS  with $D_{\alpha} = 0.011$ and $D_{\alpha,\mathrm{cr}} = 0.012$. The plots clearly demonstrate that the solution becomes unstable at the critical dispersion number, which is in agreement with the predictions made by the GSA.\par
In addition, Table~\ref{Table:Eg1} reports the $ L_{\infty} $ error and the corresponding order of convergence for each numerical scheme considered. All methods exhibit excellent agreement with the exact solution for dispersion numbers below the critical threshold. Among them, the CCS8 consistently achieves an error nearly one order of magnitude smaller than that of CNCS8, thereby demonstrating superior accuracy within the stable dispersion regime.\par
To further investigate the dispersion properties of the numerical schemes, Eq.~\eqref{Eqn:Eg1} is solved for the following initial conditions using all three schemes with $ N = 40 $:
\begin{itemize}
    \item Case (a): $ u(x,0) = \sin(x) $; $ kh = 0.1571 $ (low wavenumber),
    \item Case (b): $ u(x,0) = \sin(4x) $; $ kh = 0.6283 $ (intermediate wavenumber),
    \item Case (c): $ u(x,0) = \sin(6x) $; $ kh = 0.9425 $ (high wavenumber).
\end{itemize}
The phase error behavior is illustrated in Fig.~\ref{Fig:Eg1b}, which presents the evolution of the maximum absolute error over time for cases (a), (b), and (c). Among the three schemes, CCS consistently demonstrates the highest accuracy in capturing wave propagation. From the corresponding dispersion relation plot (see Fig.~\ref{Fig:Spectral}), it is evident that the CCS  maintains a reliable \textit{DRP} region up to $ kh \leq \pi $, whereas the CNCS  performs well only up to $ kh \leq 0.7 $ for small values of $ N_c $.\par
\begin{figure}[htbp!]  
    \centering
    \begin{minipage}[b]{0.3\linewidth}
    \includegraphics[width=\linewidth, trim={0 0 0 0}, clip]{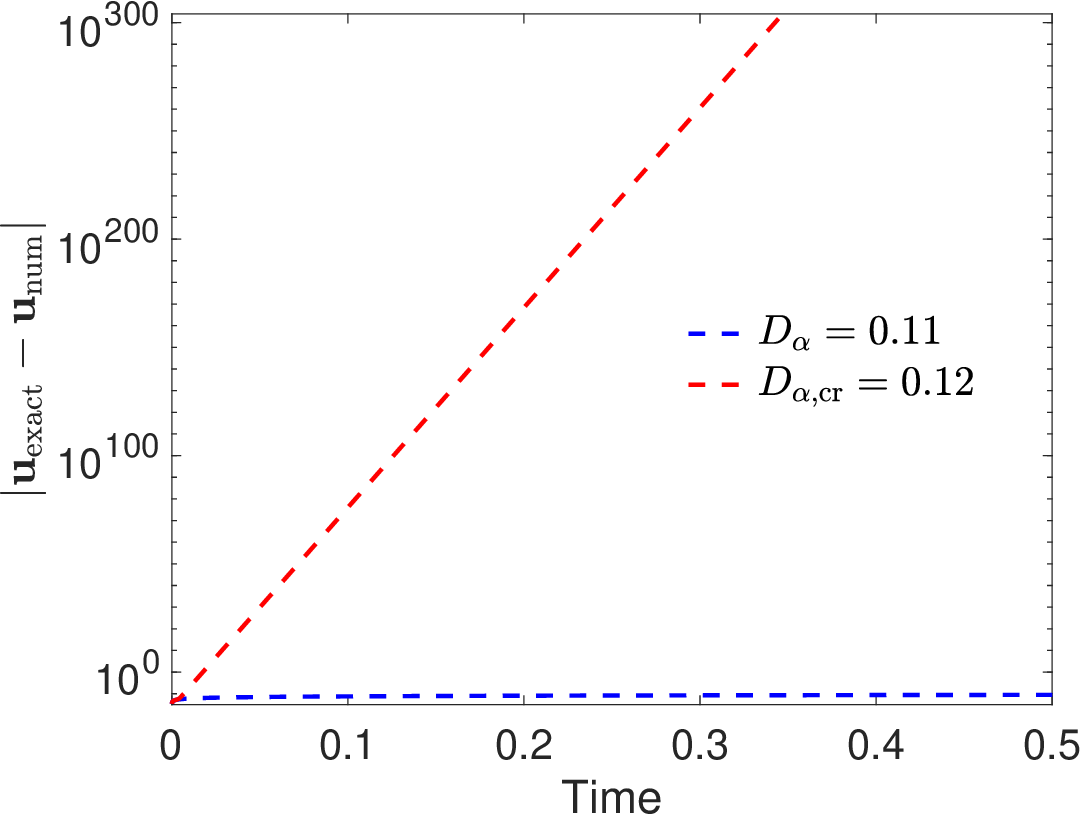}
      \subcaption{SSPRK3-CNCS6}
    \end{minipage}\hfill
    \begin{minipage}[b]{0.3\linewidth}
    \includegraphics[width=\linewidth, trim={0 0 0 0}, clip]{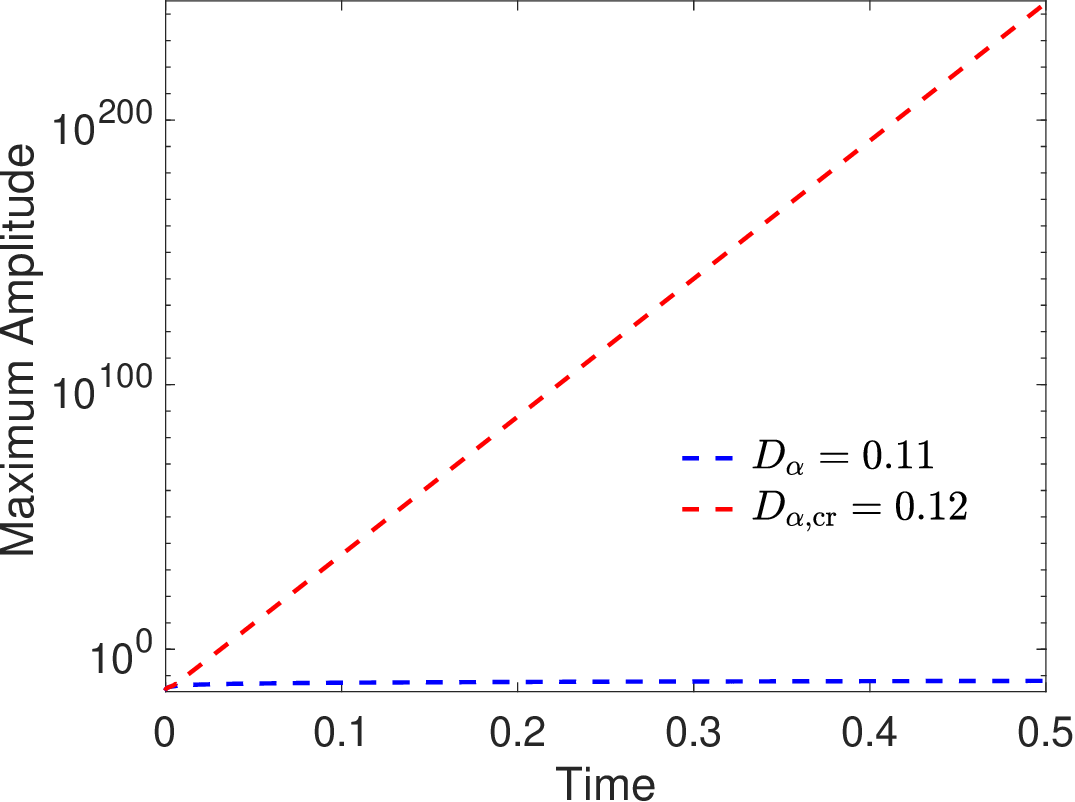}
      \subcaption{SSPRK3-CNCS8}
    \end{minipage}\hfill
    \begin{minipage}[b]{0.3\linewidth}
      \includegraphics[width=\linewidth, trim={0 0 0 0}, clip]{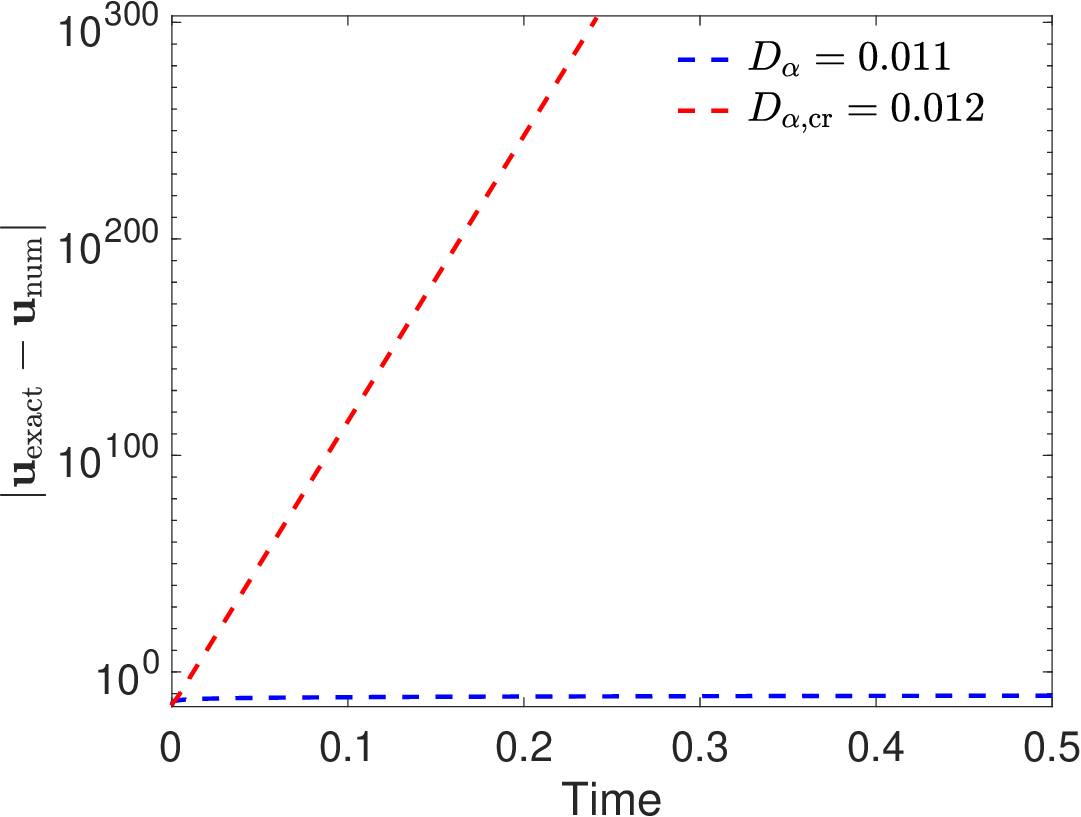}
      \subcaption{SSPRK3-CCS8}
    \end{minipage}\hfill    
    \caption{Variation of the maximum absolute difference between the exact and numerical solutions over time for example~\ref{example:1} with the initial condition $u(x,0) = \sin x$.}
    \label{Fig:Eg1a}
\end{figure}
\begin{figure}[htbp!]  
    \centering
    \begin{minipage}[b]{0.3\linewidth}
    \includegraphics[width=\linewidth, trim={0 0 0 0}, clip]{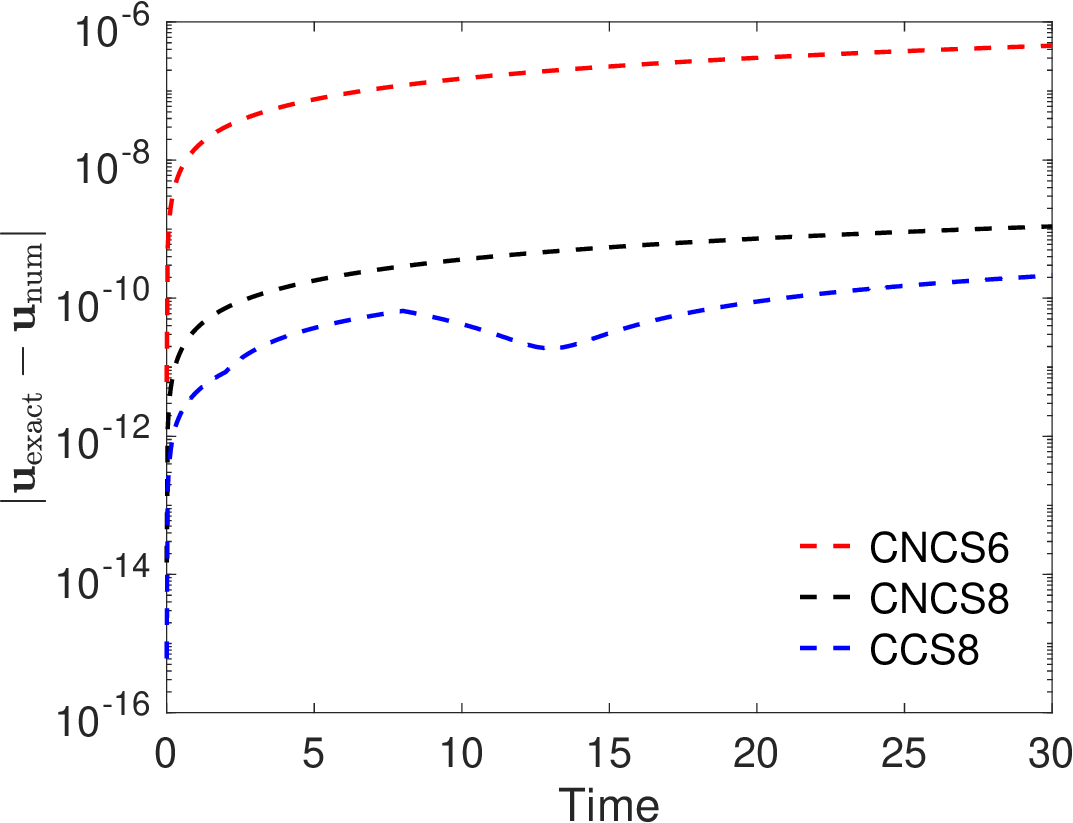}
      \subcaption{}
    \end{minipage}\hfill
    \begin{minipage}[b]{0.3\linewidth}
      \includegraphics[width=\linewidth, trim={0 0 0 0}, clip]{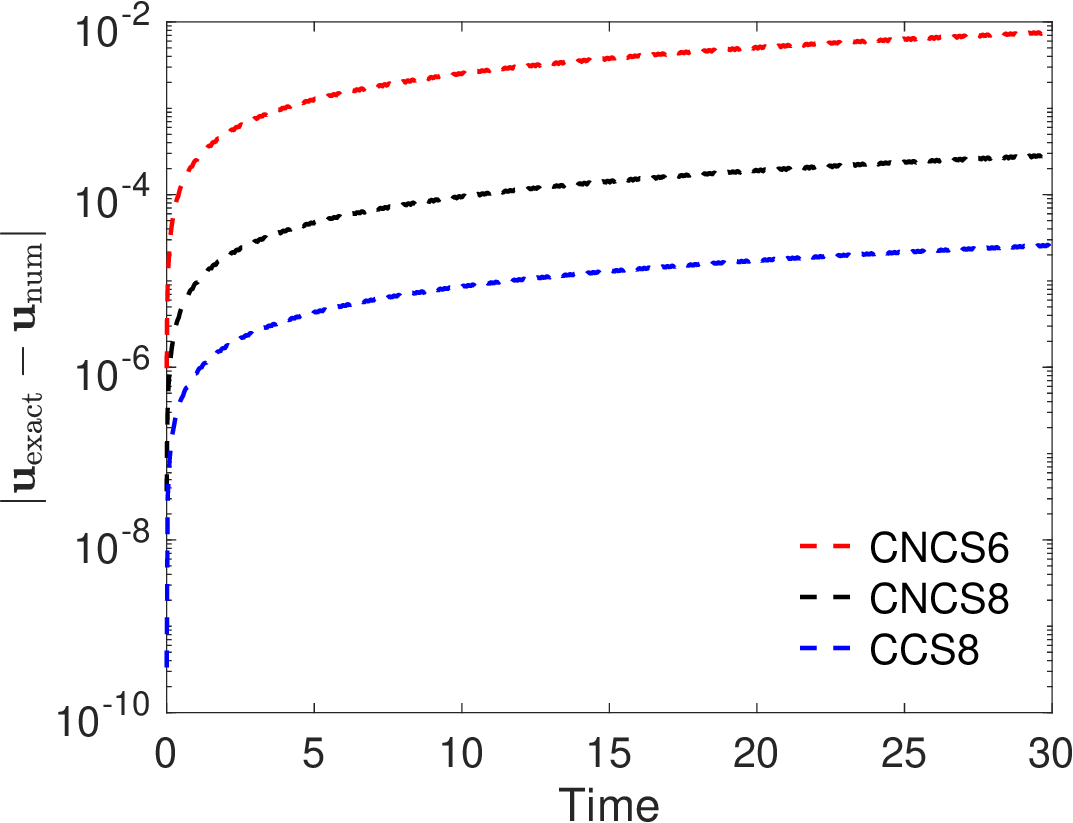}
      \subcaption{}
    \end{minipage}\hfill
    \begin{minipage}[b]{0.3\linewidth}
    \includegraphics[width=\linewidth, trim={0 0 0 0}, clip]{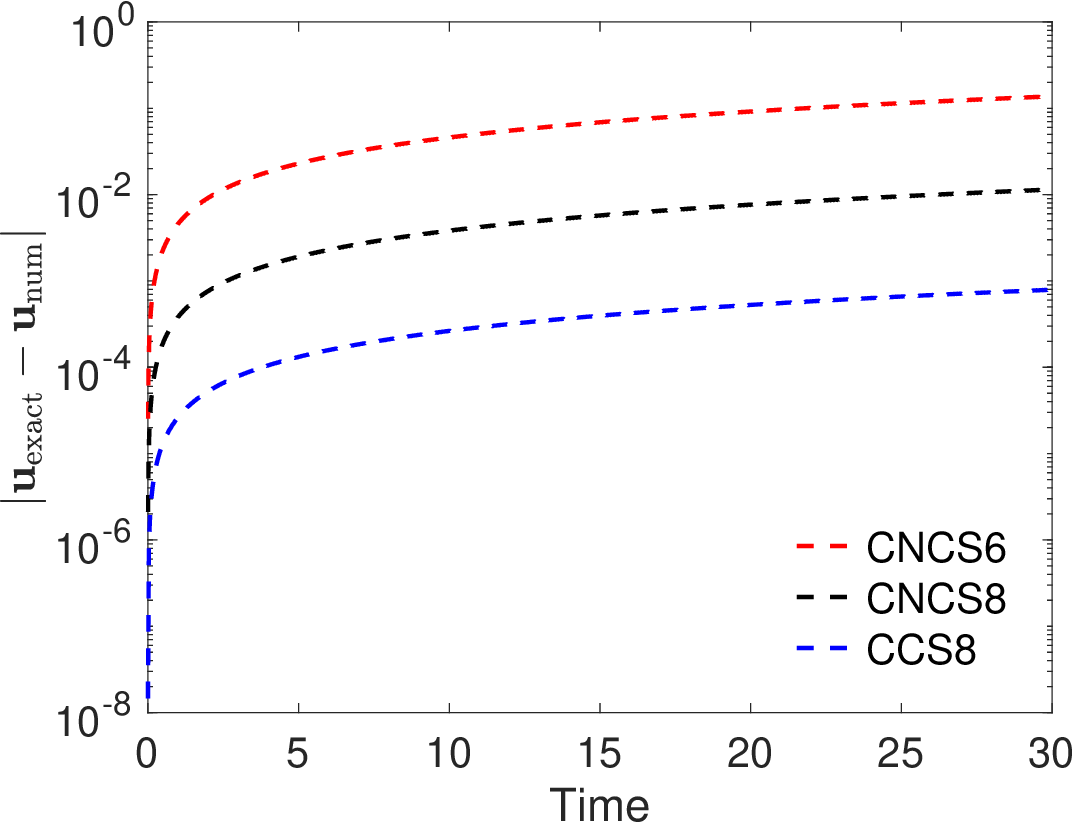}
      \subcaption{}
    \end{minipage}\hfill
    \caption{Evolution of the maximum error over time for example~\ref{example:1} with the initial conditions: (a) $u(x,0) = \sin x$; (b) $u(x,0) = \sin 4x$; (c) $u(x,0) = \sin 6x$.}
    \label{Fig:Eg1b}
\end{figure}
\begin{table}
\captionsetup{position=above} 
\caption{Error and order of convergence for example~\ref{example:1} with the initial condition $u(x,0) = \sin x$ at time $t= 0.5$.}
\label{Table:Eg1}
\centering
\begin{tabular*}{\textwidth}{@{\extracolsep{\fill}} l *{8}{c} }
\toprule
&\multicolumn{2}{c}{\textbf{SSPRK3-CNCS6}} & &
\multicolumn{2}{c}{\textbf{SSPRK3-CNCS8}} &
\multicolumn{2}{c}{\textbf{SSPRK3-CCS8}} \\
\cmidrule{2-3} \cmidrule{5-6} \cmidrule{7-8}
\textbf{N} & $\boldsymbol{L^{\infty}}$\textbf{-error} & \textbf{Rate}
&\textbf{N} & $\boldsymbol{L^{\infty}}$\textbf{-error} & \textbf{Rate}
& $\boldsymbol{L^{\infty}}$\textbf{-error} & \textbf{Rate}\\
\midrule
  10 & 3.1813e-05 & -& 10 & 1.1843e-06 & - & 1.0856e-07 & - \\
 20 & 4.8359e-07 & 6.0397&15 & 4.6512e-08 & 7.9839 & 4.9274e-09 & 7.6270 \\
 40 & 7.5647e-09 & 5.9984&20 & 4.6233e-09 & 8.0249 & 5.2131e-10 & 7.8079 \\
 60 & 6.6331e-10 & 6.0030& 25 & 7.7612e-10 & 7.9973 & 9.0868e-11 & 7.8288 \\
 80 & 1.1811e-10 & 5.9983& 30 & 1.8020e-10 & 8.0092 & 2.1415e-11 & 7.9274 \\
 100 & 3.1158e-11 & 5.9718&35 & 5.2398e-11 & 8.0129 & 6.2660e-12 & 7.9723 \\
 120 & 1.0551e-11 & 5.9390& 40 & 1.7973e-11 & 8.0133 & 2.2282e-12 & 7.7430 \\
\bottomrule
\end{tabular*}
\end{table}

\end{example}
\newpage
\begin{example}\label{example:2}
\normalfont
Consider the one-dimensional nonlinear soliton solution of the KdV equation:
\begin{equation}
    u_t + 3(u^2)_x + u_{xxx} = 0, \quad x \in [-10,12], \quad t \ge 0.
\end{equation}

\subsection*{Single Soliton Propagation}
In the case of a single soliton, the initial condition and exact solution are given by:
\begin{align}\label{IC:2a}
    u(x,0) &= 2\sech^2(x), \quad x \in [-10,12], \\
    u(x,t) &= 2\sech^2(x - 4t), \quad x \in [-10,12], \quad t \in [0, 0.5]. \notag
\end{align}
As observed in the linear case, the nonlinear case exhibits instability at the critical dispersion number. Specifically, blow-up occurs for CNCS at $D_{\alpha,\mathrm{cr}} = 0.12$ and for CCS at $D_{\alpha,\mathrm{cr}} = 0.012$, although these results are not presented here. This is consistent with the predictions made by the GSA. Henceforth, we fix $\text{CFL} = 0.11$ for CNCS and $\text{CFL} = 0.011$ for CCS. The resulting solutions and corresponding errors for $N = 140$ at different time instances $t = 0$, $0.25$, and $0.5$ are illustrated in Fig.~\ref{Fig:Eg2a}. As expected, the solution exhibits a soliton traveling at a constant velocity $c = 4$, preserving its shape and amplitude. Table~\ref{Table:Eg2a} provides a detailed comparison of the $L^{\infty}$ errors and the corresponding convergence rates for various grid sizes: $N$ ranging from 50 to 350 for CNCS6, and from 20 to 140 for CNCS8 and CCS8. The tabulated results demonstrate the superiority of the CCS  in terms of accuracy. It consistently outperforms the CNCS, producing errors at least an order of magnitude smaller across all tested grid resolutions. However, this enhanced accuracy comes at a computational cost; CCS demands significantly more resources, with runtimes exceeding those of CNCS by a factor of at least two. While CCS exhibits a convergence rate approaching the seventh order for larger values of $N$, CNCS maintains a consistent rate close to the theoretical eighth order.
\begin{table}
\captionsetup{position=above} 
\caption{Error and order of convergence for example~\ref{example:2} with the initial condition (\ref{IC:2a}) at time $t= 0.5$.}
\label{Table:Eg2a}
\centering
\begin{tabular*}{\textwidth}{@{\extracolsep{\fill}} l *{8}{c} }
\toprule
&\multicolumn{2}{c}{\textbf{SSPRK3-CNCS6}} & &
\multicolumn{2}{c}{\textbf{SSPRK3-CNCS8}} &
\multicolumn{2}{c}{\textbf{SSPRK3-CCS8}} \\
\cmidrule{2-3} \cmidrule{5-6} \cmidrule{7-8}
\textbf{N} & $\boldsymbol{L^{\infty}}$\textbf{-error} & \textbf{Rate}
&\textbf{N} & $\boldsymbol{L^{\infty}}$\textbf{-error} & \textbf{Rate}
& $\boldsymbol{L^{\infty}}$\textbf{-error} & \textbf{Rate}\\
\midrule
 50 & 5.0253e-03 & - & 20 & 5.4857e-01 & - & 2.0463e-02 & - & \\
 100 & 6.4986e-05 & 6.2729 & 40 & 1.2963e-02 & 5.4032 & 2.6303e-04 & 6.2817 & \\
 150 & 5.3198e-06 & 6.1725 & 60 & 3.3124e-04 & 9.0440 & 1.7383e-05 & 6.7003 & \\
 200 & 9.4989e-07 & 5.9887 & 80 & 3.3098e-05 & 8.0066 & 2.3507e-06 & 6.9549 & \\
 250 & 2.4744e-07 & 6.0283 & 100 & 5.7921e-06 & 7.8110 & 4.8861e-07 & 7.0399 & \\
 300 & 8.2175e-08 & 6.0461 & 120 & 1.3299e-06 & 8.0703 & 1.3223e-07 & 7.1689 & \\
 350 & 3.2686e-08 & 5.9805 & 140 & 3.7966e-07 & 8.1322 & 4.2366e-08 & 7.3836 & \\
\bottomrule
\end{tabular*}
\end{table}
\begin{figure}[htbp!]  
    \centering
    \begin{minipage}[b]{0.3\linewidth}
    \includegraphics[width=\linewidth, trim={0 0 0 0}, clip]{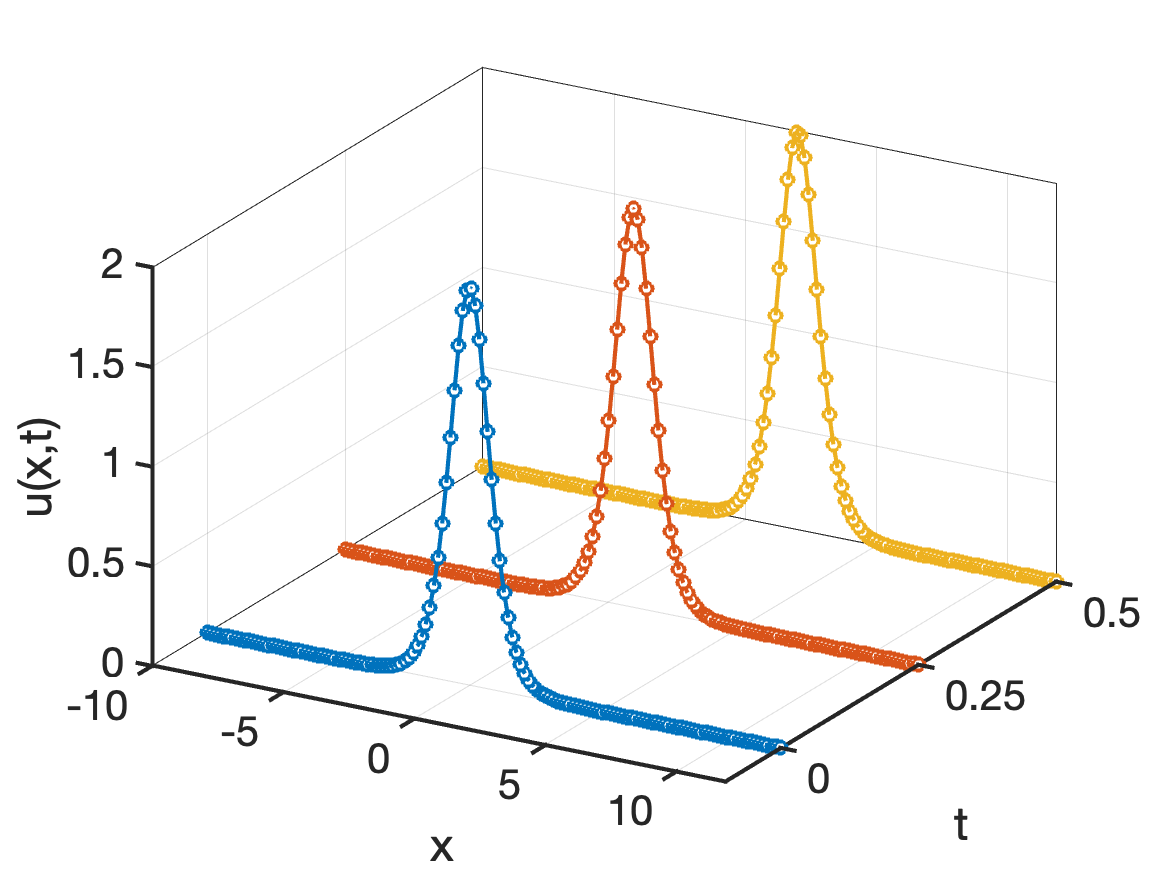}
      \subcaption{SSPRK3-CNCS6}
    \end{minipage}\hfill
    \begin{minipage}[b]{0.3\linewidth}
      \includegraphics[width=\linewidth, trim={0 0 0 0}, clip]{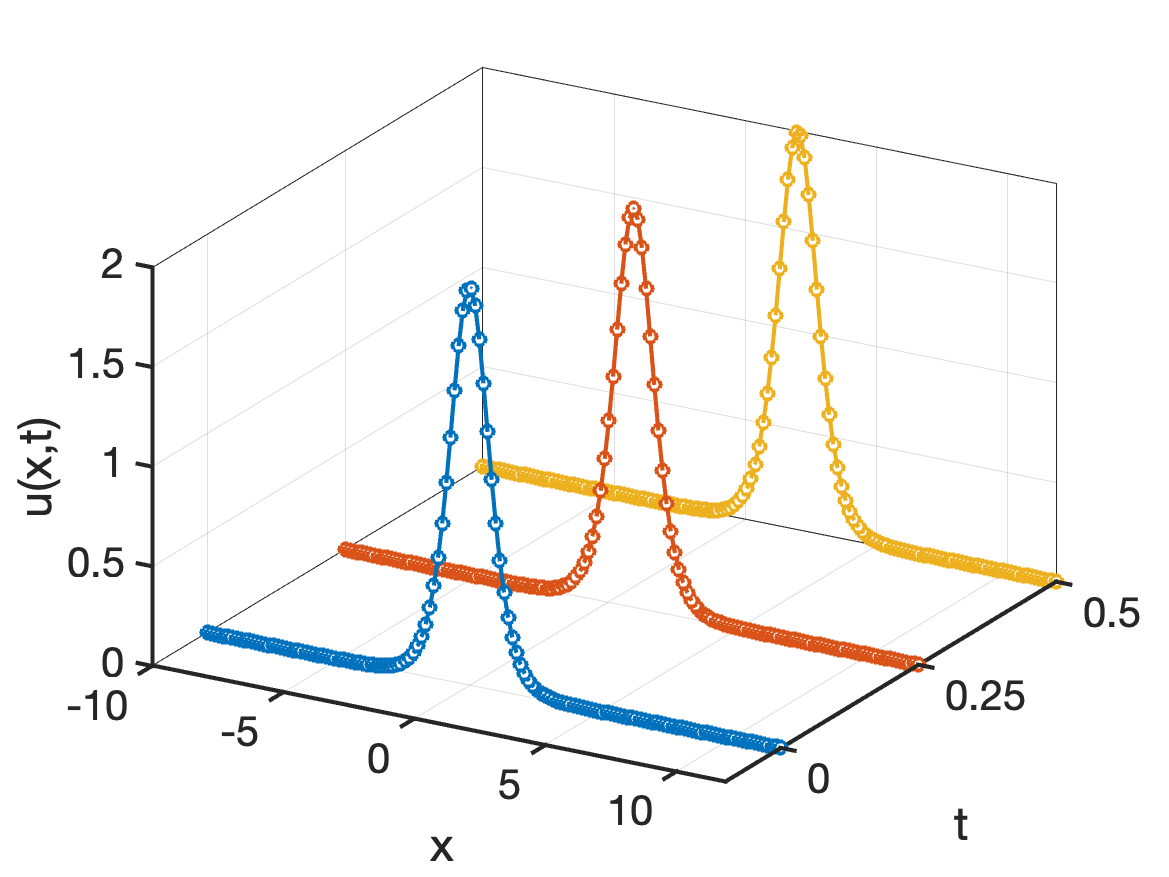}
      \subcaption{SSPRK3-CNCS8}
    \end{minipage}\hfill
    \begin{minipage}[b]{0.3\linewidth}
    \includegraphics[width=\linewidth, trim={0 0 0 0}, clip]{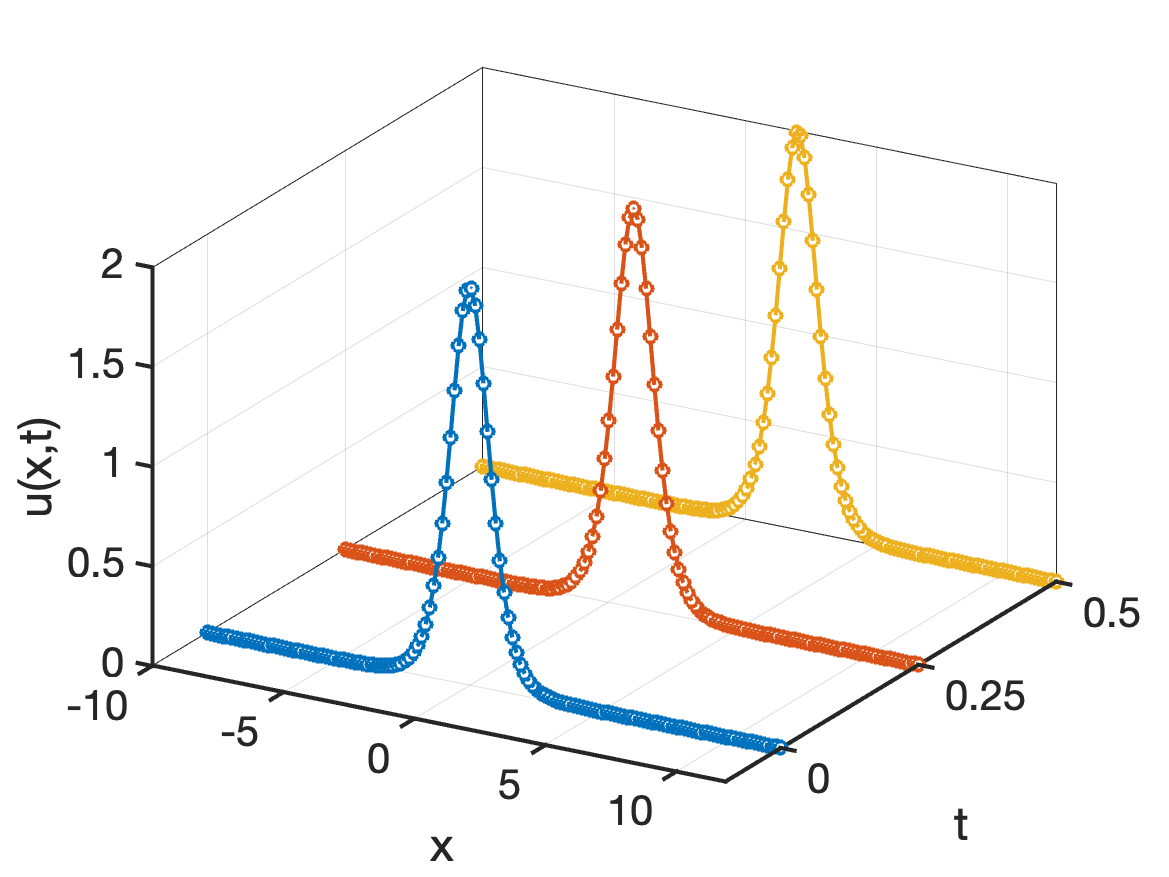}
      \subcaption{SSPRK3-CCS8}
    \end{minipage}\hfill
        \begin{minipage}[b]{0.3\linewidth}
    \includegraphics[width=\linewidth, trim={0 0 0 0}, clip]{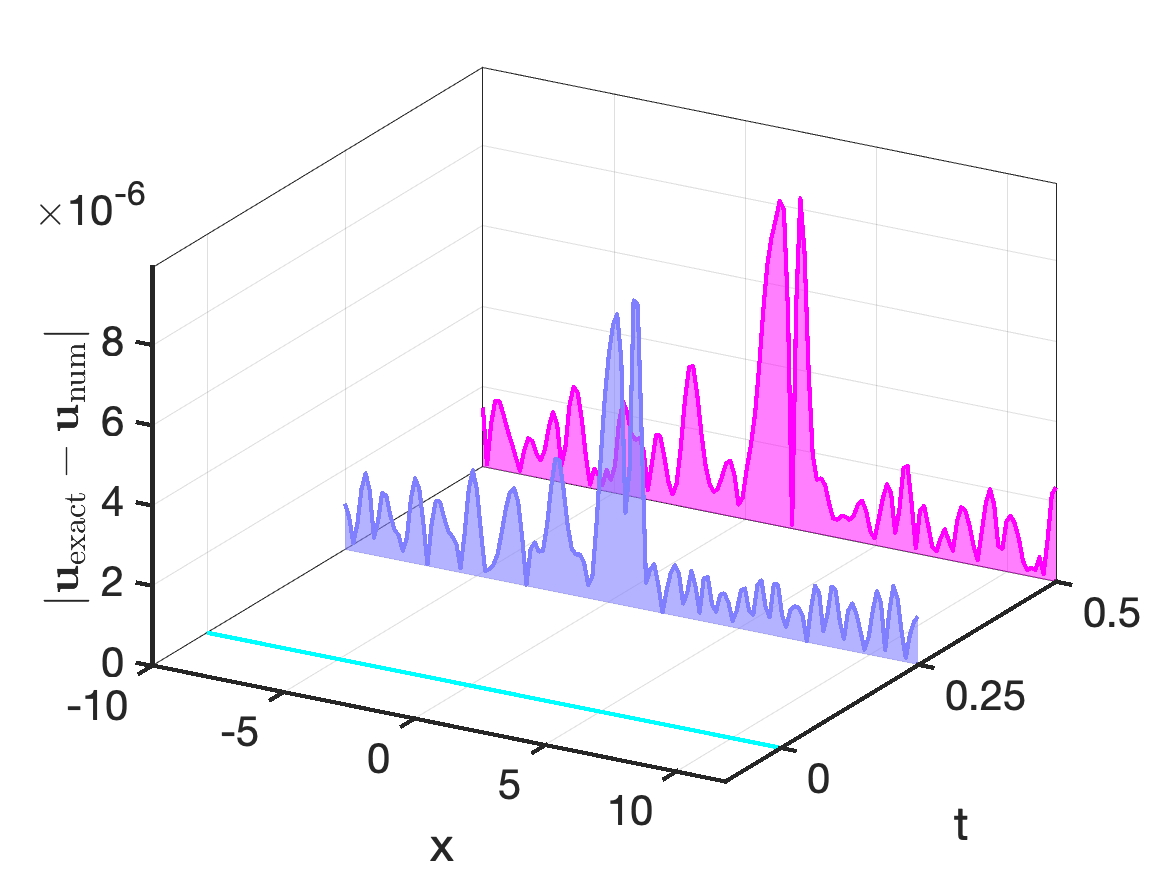}
      \subcaption{SSPRK3-CNCS6}
    \end{minipage}\hfill
    \begin{minipage}[b]{0.3\linewidth}
      \includegraphics[width=\linewidth, trim={0 0 0 0}, clip]{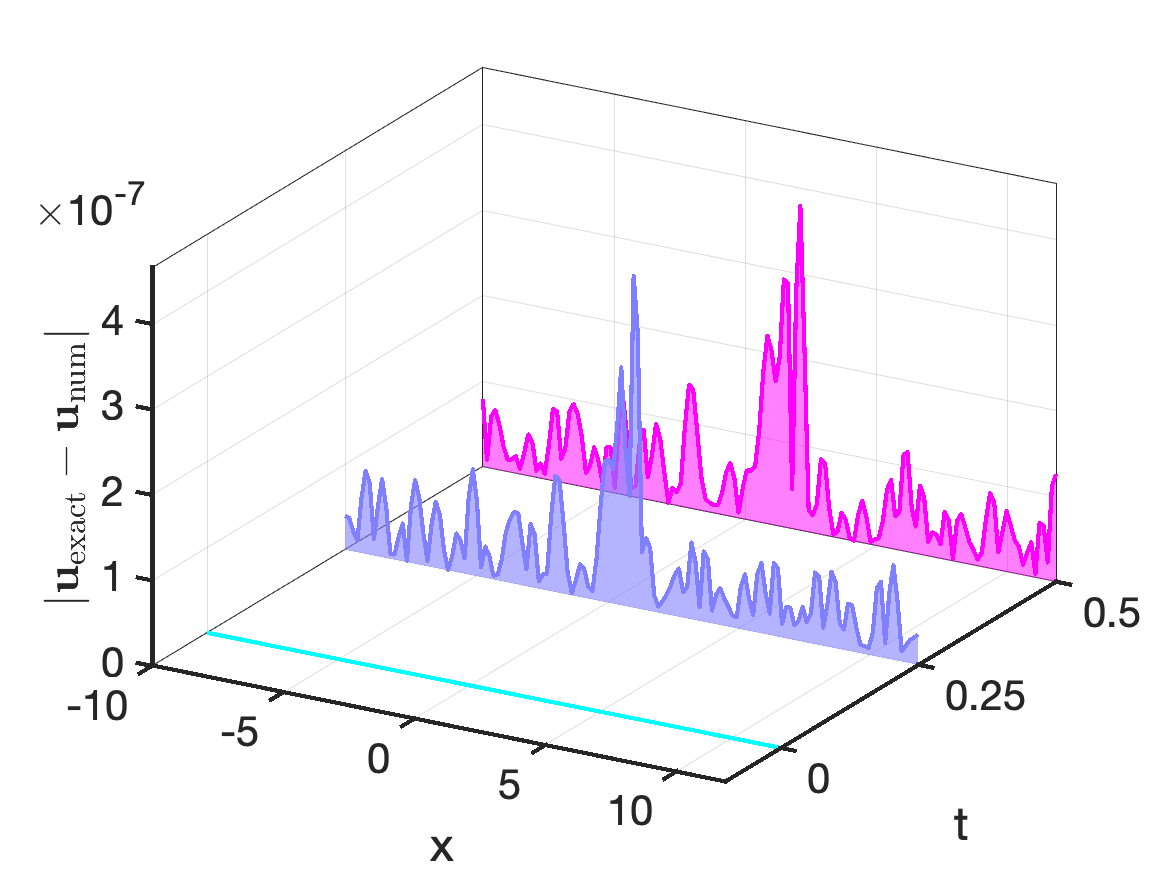}
      \subcaption{SSPRK3-CNCS8}
    \end{minipage}\hfill
    \begin{minipage}[b]{0.3\linewidth}
    \includegraphics[width=\linewidth, trim={0 0 0 0}, clip]{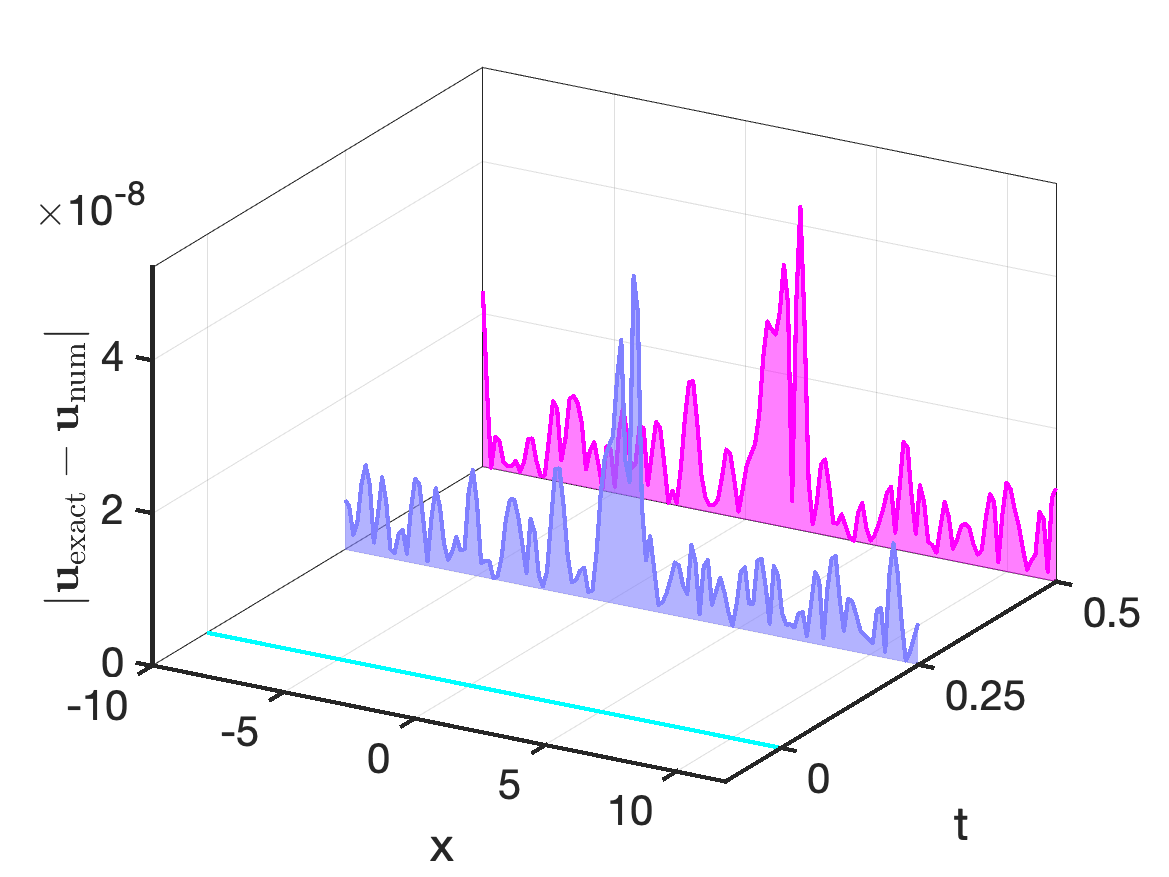}
      \subcaption{SSPRK3-CCS8}
    \end{minipage}\hfill
    \caption{Solutions and errors for the initial condition (\ref{IC:2a}) of example \ref{example:2}. The first row shows exact (solid) and numerical ($\circ$) solutions at $t=0, 0.25, 0.5$, and the second row shows errors for CNCS6, CNCS8, and CCS8 for $N=140$.}
    \label{Fig:Eg2a}
\end{figure}
\subsection*{Double soliton interaction} 
The second simulation involves an initial single wave at $t = 0$, which evolves into two distinct solitons of different amplitudes~\cite{yoneyama1984korteweg}. The problem admits the exact solution:
\begin{equation}\label{IC:2b}
    u(x,t) = 12 \dfrac{3 + 4 \cosh(2x - 8t) + \cosh(4x - 64t)}{\left(3 \cosh(x - 28t) + \cosh(3x - 36t)\right)^2}, \quad x \in [-10, 20], \quad t \in [0, 0.5].
\end{equation}
The KdV equation is numerically solved with the initial condition $u(x,0) = 6 \sech^2(x)$ at time $t = 0$, defined over the spatial domain $\Omega = (-10, 20)$. The corresponding results, summarized in Table~\ref{Table:Eg2b}, demonstrate the expected convergence trends.  Fig.~\ref{Fig:Eg2b} illustrates the numerical solution and the corresponding errors at time instances $t = 0$, $0.25$, and $0.5$ for $N = 400$.
\begin{table}
\captionsetup{position=above} 
\caption{Error and order of convergence for example~\ref{example:2} with the initial condition (\ref{IC:2b}) at time $t= 0.5$.}
\label{Table:Eg2b}
\centering
\begin{tabular*}{\textwidth}{@{\extracolsep{\fill}} l *{8}{c} }
\toprule
&\multicolumn{2}{c}{\textbf{SSPRK3-CNCS6}} & &
\multicolumn{2}{c}{\textbf{SSPRK3-CNCS8}} &
\multicolumn{2}{c}{\textbf{SSPRK3-CCS8}} \\
\cmidrule{2-3} \cmidrule{5-6} \cmidrule{7-8}
\textbf{N} & $\boldsymbol{L^{\infty}}$\textbf{-error} & \textbf{Rate}
&\textbf{N} & $\boldsymbol{L^{\infty}}$\textbf{-error} & \textbf{Rate}
& $\boldsymbol{L^{\infty}}$\textbf{-error} & \textbf{Rate}\\
\midrule
 100 & 2.7679e-01 & - & 100 & 1.1777e-01 & - & 5.1118e-03 & - & \\
 300 & 4.4723e-04 & 5.8510 & 200 & 6.9143e-04 & 7.4122 & 4.4629e-05 & 6.8397 & \\
 500 & 2.0677e-05 & 6.0178 & 300 & 2.5310e-05 & 8.1575 & 2.4851e-06 & 7.1228 & \\
 700 & 2.7176e-06 & 6.0310 & 400 & 2.4940e-06 & 8.0551 & 3.0127e-07 & 7.3347 & \\
 900 & 5.9819e-07 & 6.0227 & 500 & 4.1447e-07 & 8.0425 & 5.7577e-08 & 7.4163 & \\
 1100 & 1.7915e-07 & 6.0082 & 600 & 9.5634e-08 & 8.0433 & 1.6045e-08 & 7.0081 & \\
\bottomrule
\end{tabular*}
\end{table}
\begin{figure}[htbp!]  
    \centering
    \begin{minipage}[b]{0.3\linewidth}
    \includegraphics[width=\linewidth, trim={0 0 0 0}, clip]{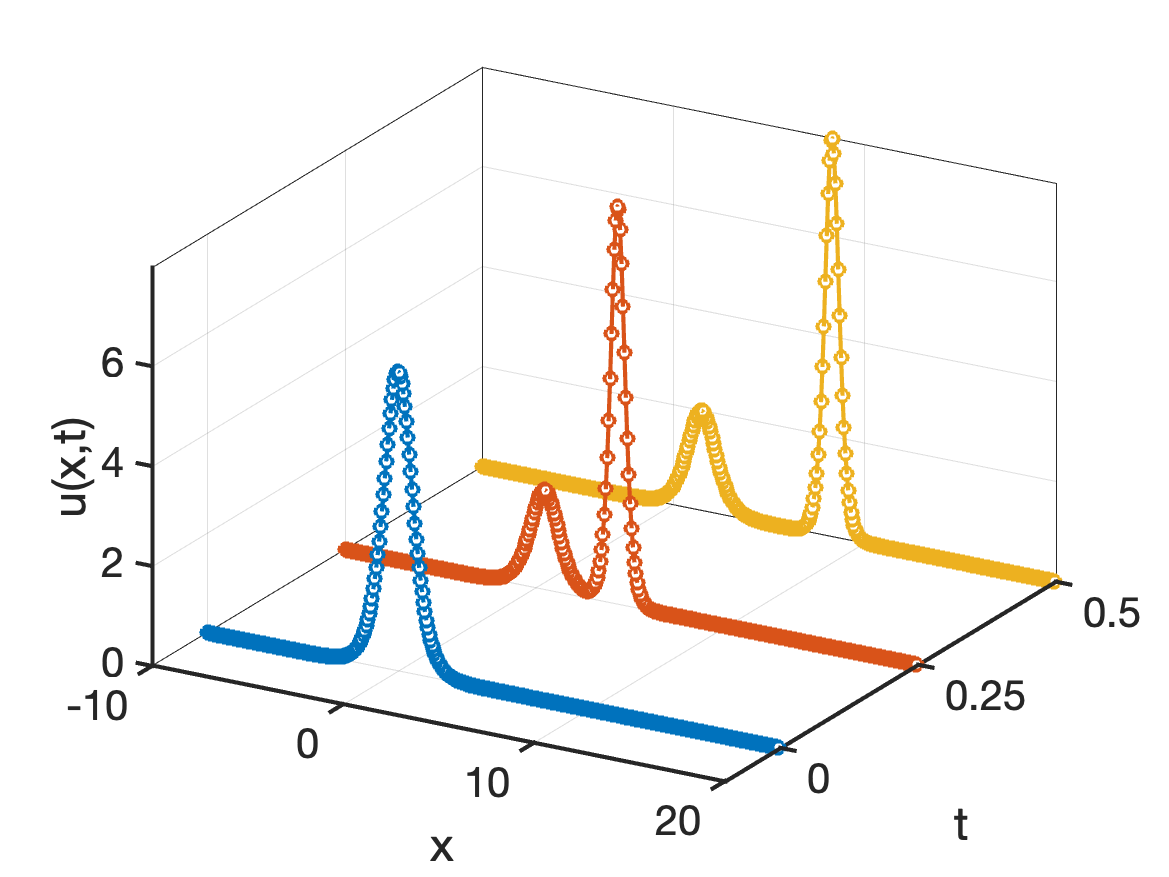}
      \subcaption{SSPRK3-CNCS6}
    \end{minipage}\hfill
    \begin{minipage}[b]{0.3\linewidth}
      \includegraphics[width=\linewidth, trim={0 0 0 0}, clip]{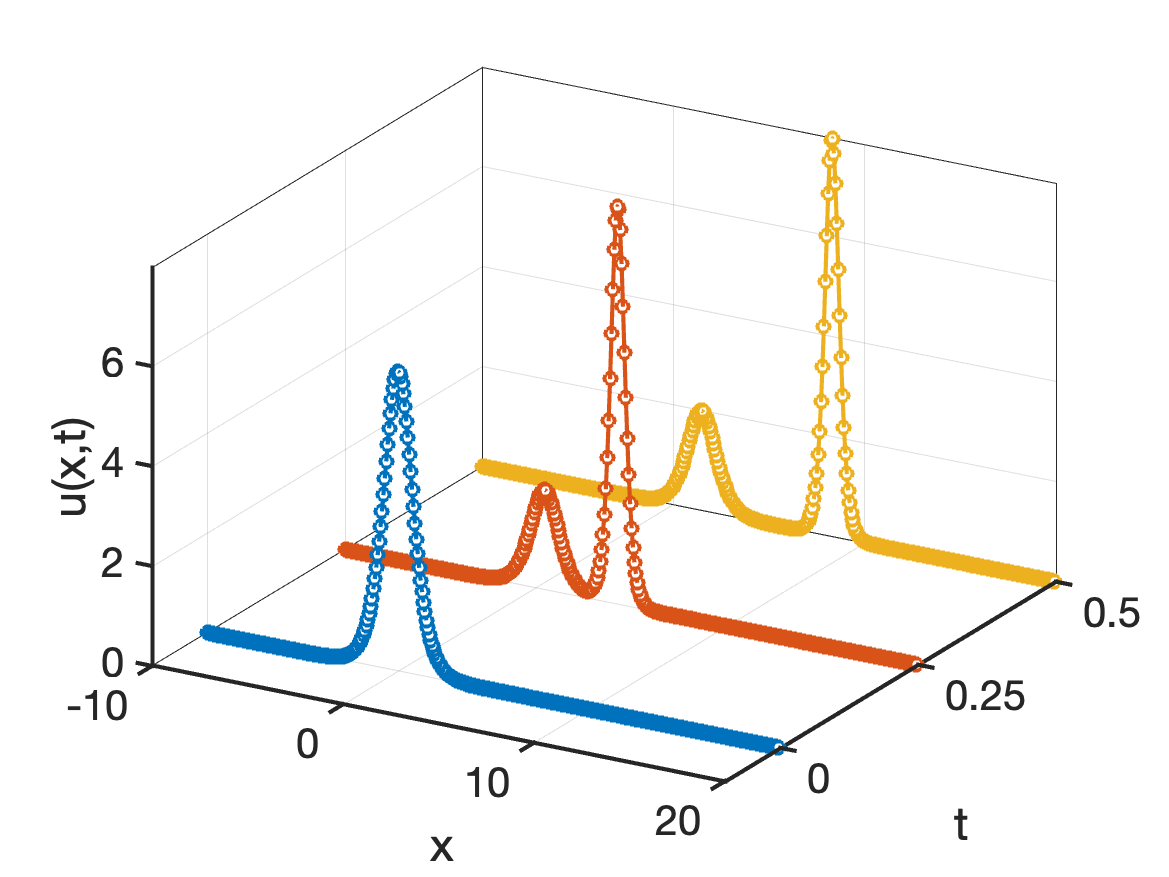}
      \subcaption{SSPRK3-CNCS8}
    \end{minipage}\hfill
    \begin{minipage}[b]{0.3\linewidth}
    \includegraphics[width=\linewidth, trim={0 0 0 0}, clip]{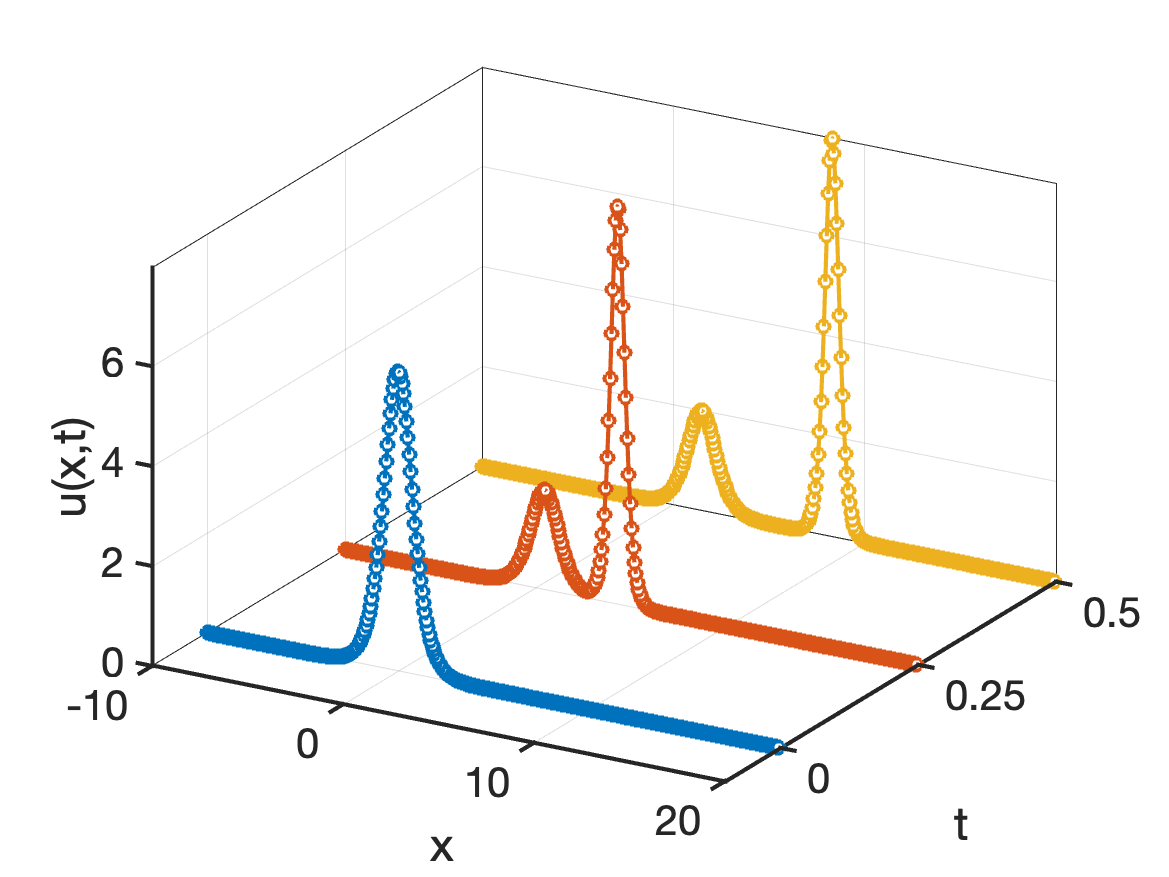}
      \subcaption{SSPRK3-CCS8}
    \end{minipage}\hfill
        \begin{minipage}[b]{0.3\linewidth}
    \includegraphics[width=\linewidth, trim={0 0 0 0}, clip]{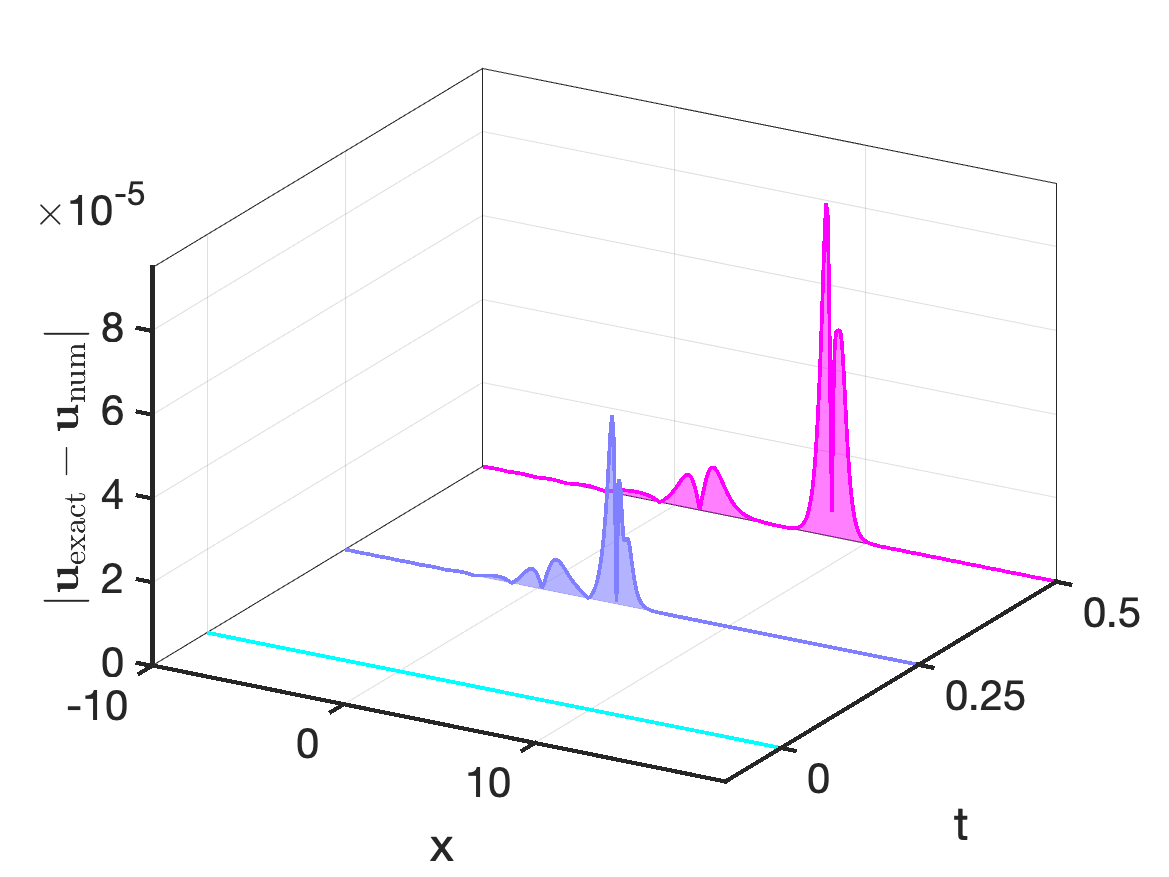}
      \subcaption{SSPRK3-CNCS6}
    \end{minipage}\hfill
    \begin{minipage}[b]{0.3\linewidth}
      \includegraphics[width=\linewidth, trim={0 0 0 0}, clip]{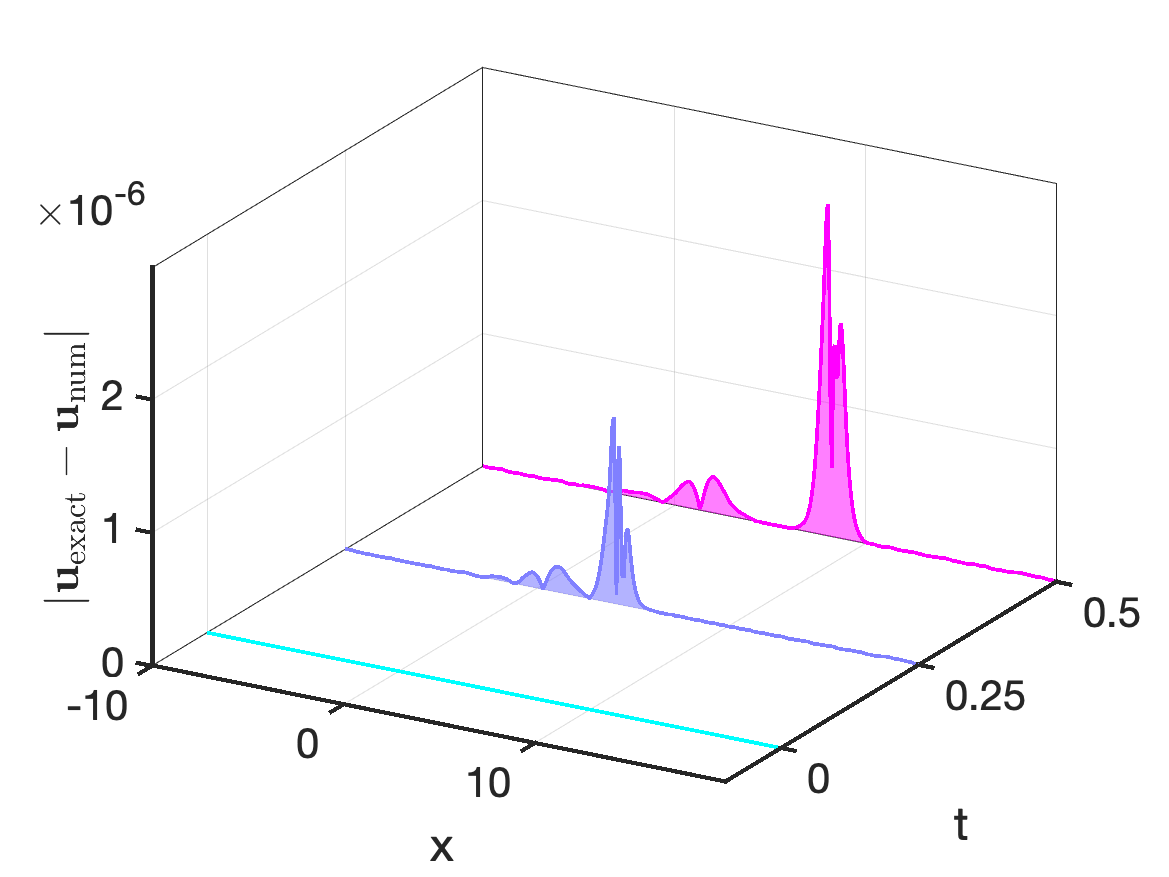}
      \subcaption{SSPRK3-CNCS8}
    \end{minipage}\hfill
    \begin{minipage}[b]{0.3\linewidth}
    \includegraphics[width=\linewidth, trim={0 0 0 0}, clip]{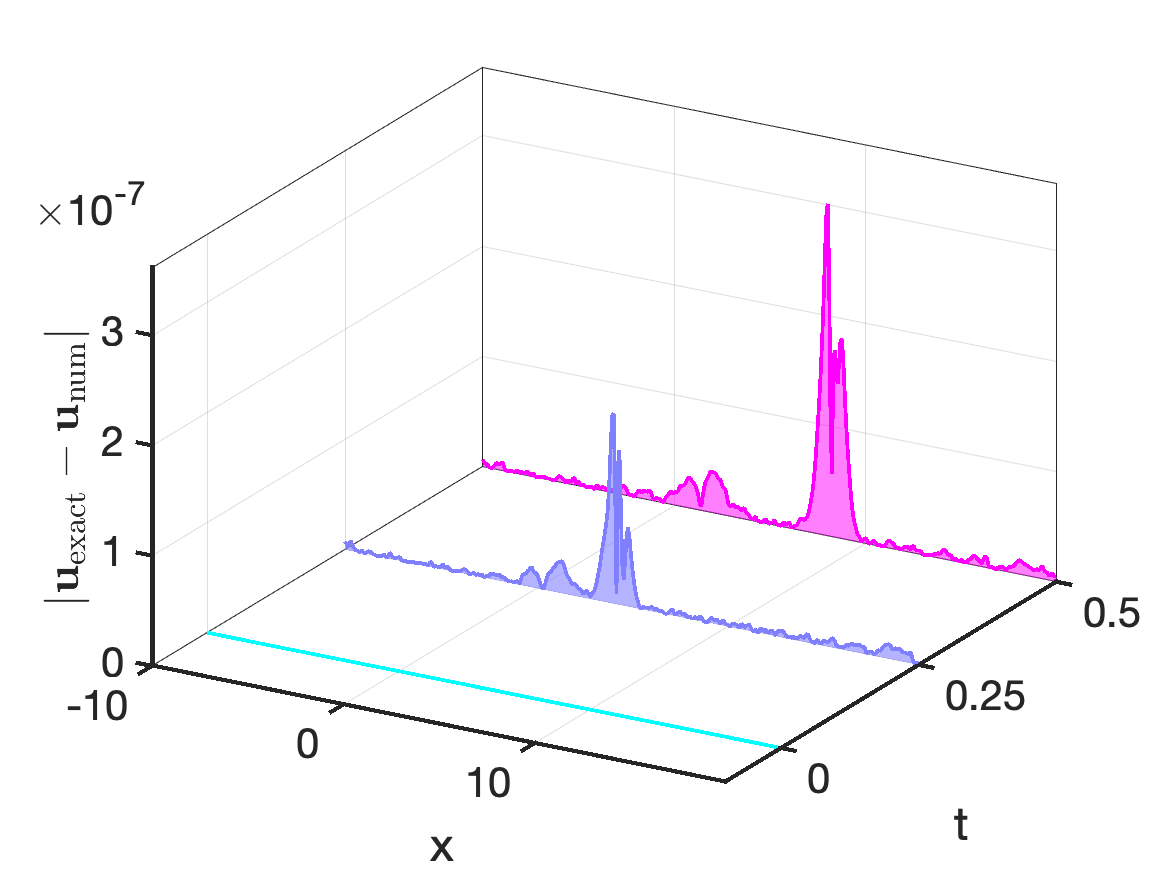}
      \subcaption{SSPRK3-CCS8}
    \end{minipage}\hfill
    \caption{Solutions and errors for the initial condition (\ref{IC:2b}) of example \ref{example:2}. The first row shows exact (solid) and numerical ($\circ$) solutions at $t=0, 0.25, 0.5$, and the second row shows errors for CNCS6, CNCS8, and CCS8 for $N=400$.}
    \label{Fig:Eg2b}
\end{figure}
\end{example}
\begin{example}\label{example:3}
\normalfont To analyze the dispersive features of solutions to nonlinear PDEs, we consider the classical KdV equation in the small-dispersion limit:
\begin{equation}
    u_t + \left(\dfrac{u^2}{2}\right)_x + \epsilon\, u_{xxx} = 0, \quad x \in [0,1], \quad t \ge 0,
\end{equation}
subject to two different types of initial data. The first is a smooth periodic profile:
\begin{equation}
    u(x,0) = 2 + 0.5 \sin(2\pi x), \quad x \in [0,1],
    \label{IC:3a}
\end{equation}
and the second is a discontinuous top-hat function:
\begin{equation}
    u(x,0) = 
    \begin{cases}
        1, & \text{if } 0.25 < x < 0.4, \\
        0, & \text{otherwise}.
    \end{cases}
    \label{IC:3b}
\end{equation}
For the smooth initial condition~\eqref{IC:3a}, the solution develops a dispersive shock structure over time. When the dispersion parameter $ \epsilon $ is small, high-frequency wavelets emerge near the steep gradients. Numerical simulations, shown in Fig.~\ref{Fig:Eg3a} for $ \epsilon = 10^{-4} (N=200), 10^{-5} (N=800), 10^{-6} (N=1500)$ at $ t= 0.5 $, demonstrate that the resulting oscillations are well resolved and remain free from spurious numerical artifacts. These dispersive features are captured, even on relatively coarse spatial grids, with the solution remaining stable and physically consistent before and after the formation of the dispersive shock.\par
The solution exhibits a distinct behaviour for the discontinuous initial profile~\eqref{IC:3b}. The initial jump evolves into a train of smooth, left-propagating dispersive waves. Fig.~\ref{Fig:Eg3b} illustrates this transition, showing the formation of continuous fine-scale structures that gradually evolve into solitary waves rather than classical shocks. This behavior is captured at time instances $ t = 0.01, 0.05, 0.1 $ using a spatial grid of $ N = 1500 $, with $ \epsilon = 10^{-4} $. Minor spurious oscillations appearing to the left of the discontinuity at early times are effectively removed by applying the filtering procedure described in~\cite{salian2024novel}, which enhances the clarity of the solution by suppressing non-physical artifacts while preserving the essential wave dynamics.
\begin{figure}[htbp!]  
\centering
\begin{minipage}[b]{0.3\linewidth}
\includegraphics[width=\linewidth, trim={0 0 0 0}, clip]{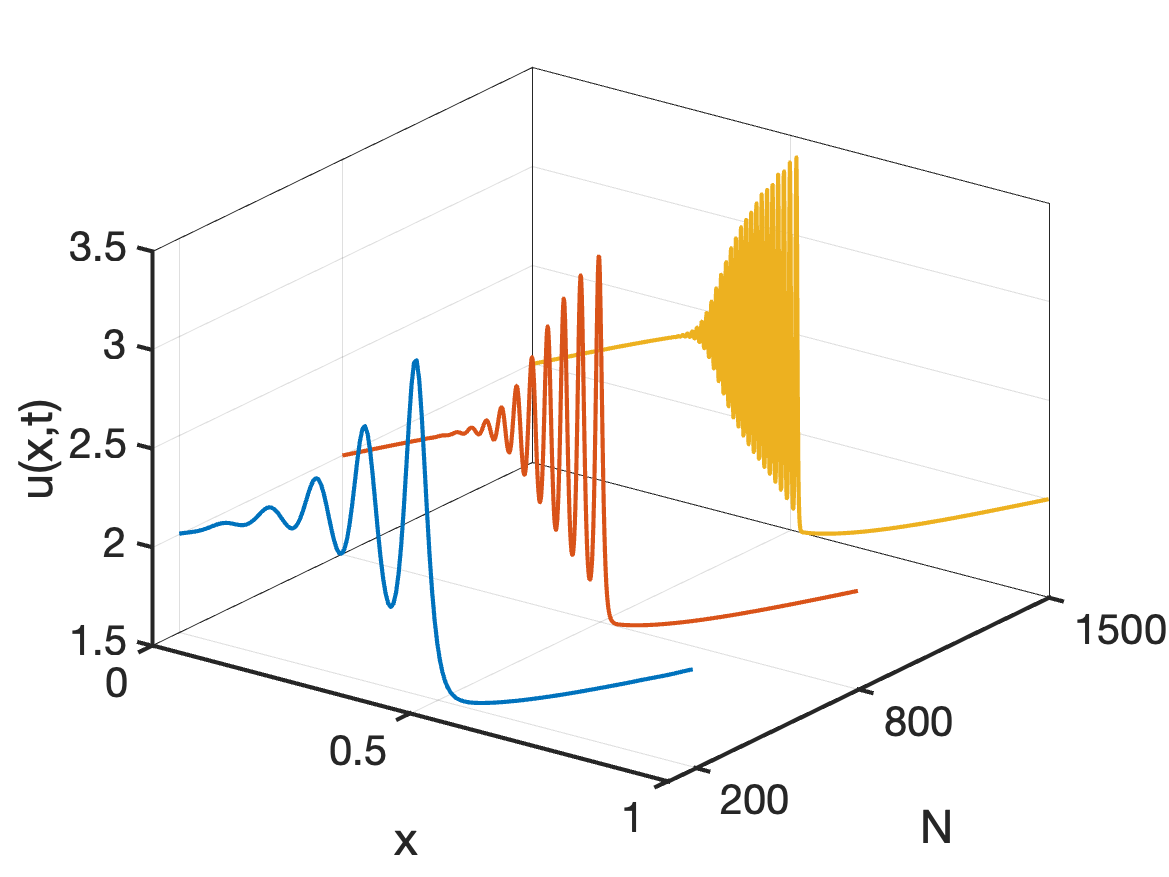}
  \subcaption{SSPRK3-CNCS6}
\end{minipage}\hfill
\begin{minipage}[b]{0.3\linewidth}
  \includegraphics[width=\linewidth, trim={0 0 0 0}, clip]{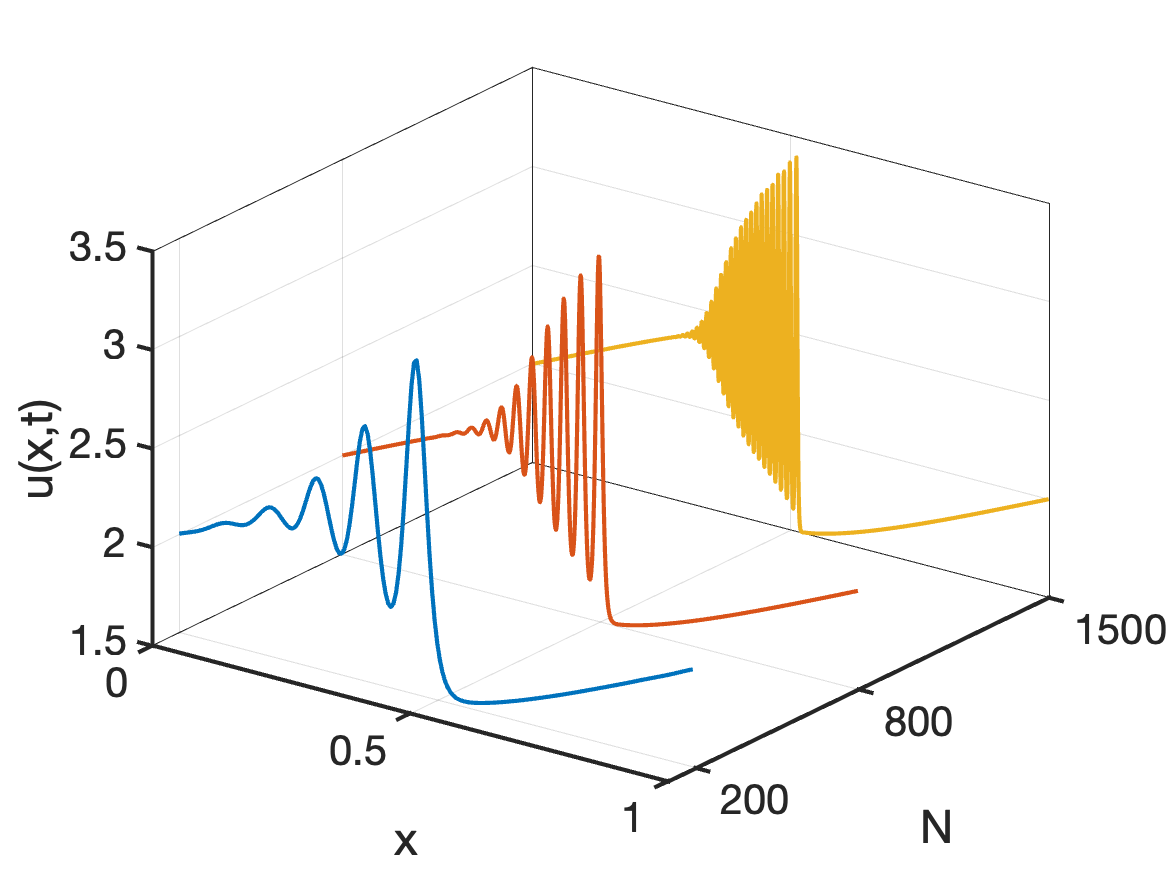}
  \subcaption{SSPRK3-CNCS8}
\end{minipage}\hfill
\begin{minipage}[b]{0.3\linewidth}
\includegraphics[width=\linewidth, trim={0 0 0 0}, clip]{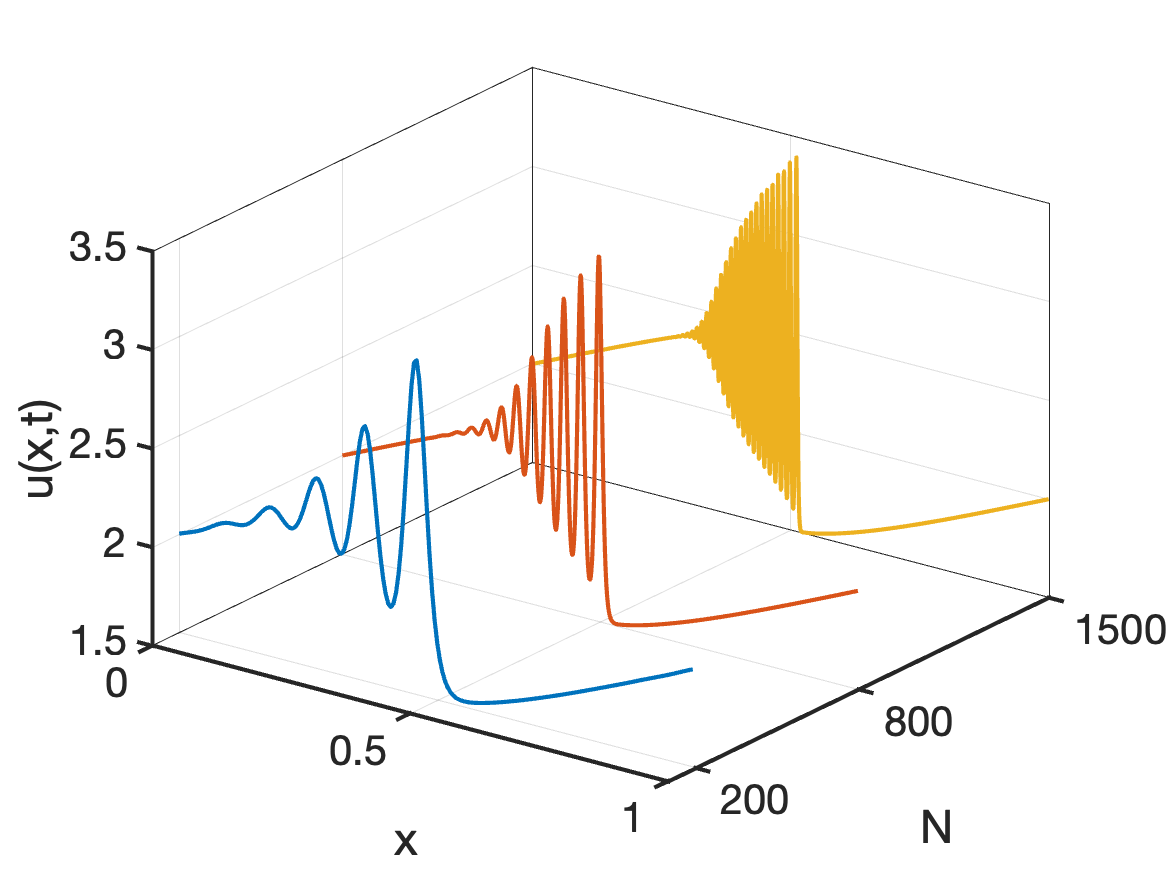}
  \subcaption{SSPRK3-CCS8}
\end{minipage}\hfill
\caption{Numerical solutions for initial condition (\ref{IC:3a}) of example \ref{example:3}.}
\label{Fig:Eg3a}
\end{figure}
\begin{figure}[htbp!]  
\centering
\begin{minipage}[b]{0.3\linewidth}
\includegraphics[width=\linewidth, trim={0 0 0 0}, clip]{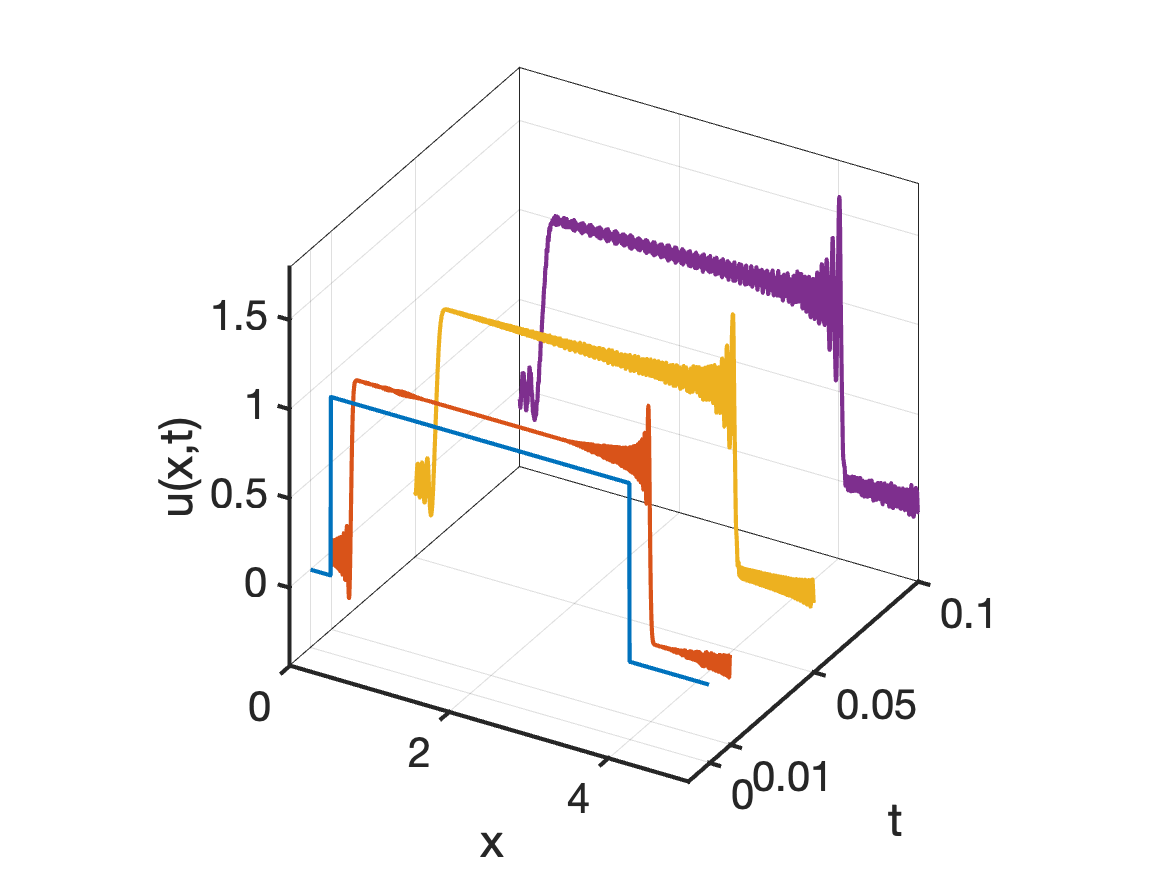}
  \subcaption{SSPRK3-CNCS6}
\end{minipage}\hfill
\begin{minipage}[b]{0.3\linewidth}
  \includegraphics[width=\linewidth, trim={0 0 0 0}, clip]{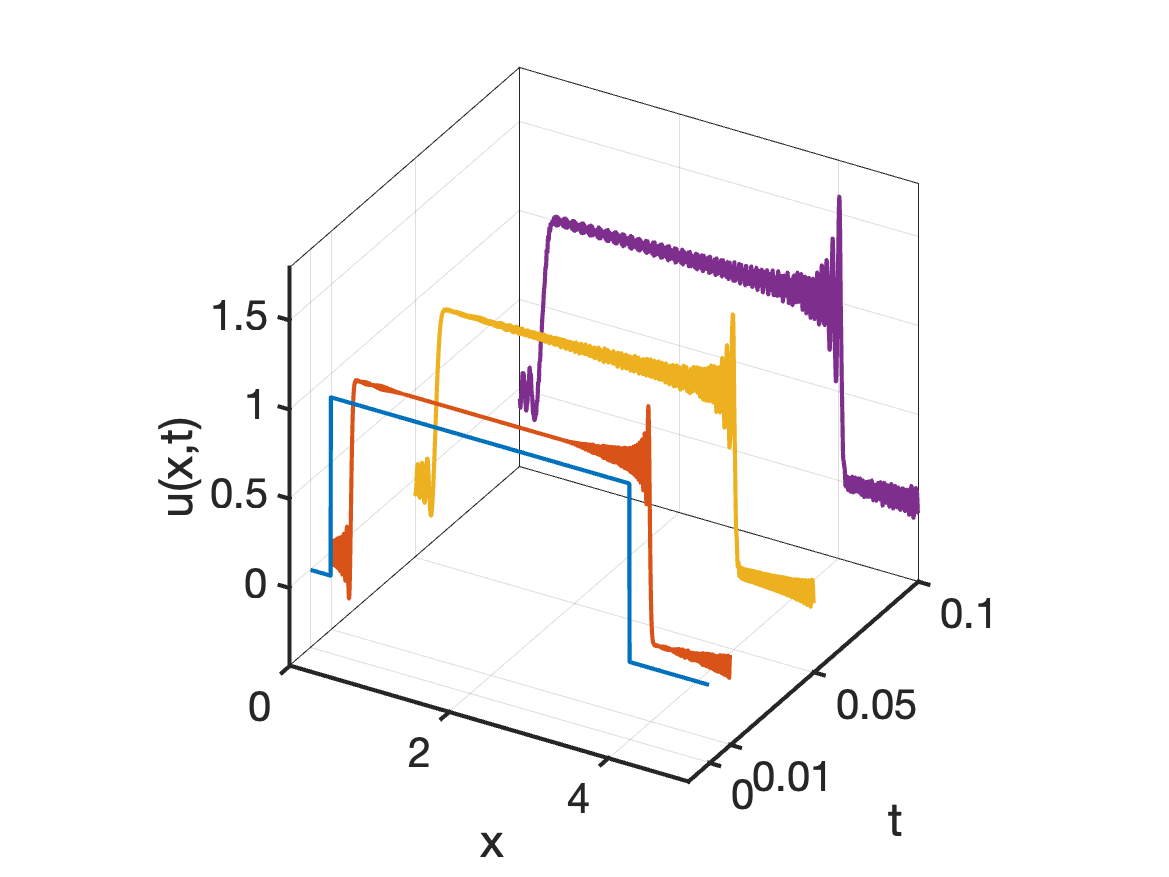}
  \subcaption{SSPRK3-CNCS8}
\end{minipage}\hfill
\begin{minipage}[b]{0.3\linewidth}
\includegraphics[width=\linewidth, trim={0 0 0 0}, clip]{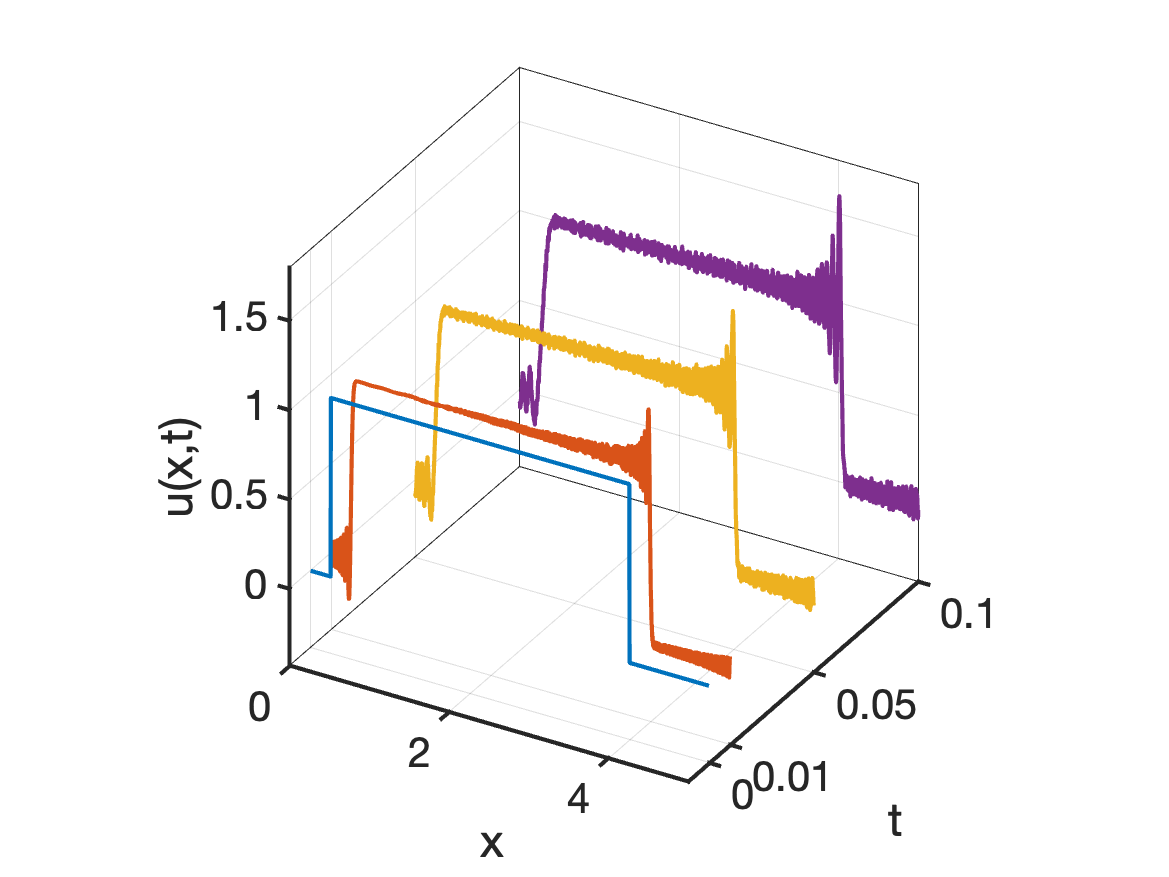}
  \subcaption{SSPRK3-CCS8}
\end{minipage}\hfill
\caption{Numerical solutions for initial condition (\ref{IC:3b}) of example \ref{example:3}.}
\label{Fig:Eg3b}
\end{figure}
\end{example}
\begin{example}\label{example:4}
\normalfont We consider the following nonlinear mKdV equation~\cite{ismail2020numerical}:
\begin{equation}
u_t + \mu\, u^2 u_x + \epsilon\, u_{xxx} = 0.
\end{equation}
This equation can be rewritten in the conservative form as:
\begin{equation}\label{example:4a}
u_t + \dfrac{\mu}{3}\, (u^3)_x + \epsilon\, u_{xxx} = 0.
\end{equation}
The initial condition of the single solitary wave is given by:
\begin{equation}\label{IC:4a}
u(x,0) = A\,\sech(k(x - x_0)), \quad x \in [0, 80], 
\end{equation}
where $ A = \sqrt{\dfrac{6c}{\mu}} $, $ k = \sqrt{\dfrac{c}{\epsilon}} $, $ c = 0.845 $, $ x_0 = 20 $, $ \mu = 3 $, and $ \epsilon = 1 $. Here, $ A $ represents the amplitude, and $ k $ denotes the width of the solitary wave. The exact solution is given by $u(x,t) = A\,\sech(k[(x - x_0) - ct])$. The comparison of error norms is presented in Table~\ref{Table:Eg4a}. It is evident from the table that the error norms obtained using the CCS are significantly lower than those obtained with the CNCS . Fig.~\ref{Fig:Eg4a} displays the solitary wave profiles at different time levels, showing that the soliton propagates to the right with nearly constant speed and amplitude over time, as expected.\par
\begin{table}
\captionsetup{position=above} 
\caption{Error and order of convergence for example~\ref{example:4} with the initial condition (\ref{IC:4a}) at time $t= 20$.}
\label{Table:Eg4a}
\centering
\begin{tabular*}{\textwidth}{@{\extracolsep{\fill}} l *{8}{c} }
\toprule
&\multicolumn{2}{c}{\textbf{SSPRK3-CNCS6}} & &
\multicolumn{2}{c}{\textbf{SSPRK3-CNCS8}} &
\multicolumn{2}{c}{\textbf{SSPRK3-CCS8}} \\
\cmidrule{2-3} \cmidrule{5-6} \cmidrule{7-8}
\textbf{N} & $\boldsymbol{L^{\infty}}$\textbf{-error} & \textbf{Rate}
&\textbf{N} & $\boldsymbol{L^{\infty}}$\textbf{-error} & \textbf{Rate}
& $\boldsymbol{L^{\infty}}$\textbf{-error} & \textbf{Rate}\\
\midrule
 200 & 1.3361e-03 & - & 100 & 1.0444e-02 & - & 1.3586e-03 & - & \\
 300 & 1.1040e-04 & 6.1496 & 150 & 1.4306e-03 & 4.9027 & 8.2887e-05 & 6.8975 & \\
 400 & 1.9275e-05 & 6.0666 & 200 & 1.5093e-04 & 7.8177 & 1.2108e-05 & 6.6865 & \\
 500 & 5.0076e-06 & 6.0403 & 250 & 2.4611e-05 & 8.1276 & 2.4582e-06 & 7.1455 & \\
 600 & 1.6684e-06 & 6.0283 & 300 & 5.5839e-06 & 8.1357 & 6.8040e-07 & 7.0452 & \\
 700 & 6.5933e-07 & 6.0228 & 350 & 1.6289e-06 & 7.9920 & 2.1899e-07 & 7.3541 & \\
 800 & 2.9544e-07 & 6.0118 & 400 & 5.5220e-07 & 8.1012 & 8.1742e-08 & 7.3801 & \\
\bottomrule
\end{tabular*}
\end{table}
\begin{figure}[htbp!]  
    \centering
    \begin{minipage}[b]{0.3\linewidth}
    \includegraphics[width=\linewidth, trim={0 0 0 0}, clip]{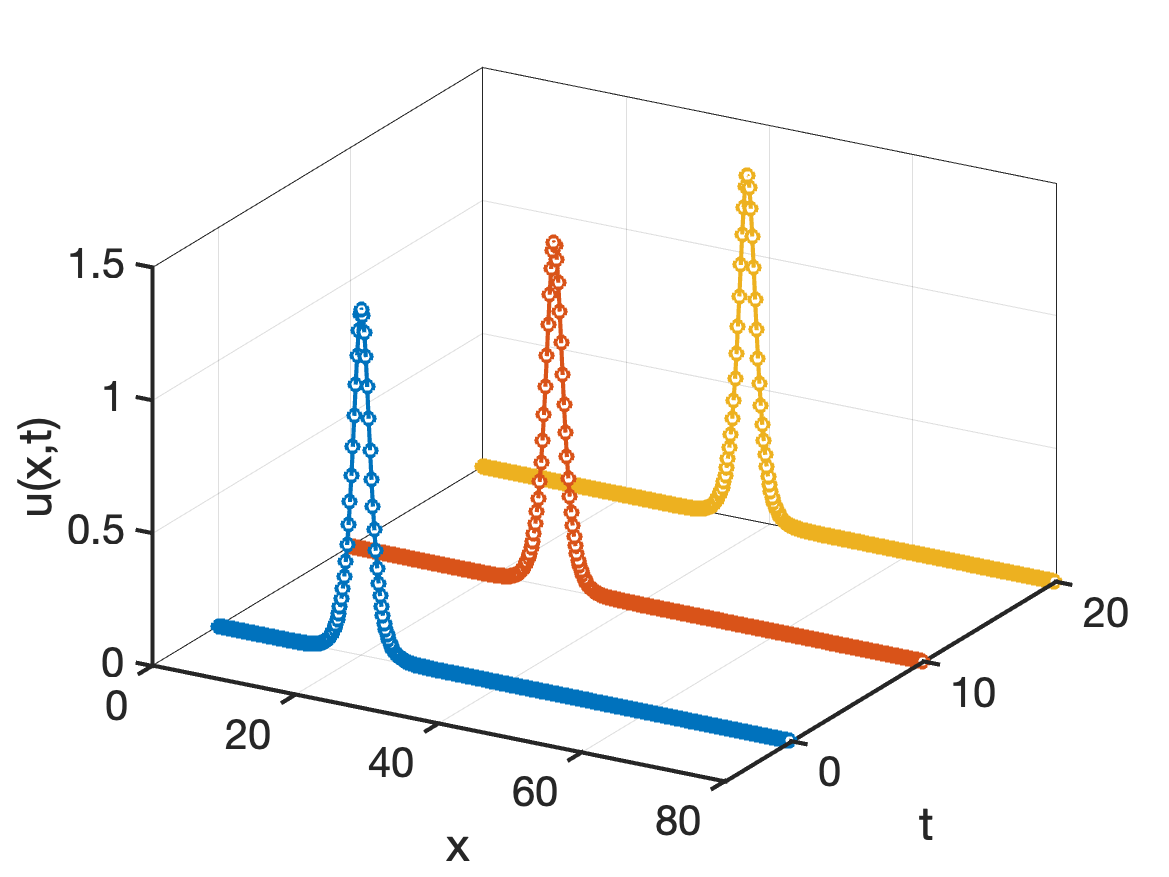}
      \subcaption{SSPRK3-CNCS6}
    \end{minipage}\hfill
    \begin{minipage}[b]{0.3\linewidth}
      \includegraphics[width=\linewidth, trim={0 0 0 0}, clip]{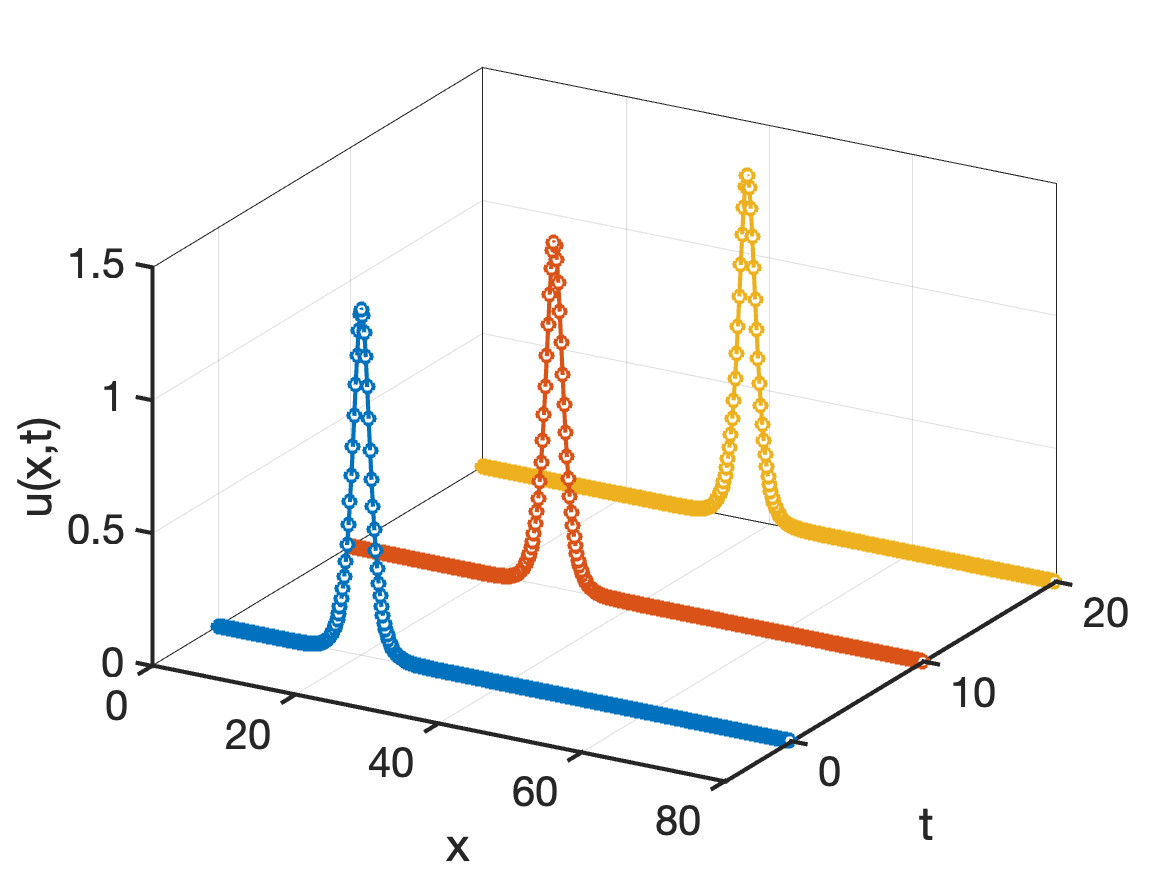}
      \subcaption{SSPRK3-CNCS8}
    \end{minipage}\hfill
    \begin{minipage}[b]{0.3\linewidth}
    \includegraphics[width=\linewidth, trim={0 0 0 0}, clip]{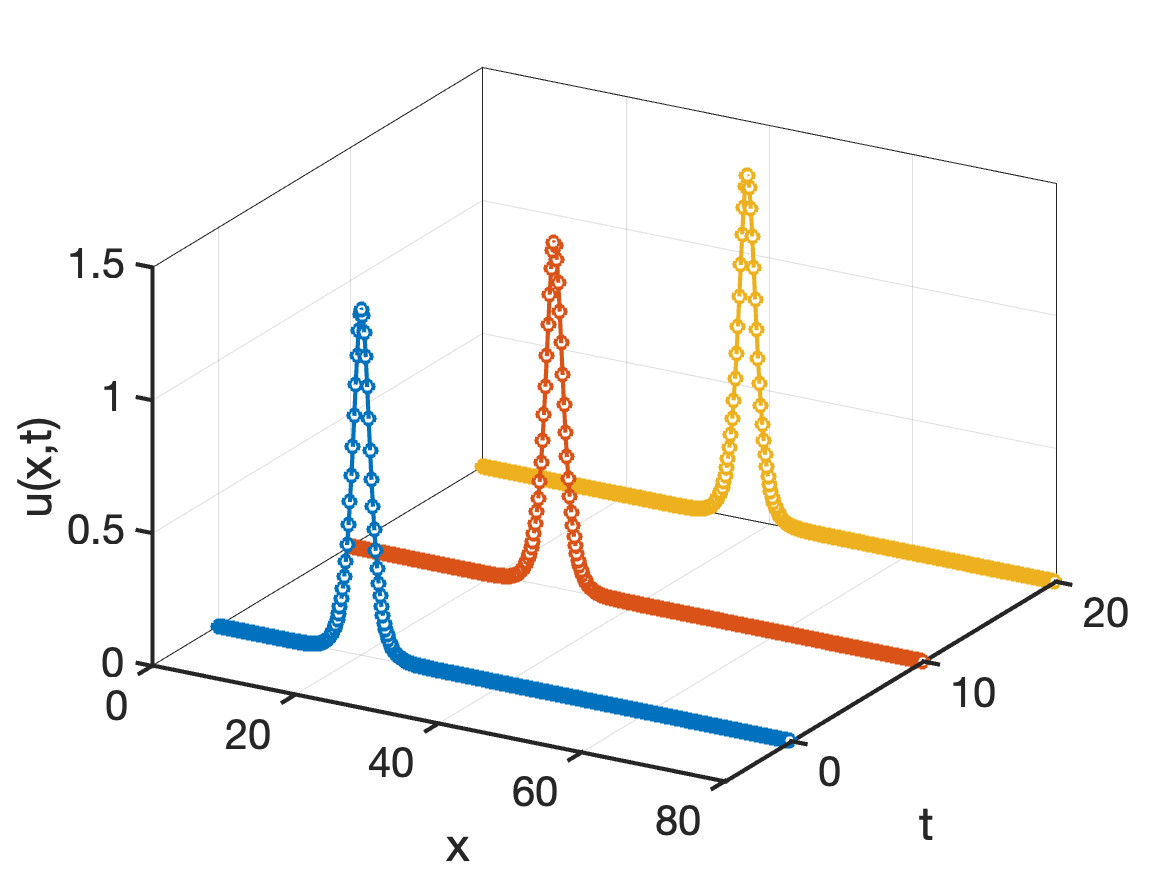}
      \subcaption{SSPRK3-CCS8}
    \end{minipage}\hfill
        \begin{minipage}[b]{0.3\linewidth}
    \includegraphics[width=\linewidth, trim={0 0 0 0}, clip]{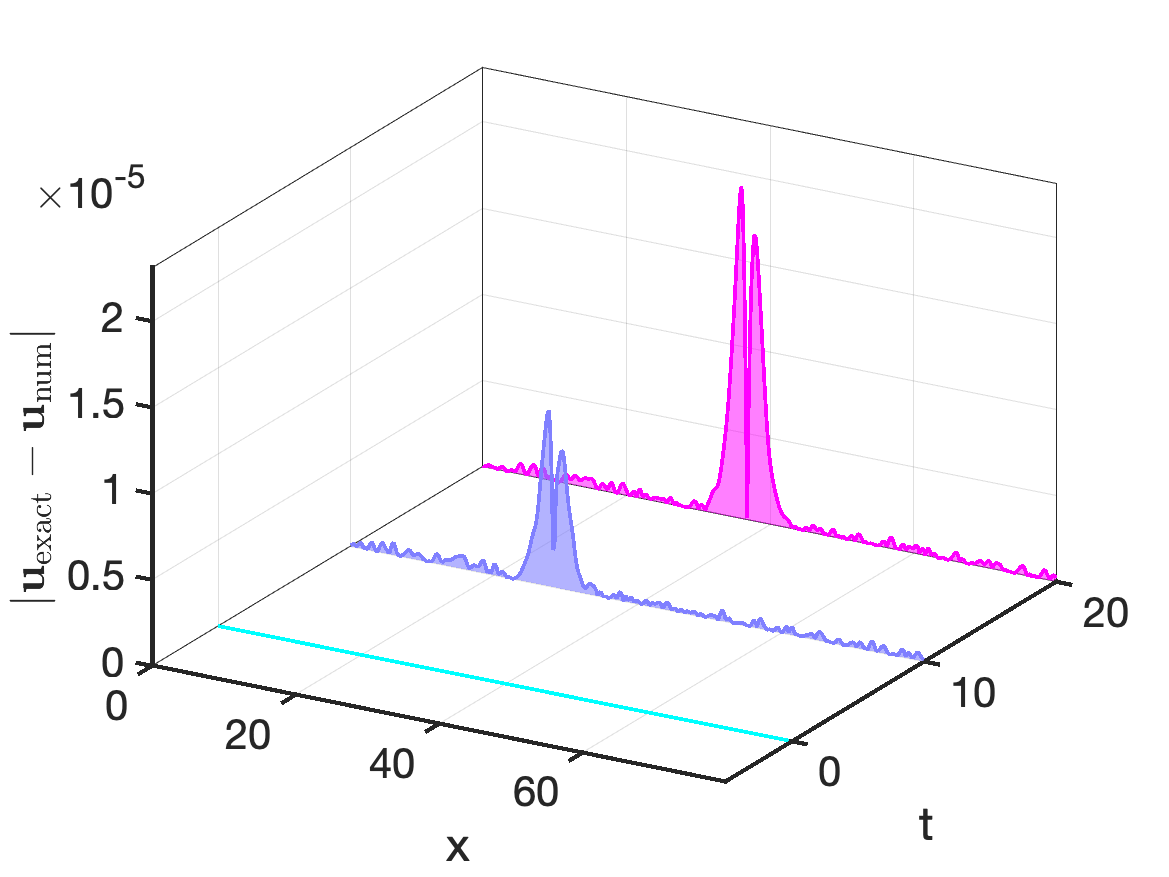}
      \subcaption{SSPRK3-CNCS6}
    \end{minipage}\hfill
    \begin{minipage}[b]{0.3\linewidth}
      \includegraphics[width=\linewidth, trim={0 0 0 0}, clip]{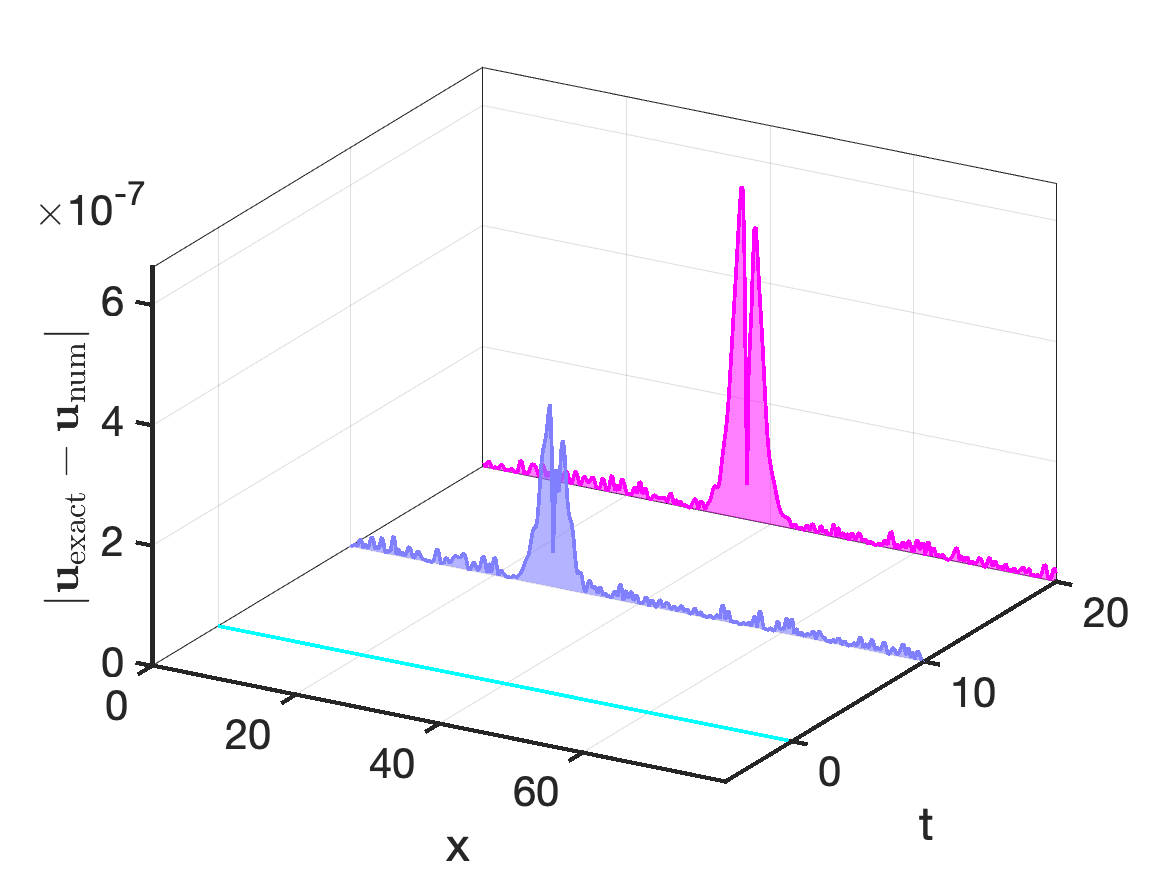}
      \subcaption{SSPRK3-CNCS8}
    \end{minipage}\hfill
    \begin{minipage}[b]{0.3\linewidth}
    \includegraphics[width=\linewidth, trim={0 0 0 0}, clip]{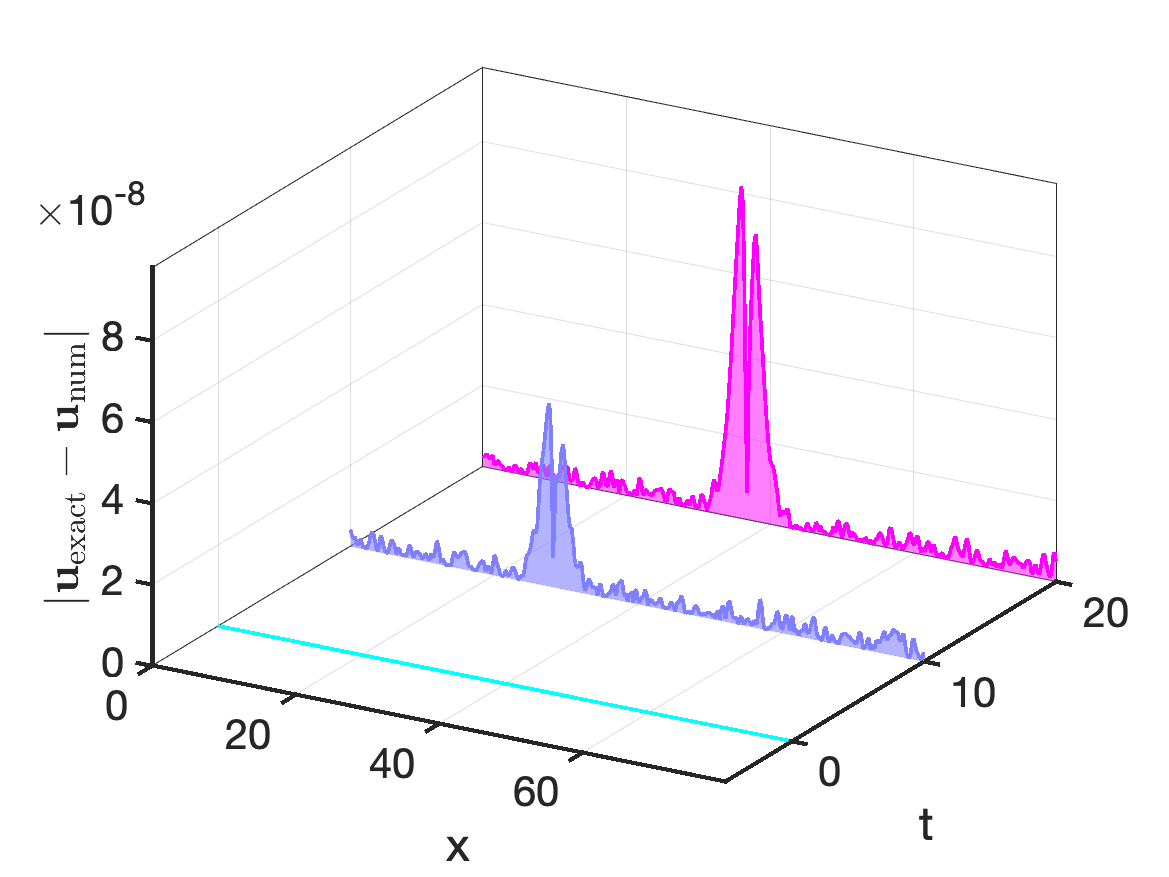}
      \subcaption{SSPRK3-CCS8}
    \end{minipage}\hfill
    \caption{Solutions and errors for the initial condition (\ref{IC:4a}) of example \ref{example:4a}. The first row shows exact (solid) and numerical ($\circ$) solutions at $t=0, 10, 20$, and the second row shows errors for CNCS6, CNCS8, and CCS8 for $N=400$.}
    \label{Fig:Eg4a}
\end{figure}
As a second case, we examine the interaction dynamics of two solitary waves with different amplitudes traveling in the same direction. The initial condition consists of two well-separated solitary waves and is defined as follows:
\begin{equation}\label{IC:4b}
u(x,0) = A_1\,\sech(k_1(x - x_1)) + A_2\,\sech(k_2(x - x_2)), \quad x \in [0, 80].
\end{equation}
To ensure interaction between the two solitary waves, the parameters are chosen as
$ A_i = \sqrt{\dfrac{6c_i}{\mu}} $, $ k_i = \sqrt{\dfrac{c_i}{\epsilon}}, \, i=1,2 $, with $ c_1 = 2 $, $ c_2 = 1 $, $ x_1 = 15 $, $ x_2 = 25 $, $ \mu = 3 $, and $ \epsilon = 1 $,
consistent with previous studies~\cite{biswas2011numerical}.
The evolution of the interaction is illustrated at various time levels in Fig.~\ref{Fig:Eg4b} for $N=500$. At time $ t = 0 $, the taller soliton with the amplitude $ 2.0 $ is positioned to the left of the smaller soliton, which has an initial amplitude of $ \sqrt2 $. Due to its larger amplitude, the taller soliton propagates faster and eventually catches up with the smaller one, leading to a collision around $ t = 6 $. As time progresses, the taller soliton overtakes and moves ahead of the smaller one.\par
By $ t = 16 $, the interaction is essentially complete, and both solitons regain their individual shapes and speeds, similar to their pre-interaction states. At $ t = 20 $, the absolute change in amplitude is observed to be $ 6.6779\times 10^{-3} $ for CNCS6, $ 6.6013\times 10^{-3} $ for CNCS8, and $ 6.5979\times 10^{-3} $ for CCS8,  for the larger soliton. The absolute change in amplitude is observed to be $ 1.6152\times 10^{-4} $ for CNCS6, $ 1.6219\times 10^{-4} $ for CNCS8, and $ 1.6224\times 10^{-4} $ for CCS8, for the smaller soliton, indicating the effectiveness of the numerical scheme in preserving the soliton characteristics.
\begin{figure}[htbp!]  
\centering
\begin{minipage}[b]{0.3\linewidth}
\includegraphics[width=\linewidth, trim={0 0 0 0}, clip]{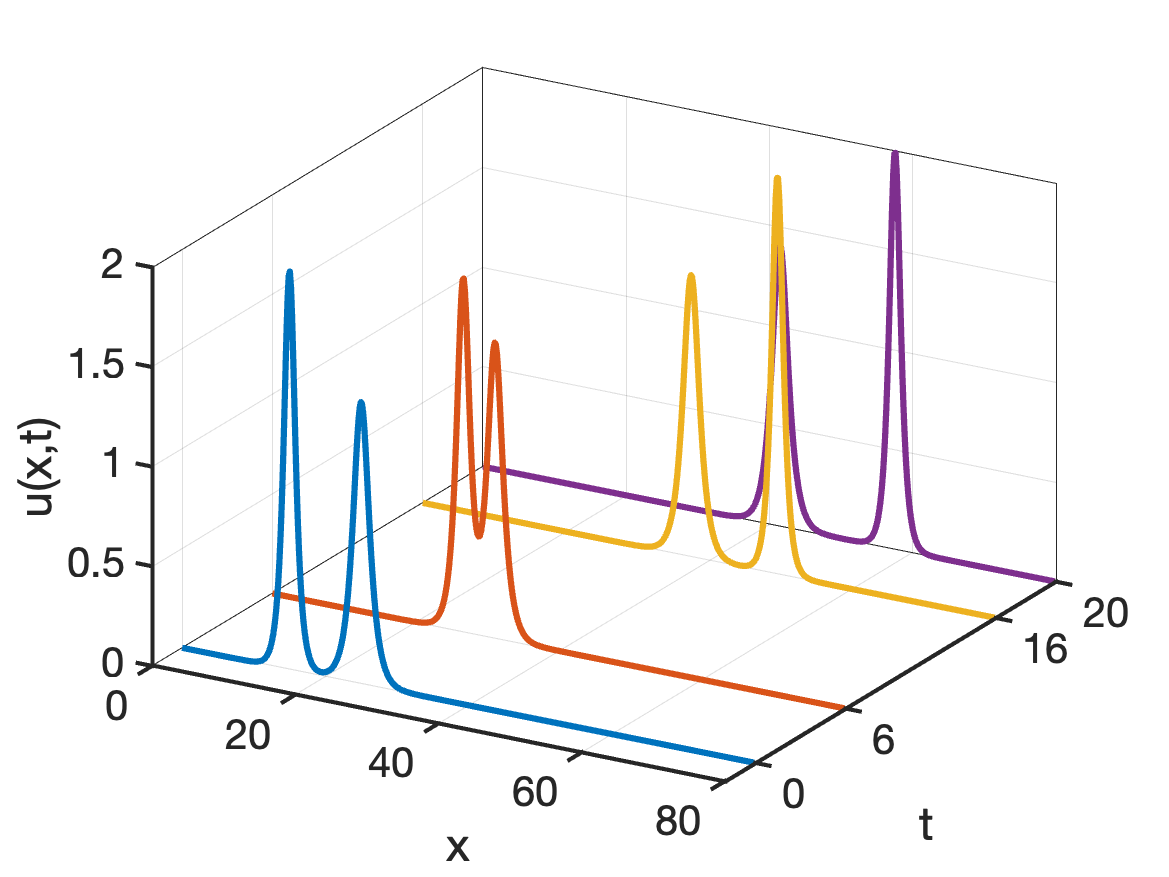}
  \subcaption{SSPRK3-CNCS6}
\end{minipage}\hfill
\begin{minipage}[b]{0.3\linewidth}
  \includegraphics[width=\linewidth, trim={0 0 0 0}, clip]{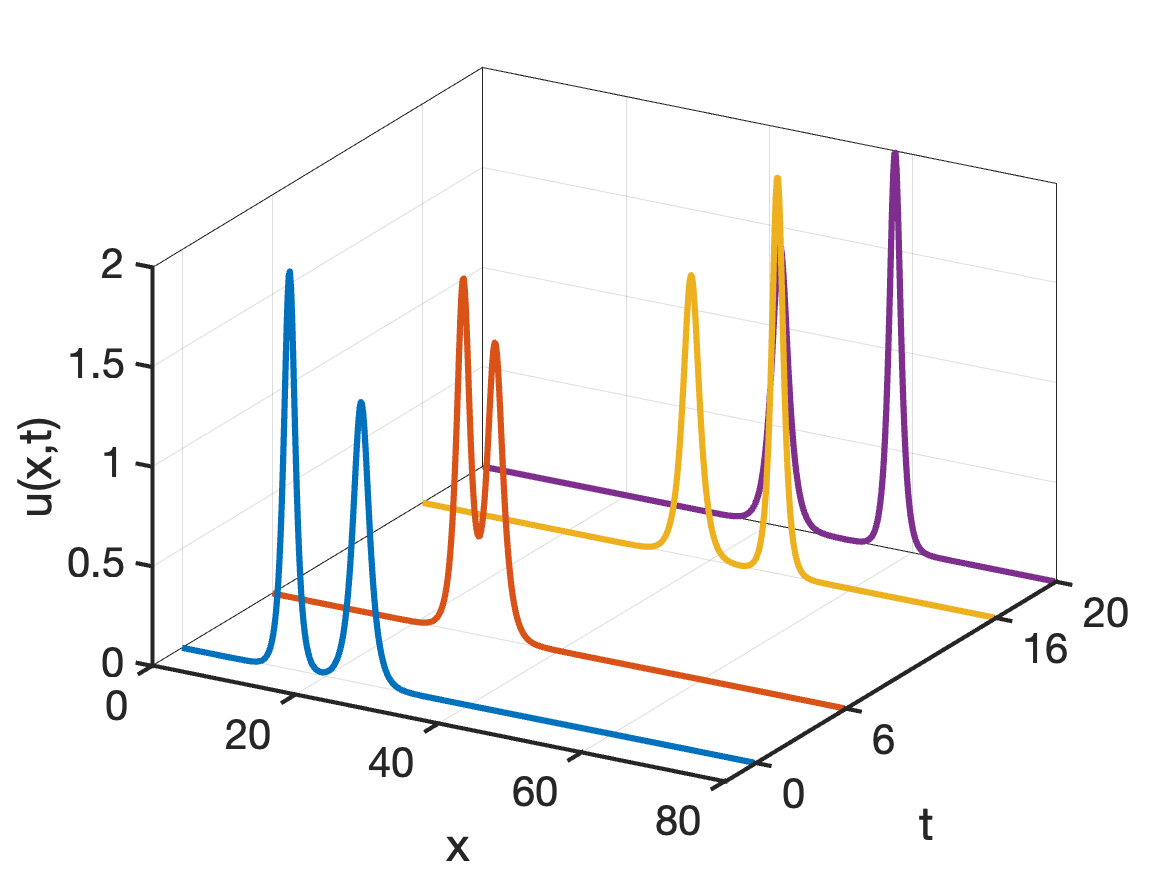}
  \subcaption{SSPRK3-CNCS8}
\end{minipage}\hfill
\begin{minipage}[b]{0.3\linewidth}
\includegraphics[width=\linewidth, trim={0 0 0 0}, clip]{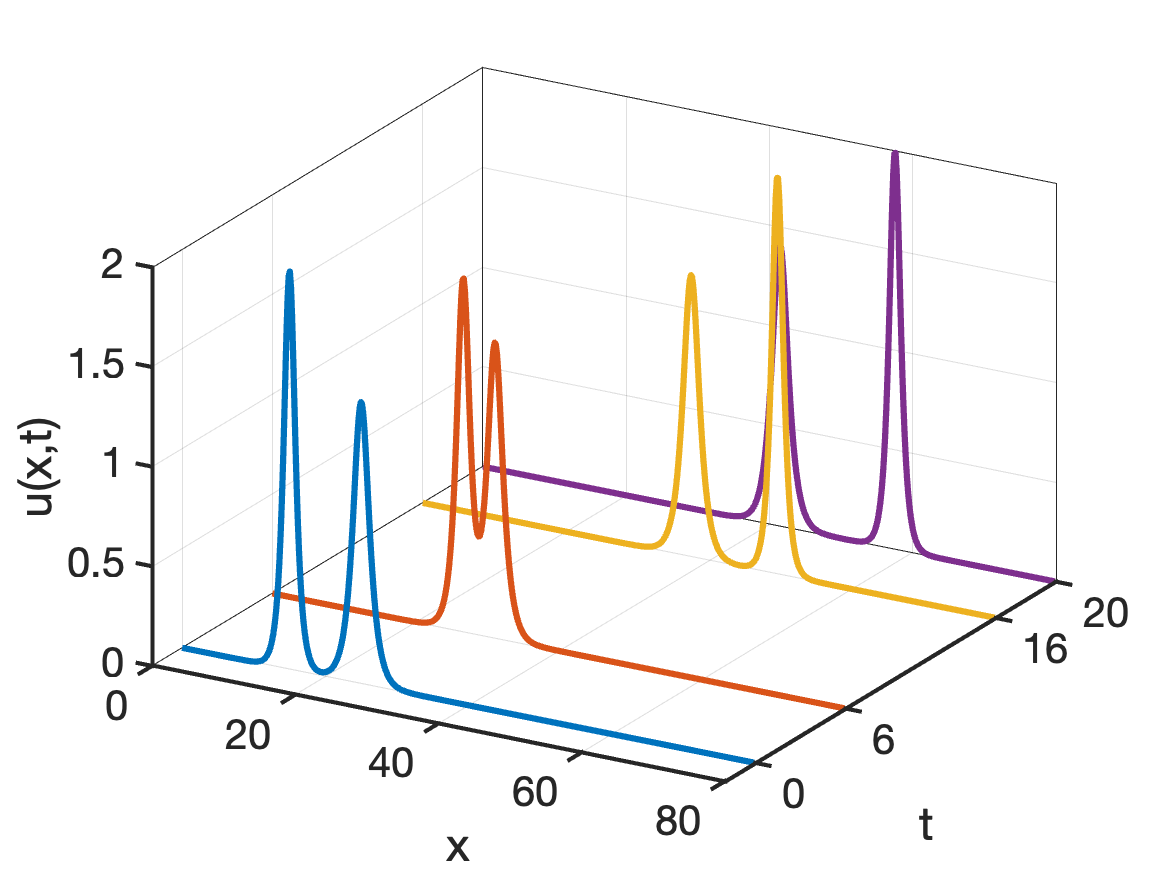}
  \subcaption{SSPRK3-CCS8}
\end{minipage}\hfill
\caption{Numerical solutions for initial condition (\ref{IC:4b}) of example \ref{example:4}.}
\label{Fig:Eg4b}
\end{figure}
\end{example}
\begin{example}\label{2D_example:1}
\normalfont
We examine the two-dimensional linear convection-dispersion equation
\begin{equation}
    u_t + 2 (u_x + u_y) + u_{xxx} + u_{yyy} = 0, \quad (x,y,t)\in [0, 2\pi] \times [0, 2\pi] \times [0, T],
\end{equation}
with the initial condition $ u(x,y,0) = \sin(x+y) $. The exact corresponding solution is \[ u (x, y, t) = \sin (x + y - 2t).
\]
Fig.~\ref{Fig:Eg2DA} captures the time evolution of the maximum amplitude of the solution up to $ t = 0.5 $, calculated using a grid resolution of $ N_x \times N_y = 100 \times 100 $. The stability characteristics of the CNCS for both sixth- and eighth-order accuracy levels are illustrated in Figs.~\ref{Fig:Eg2DA}(a) and~\ref{Fig:Eg2DA}(b), which correspond to dispersion numbers $ D_{\alpha} = 0.11 $ and the critical value $ D_{\alpha,\mathrm{cr}} = 0.12 $, respectively. Fig.~\ref{Fig:Eg2DA}(c) showcases the performance of the eighth-order CCS under analogous conditions, where $ D_{\alpha} = 0.011 $ and $ D_{\alpha,\mathrm{cr}} = 0.012 $. As anticipated, numerical instability arises beyond the critical threshold, in alignment with predictions from the GSA.\par
Furthermore, Table~\ref{Table:Eg2D} presents the computed $ L_{\infty}$-error and associated convergence rates in various schemes. All methods demonstrate high accuracy within the stable regime, with the CCS8 producing notably lower error levels, approximately an order of magnitude smaller than those from CNCS8. A visual comparison of the numerical and exact solutions, along with their error distribution at $ t = 0.5 $ for a coarser resolution of $ 40 \times 40 $, is provided in Fig.~\ref{Fig:Eg2DB}.
\begin{figure}[htbp!]  
    \centering
    \begin{minipage}[b]{0.3\linewidth}
    \includegraphics[width=\linewidth, trim={0 0 0 0}, clip]{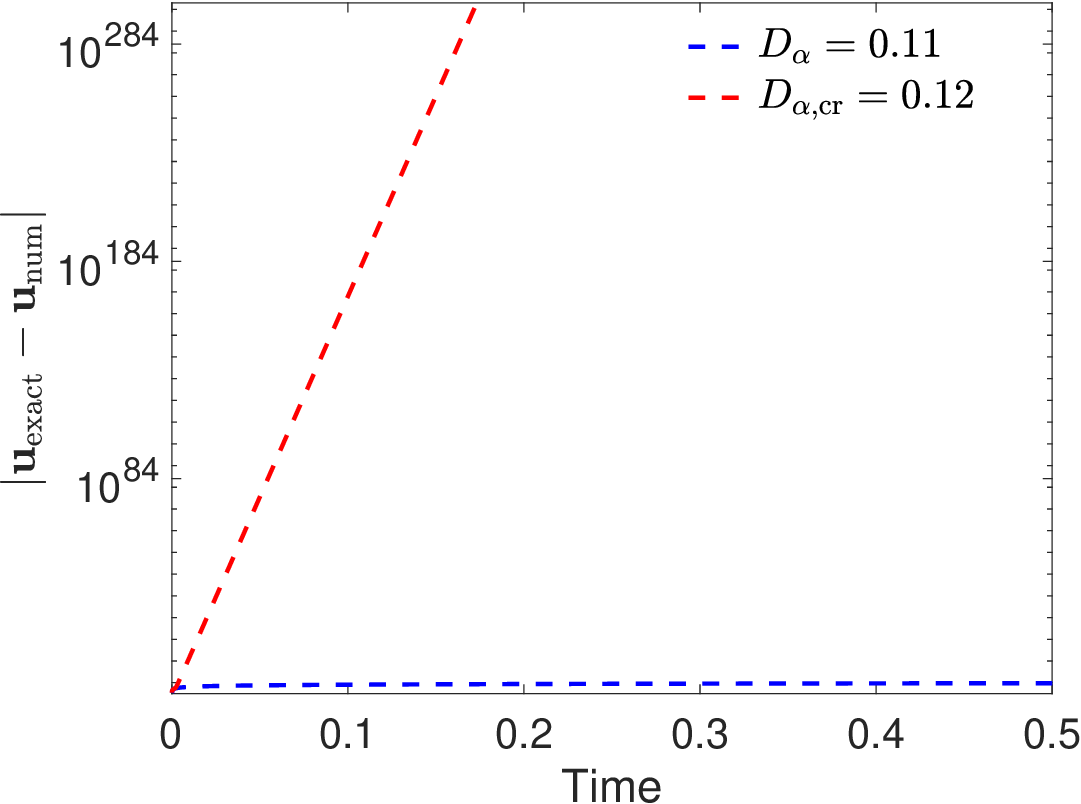}
      \subcaption{SSPRK3-CNCS6}
    \end{minipage}\hfill
    \begin{minipage}[b]{0.3\linewidth}
    \includegraphics[width=\linewidth, trim={0 0 0 0}, clip]{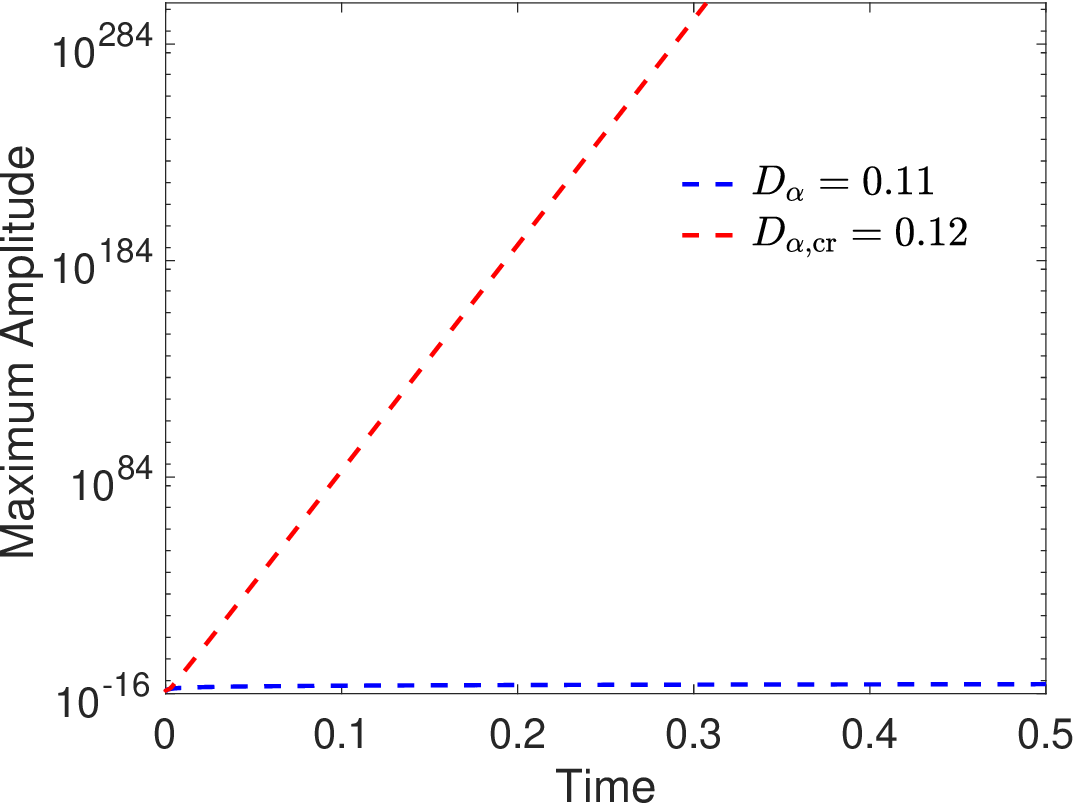}
      \subcaption{SSPRK3-CNCS8}
    \end{minipage}\hfill
    \begin{minipage}[b]{0.3\linewidth}
      \includegraphics[width=\linewidth, trim={0 0 0 0}, clip]{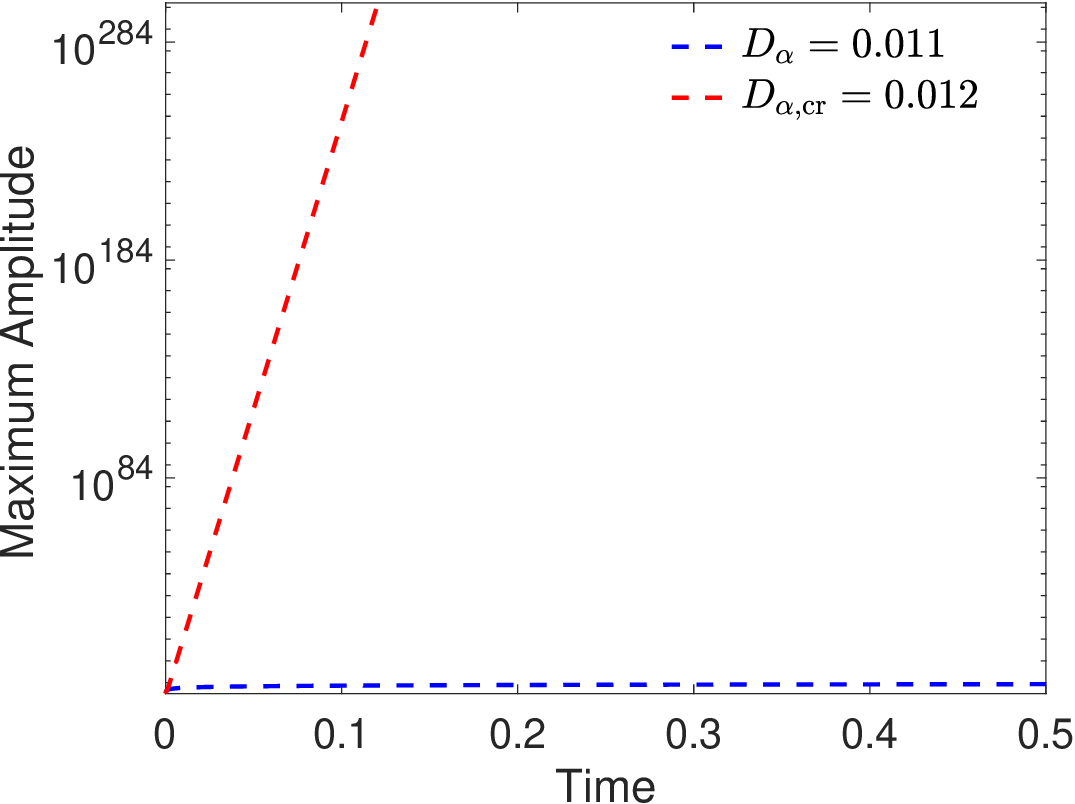}
      \subcaption{SSPRK3-CCS8}
    \end{minipage}\hfill    
    \caption{Variation of the maximum amplitude of $u$ over time for example~\ref{2D_example:1}.}
    \label{Fig:Eg2DA}
\end{figure}
\begin{table}
\captionsetup{position=above} 
\caption{Error and order of convergence for example~\ref{2D_example:1} at time $t= 0.5$.}
\label{Table:Eg2D}
\centering
\begin{tabular*}{\textwidth}{@{\extracolsep{\fill}} l *{8}{c} }
\toprule
&\multicolumn{2}{c}{\textbf{SSPRK3-CNCS6}} & &
\multicolumn{2}{c}{\textbf{SSPRK3-CNCS8}} &
\multicolumn{2}{c}{\textbf{SSPRK3-CCS8}} \\
\cmidrule{2-3} \cmidrule{5-6} \cmidrule{7-8}
$\boldsymbol{\textbf{N}_x \times \textbf{N}_y}$& $\boldsymbol{L^{\infty}}$\textbf{-error} & \textbf{Rate}
&$\boldsymbol{\textbf{N}_x \times \textbf{N}_y}$ & $\boldsymbol{L^{\infty}}$\textbf{-error} & \textbf{Rate}
& $\boldsymbol{L^{\infty}}$\textbf{-error} & \textbf{Rate}\\
\midrule
 10 $\times$ 10 & 6.2032e-05 & - & 10 $\times$ 10 & 2.2960e-06 & - & 2.1173e-07 & - & \\
 30 $\times$ 30 & 8.5096e-08 & 6.0000 & 15 $\times$ 15 & 9.2835e-08 & 7.9121 & 9.8751e-09 & 7.5599 & \\
 50 $\times$ 50 & 3.9644e-09 & 6.0029 & 20 $\times$ 20 & 9.2681e-09 & 8.0097 & 1.0497e-09 & 7.7917 & \\
 70 $\times$ 70 & 5.2590e-10 & 6.0035 & 25 $\times$ 25 & 1.5523e-09 & 8.0074 & 1.8189e-10 & 7.8552 & \\
 90 $\times$ 90 & 1.1674e-10 & 5.9891 & 30 $\times$ 30 & 3.5959e-10 & 8.0219 & 4.3087e-11 & 7.8991 & \\
 110 $\times$ 110 & 3.4132e-11 & 6.1280 & 35 $\times$ 35 & 1.0484e-10 & 7.9957 & 1.2512e-11 & 8.0217 & \\
\bottomrule
\end{tabular*}
\end{table}
\begin{figure}[htbp!]  
    \centering
    \begin{minipage}[b]{0.3\linewidth}
    \includegraphics[width=\linewidth, trim={0 0 0 0}, clip]{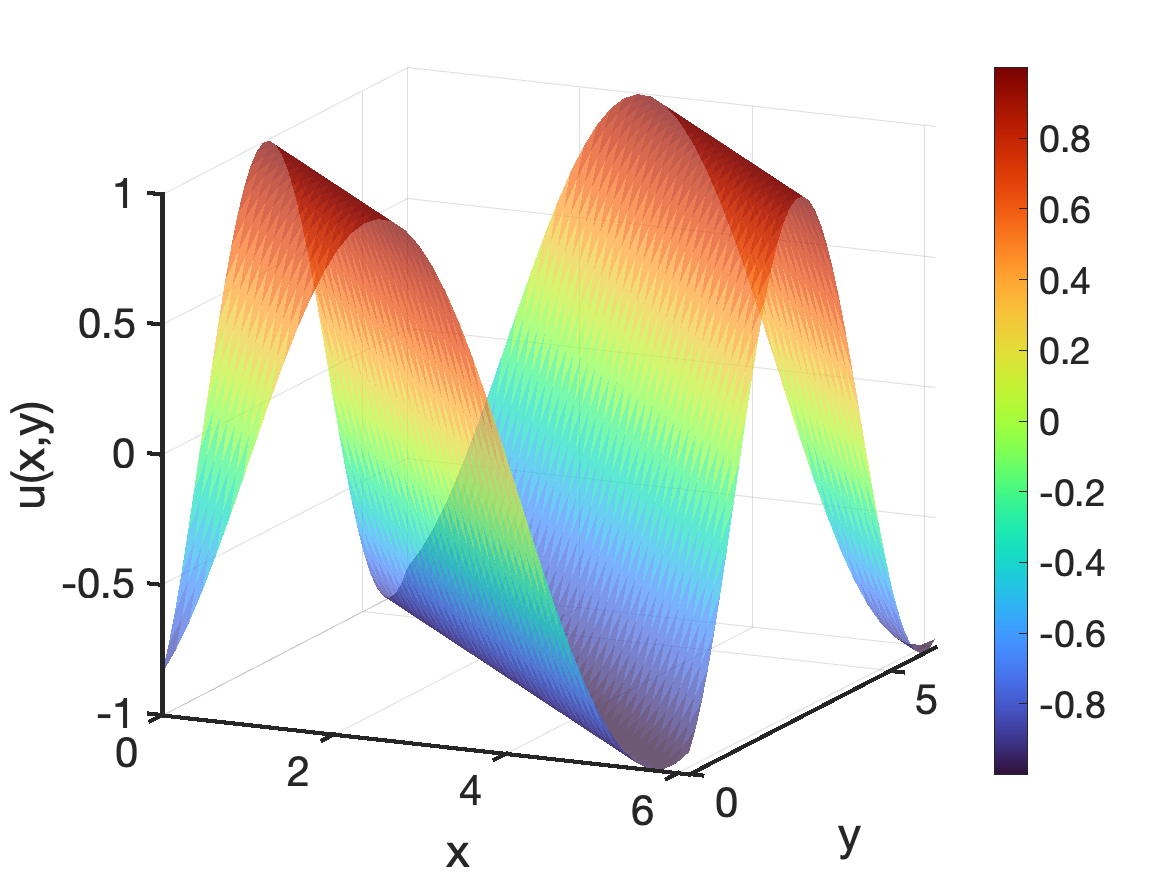}
      \subcaption{SSPRK3-CNCS6}
    \end{minipage}\hfill
    \begin{minipage}[b]{0.3\linewidth}
      \includegraphics[width=\linewidth, trim={0 0 0 0}, clip]{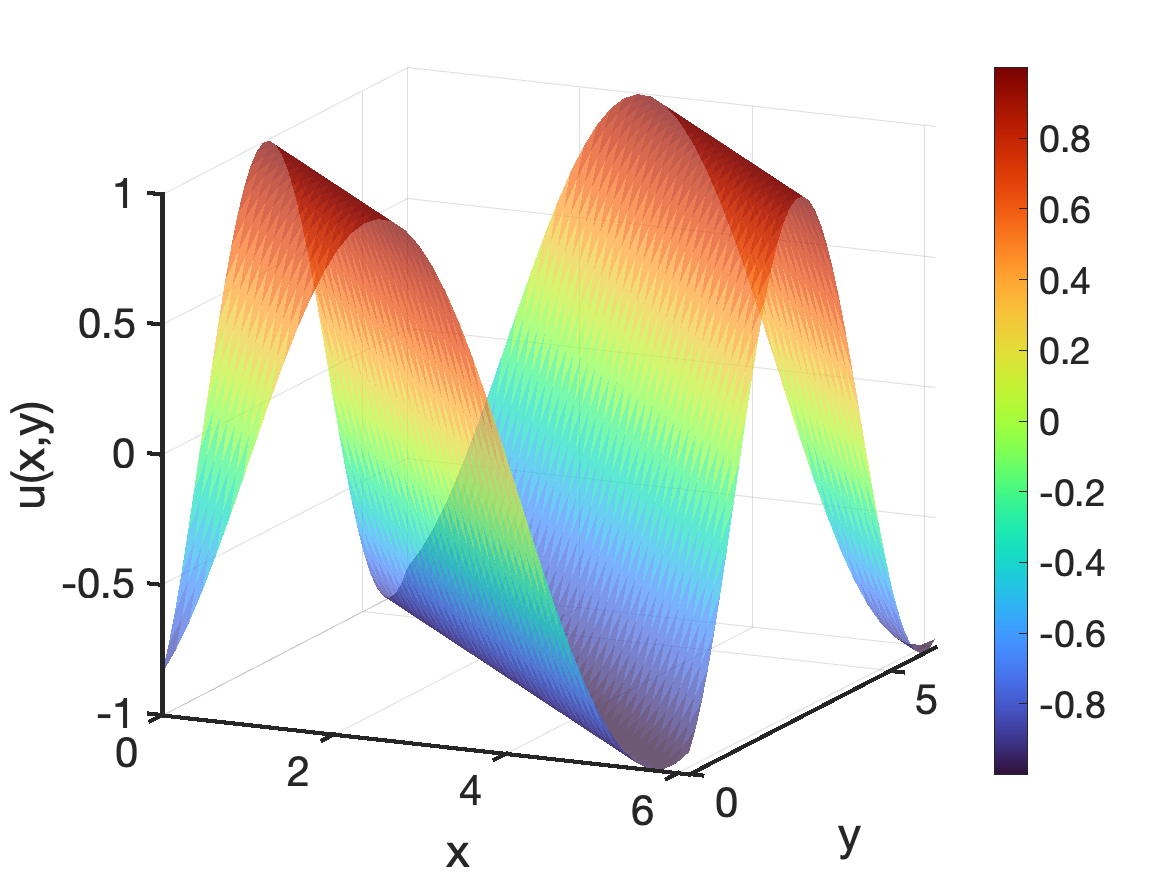}
      \subcaption{SSPRK3-CNCS8}
    \end{minipage}\hfill
    \begin{minipage}[b]{0.3\linewidth}
    \includegraphics[width=\linewidth, trim={0 0 0 0}, clip]{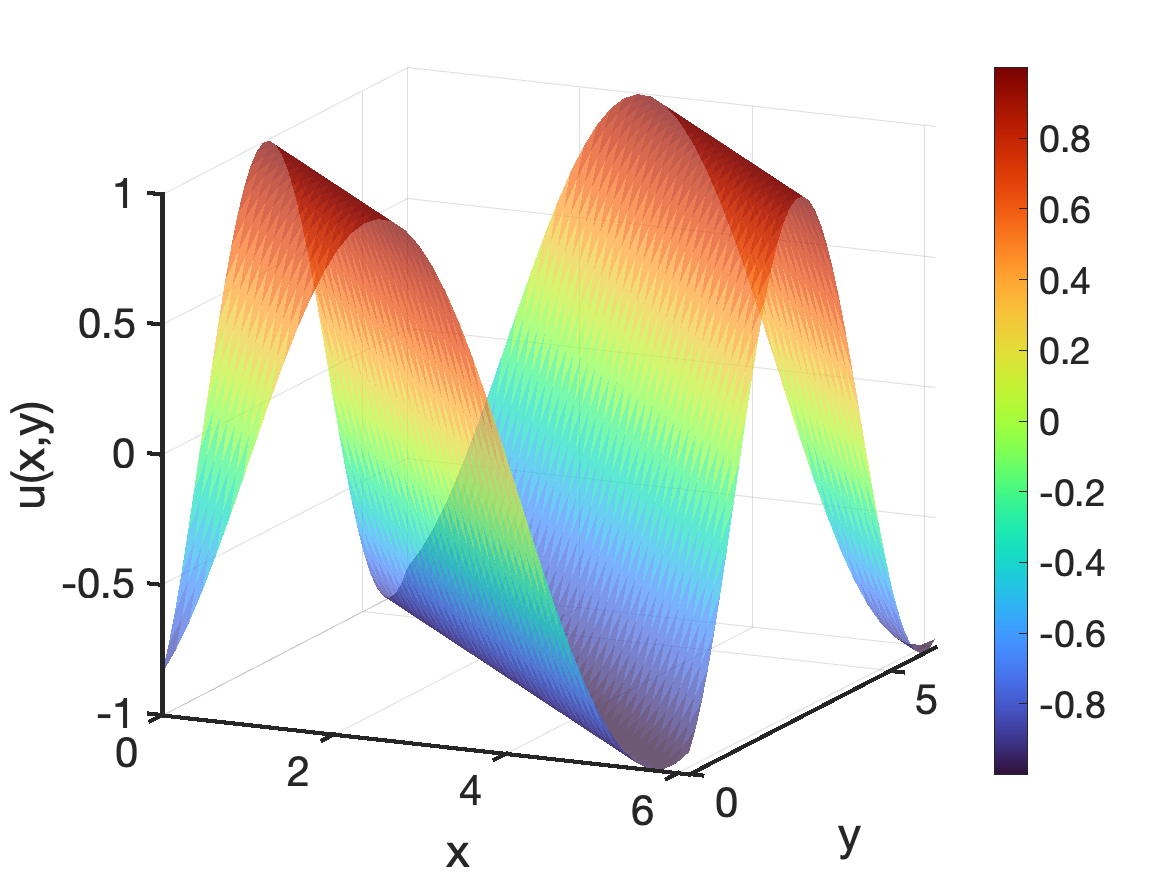}
      \subcaption{SSPRK3-CCS8}
    \end{minipage}\hfill
        \begin{minipage}[b]{0.3\linewidth}
    \includegraphics[width=\linewidth, trim={0 0 0 0}, clip]{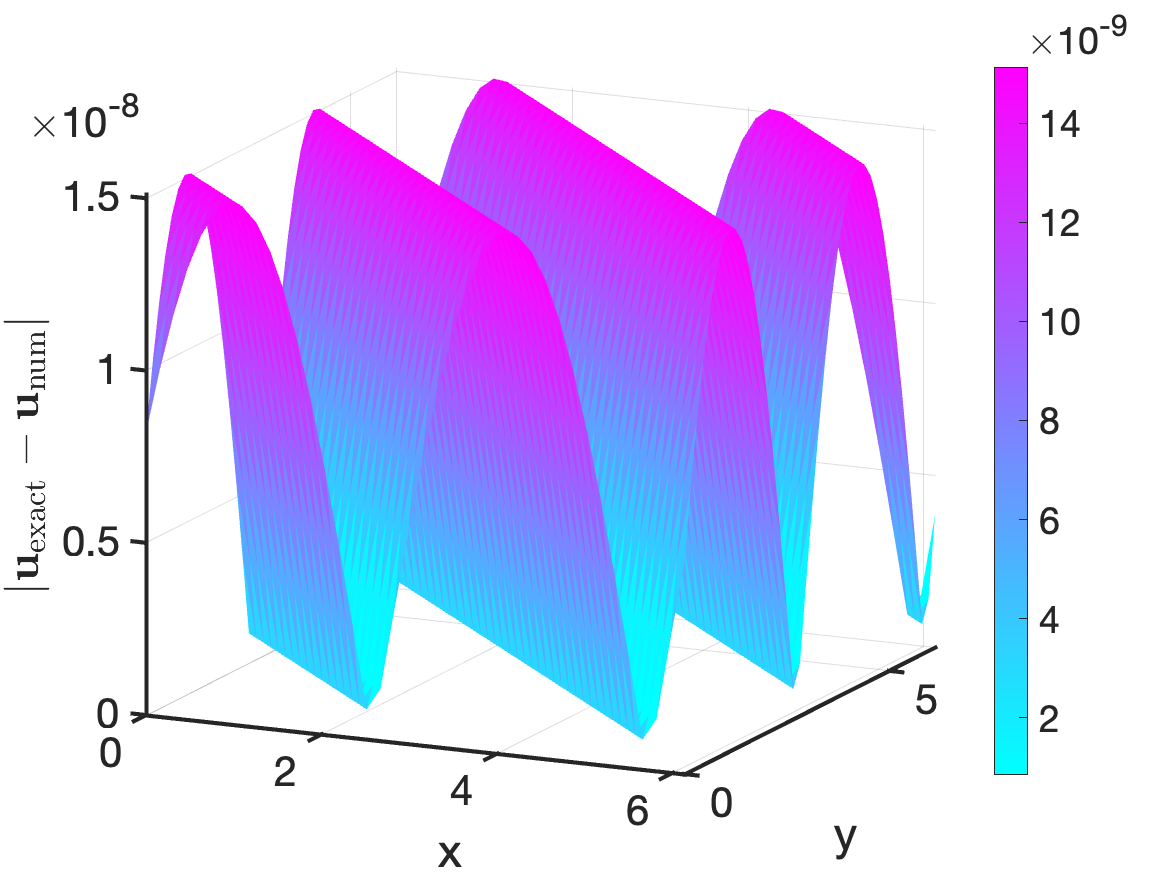}
      \subcaption{SSPRK3-CNCS6}
    \end{minipage}\hfill
    \begin{minipage}[b]{0.3\linewidth}
      \includegraphics[width=\linewidth, trim={0 0 0 0}, clip]{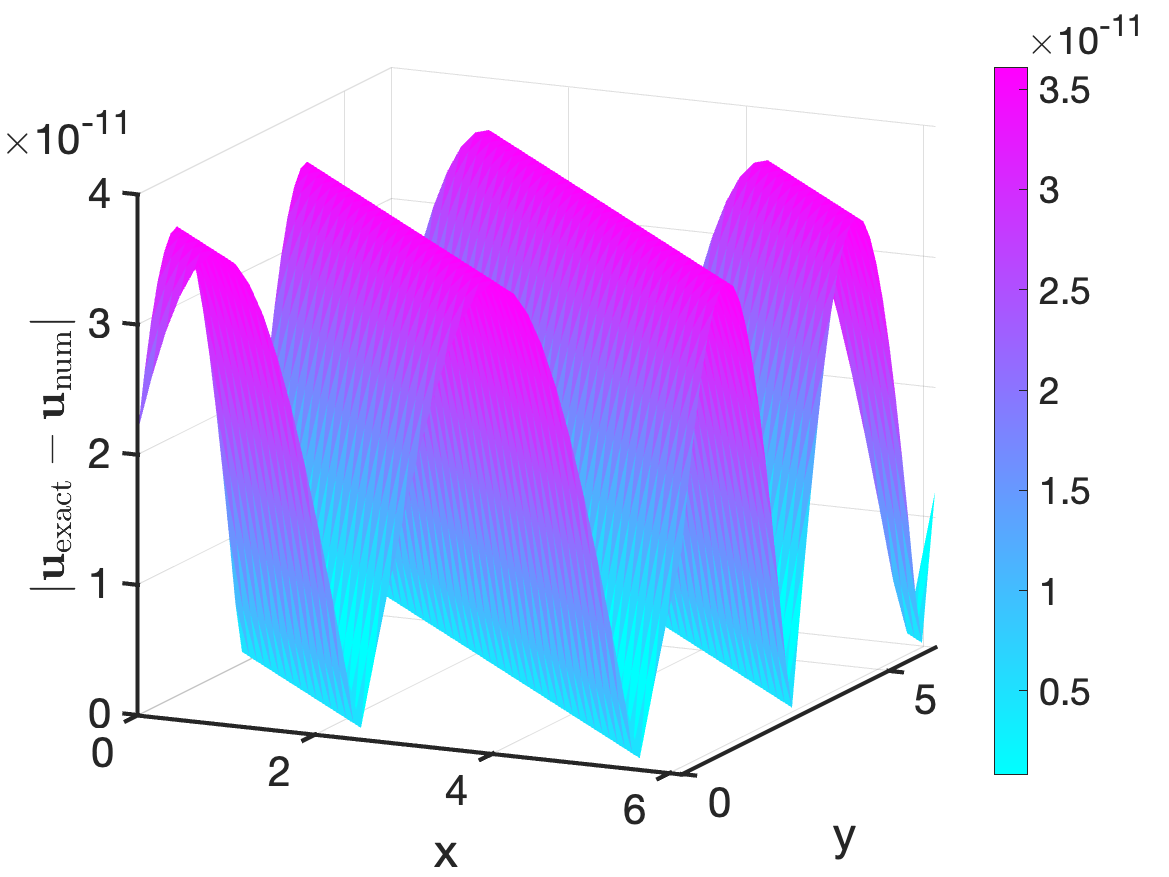}
      \subcaption{SSPRK3-CNCS8}
    \end{minipage}\hfill
    \begin{minipage}[b]{0.3\linewidth}
    \includegraphics[width=\linewidth, trim={0 0 0 0}, clip]{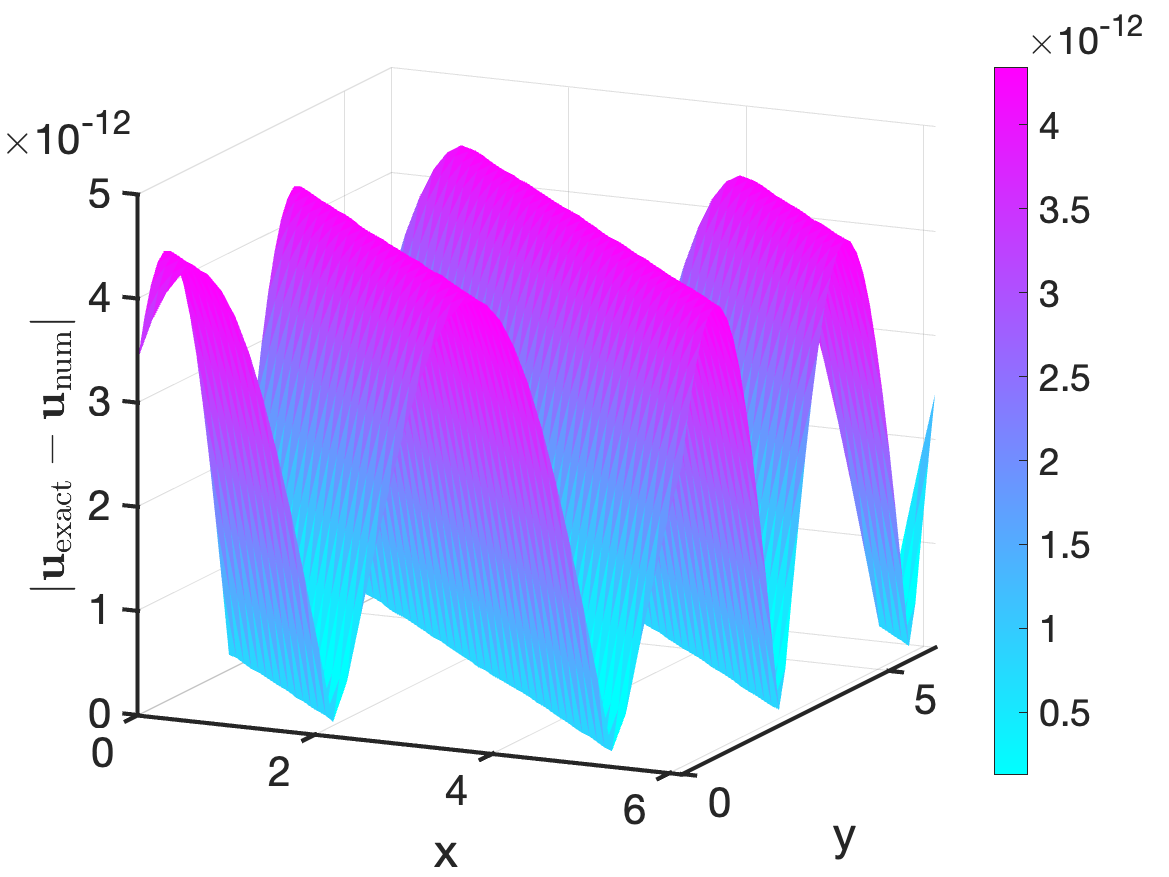}
      \subcaption{SSPRK3-CCS8}
    \end{minipage}\hfill
    \caption{Solutions and errors for example \ref{2D_example:1}. The first row shows exact (solid) and numerical ($\circ$) solutions at $t=0.5$, and the second row shows errors for CNCS6, CNCS8, and CCS8 for $N=40$.}
    \label{Fig:Eg2DB}
\end{figure}
\end{example}
\section{Conclusions}\label{Sec:Conclusion}
In this article, we have performed the GSA for the 1D and 2D convection-dispersion equations. For spatial discretization, we have used three compact schemes, CNCS6, CNCS8 and CCS8 with SSPRK3 for time discretization. The numerical amplification factor, normalized phase speed, and normalized group velocity were evaluated across the full spectral plane to assess stability and dispersive accuracy. We have observed that the node-centered schemes (CNCS6 and CNCS8) offered improved stability at moderate Courant and dispersion numbers, but exhibited spurious waves and dispersion errors at higher wave numbers. In contrast, the cell-centered scheme (CCS8) demonstrated enhanced resolution of high-frequency modes, leveraging its ability to extend the spectral domain up to \( kh = 2\pi \). However, tighter dispersion parameter tuning is still required to avoid instability. The theoretical predictions were validated through numerical experiments involving linear and nonlinear one- and two-dimensional test problems, showing strong agreement between the predicted stability thresholds and observed numerical behavior. Moreover, the test cases considered include the 1D and 2D linear convection-dispersion equations and the nonlinear KdV equation for simulating the motion of a single solitary wave and double soliton interactions with known analytical solutions. The study also covered problems with small third-derivative coefficients, focusing on dispersive shock formations under continuous and discontinuous initial conditions, where analytical solutions during the interaction are unknown. Additionally, we used the mKdV equation to examine the propagation of a single solitary wave with known analytical solutions and the interaction of two solitary waves in scenarios where analytical solutions are unavailable. This confirms the robustness of the spectral analysis approach in guiding scheme selection and parameter tuning. Overall, the results highlight the importance of comprehensive spectral characterization when designing high-order numerical methods for convection–dispersion systems.
\appendix
\section{Spatial discretization in the y-direction}\label{Appendix1}
Compact finite difference schemes for the $y$-direction follow the formulations in~\cite{lele1992compact,li2006high,liu2013new,salian2024novel}.

\subsection*{CNCS}
The CNCS formulations for first- and third-order derivatives are:
\begin{equation}\label{FDCNCS1}
\begin{split}   
    \alpha_1 (g_{2,y}^{\prime})_{i,j-1} + (g_{2,y}^{\prime})_{i,j} + \alpha_1 (g_{2,y}^{\prime})_{i,j+1} &= a_1\dfrac{ (g_{2,y})_{i,j+1} - (g_{2,y})_{i,j-1} }{2h_y} + b_1\dfrac{ (g_{2,y})_{i,j+2} - (g_{2,y})_{i,j-2} }{4h_y} \\
    &\quad + c_1\dfrac{ (g_{2,y})_{i,j+3} - (g_{2,y})_{i,j-3} }{6h_y},
\end{split}    
\end{equation}
\begin{equation}\label{TDCNCS1}
\begin{split}
    \alpha_1 (f_{2,y}^{\prime\prime\prime})_{i,j-1} + (f_{2,y}^{\prime\prime\prime})_{i,j} + \alpha_1 (f_{2,y}^{\prime\prime\prime})_{i,j+1} &= a_1\frac{ (f_{2,y})_{i,j+2} - 2(f_{2,y})_{i,j+1} + 2(f_{2,y})_{i,j-1} - (f_{2,y})_{i,j-2} }{2h_y^3} \\
    &\quad + b_1\frac{ (f_{2,y})_{i,j+3} - 3(f_{2,y})_{i,j+1} + 3(f_{2,y})_{i,j-1} - (f_{2,y})_{i,j-3} }{8h_y^3} \\
    &\quad + c_1\frac{ (f_{2,y})_{i,j+4} - 4(f_{2,y})_{i,j+1} + 4(f_{2,y})_{i,j-1} - (f_{2,y})_{i,j-4} }{20h_y^3}.
\end{split}
\end{equation}
\subsection*{CCS}
Similarly, for CCS, the discretizations are:
\begin{equation}\label{FDCCS1}
\begin{split}
    \alpha_1 (g_{2,y}^{\prime})_{i,j-1} + (g_{2,y}^{\prime})_{i,j} + \alpha_1 (g_{2,y}^{\prime})_{i,j+1} &= a_1\dfrac{(g_{2,y})_{i,j+\frac{1}{2}} - (g_{2,y})_{i,j-\frac{1}{2}}}{h_y} + b_1\dfrac{ (g_{2,y})_{i,j+1} - (g_{2,y})_{i,j-1} }{2h_y} \\
    &\quad + c_1\dfrac{ (g_{2,y})_{i,j+\frac{3}{2}} - (g_{2,y})_{i,j-\frac{3}{2}} }{3h_y},
\end{split}
\end{equation}
\begin{equation}\label{TDCCS1}
\begin{split}
    \alpha_1 (f_{2,y}^{\prime\prime\prime})_{i,j-1} + (f_{2,y}^{\prime\prime\prime})_{i,j} + \alpha_1 (f_{2,y}^{\prime\prime\prime})_{i,j+1} &= a_1\dfrac{ 4(f_{2,y})_{i,j+1} -8(f_{2,y})_{i,j+\frac{1}{2}} +8(f_{2,y})_{i,j-\frac{1}{2}} -4(f_{2,y})_{i,j-1} }{h_y^3} \\
    &\quad + b_1\dfrac{ 8(f_{2,y})_{i,j+\frac{3}{2}} -12(f_{2,y})_{i,j+1} +12(f_{2,y})_{i,j-1} -8(f_{2,y})_{i,j-\frac{3}{2}} }{5h_y^3} \\
    &\quad + c_1\dfrac{ 8(f_{2,y})_{i,j+\frac{5}{2}} -20(f_{2,y})_{i,j+1} +20(f_{2,y})_{i,j-1} -8(f_{2,y})_{i,j-\frac{5}{2}} }{35h_y^3}.
\end{split}
\end{equation}
\section{Equivalent Wavenumber Expressions}\label{Appendix2}

The expressions for the equivalent wavenumbers associated with the convection and dispersion operators are provided below for both CNCS and CCS discretization schemes.

\subsection*{CNCS}

In the CNCS, the equivalent wavenumber corresponding to the convection operator is given by~\cite{lele1992compact}:
\begin{equation}
    k^{[1]}_{\text{eq}} h = \frac{a_1 \sin(kh) + \dfrac{b_1}{2} \sin(2kh) + \dfrac{c_1}{3} \sin(3kh)}{1 + 2\alpha_1 \cos(kh)},
\end{equation}
while the expression for the dispersion operator is~\cite{lele1992compact, salian2024novel}:
\begin{equation}
    (k^{[2]}_{\text{eq}})^3 h^3 = \frac{a_1 \left[2 \sin(kh) - \sin(2kh)\right] + \dfrac{b_1}{4} \left[3 \sin(kh) - \sin(3kh)\right] + \dfrac{c_1}{10} \left[4 \sin(kh) - \sin(4kh)\right]}{1 + 2\alpha_1 \cos(kh)}.
\end{equation}

\subsection*{CCS}

For the CCS, the equivalent wavenumber for the convection operator is expressed as~\cite{liu2013new}:
\begin{equation}
    k^{[1]}_{\text{eq}} h = 2 \cdot \dfrac{a_1 \sin\left(\dfrac{kh}{2}\right) + \dfrac{b_1}{2} \sin(kh) + \dfrac{c_1}{3} \sin\left(\dfrac{3kh}{2}\right)}{1 + 2\alpha_1 \cos(kh)},
\end{equation}
and the equivalent wavenumber for the dispersion operator is given by~\cite{salian2024novel}:
\begin{equation}
(k^{[2]}_{\text{eq}})^3 h^3 = \frac{P + Q + R}{1 + 2\alpha_1 \cos(kh)},
\end{equation}
where,
\begin{align*}
P &= 2a_1 \left[8 \sin\left(\frac{kh}{2}\right) - 4 \sin(kh)\right], \\
Q &= \frac{2b_1}{5} \left[12 \sin(kh) - 8 \sin\left(\frac{3kh}{2}\right)\right], \\
R &= \frac{2c_1}{35} \left[20 \sin(kh) - 8 \sin\left(\frac{5kh}{2}\right)\right].
\end{align*}

\subsection*{\textbf{Data availability}}
 All data generated or analyzed during this study are included in this article.\\\\

\subsection*{\textbf{Conflict of interest} }
The authors declare that they have no conflict of interest.

\bibliography{references.bib}
\end{document}